\documentclass[10pt, a4paper]{article}
\usepackage[utf8]{inputenc}

\usepackage[T1]{fontenc}
\usepackage{amsmath, amssymb, amsthm}

\theoremstyle{remark}
\newtheorem*{remark}{Remark}
\numberwithin{equation}{subsection}
\usepackage{mathtools}
\usepackage{array}
\usepackage{hyperref}
\usepackage[english]{babel}
\usepackage{booktabs}
\usepackage{latexsym}
\usepackage{bm}
\usepackage[pdftex]{graphicx}
\usepackage{geometry}
\usepackage{breqn}
\usepackage{multicol}
\usepackage[table]{xcolor}
\usepackage{algorithm}
\usepackage{algpseudocode}

\pdfoutput=1
\usepackage{etoolbox}
\AtBeginEnvironment{tabular}{\linespread{0.8}\selectfont}

\title{Conjugate Power Series of Ultralogarithm Branches}
\author{Sergey V. Berezin}
\date{Version 1.1, \today}
\begin{document}

\maketitle

\begin{center}
\textbf{Affiliation}: Municipal Budgetary Institution, Iglino, Republic of Bashkortostan, Russia \\
\textbf{ORCID}: \href{https://orcid.org/0000-0001-8086-8288}{0000-0001-8086-8288} \\
\textbf{Email}: \href{mailto:bkcru@bk.ru}{bkcru@bk.ru} \\
\textbf{License}: Creative Commons Attribution 4.0 International (CC BY 4.0)
\end{center}

\begin{abstract}
We investigate multivalued solutions of the transcendental ultralogarithm equation $e^{aw} = 1 + a x e^{bw}$ for arbitrary real parameters $a > b > 0$. A sector-based conjugation criterion is formulated and rigorously proved, reducing the technical two-dimensional indexing of conjugate master-series expansions inside ($\vert{}x\vert{} < R$) and outside ($\vert{}x\vert{} > R$) the disk of convergence to a unified one-dimensional branch index $n \in \mathbb{Z}$. For integer exponents, a logarithmic analog of Vieta's theorem is established. Using multi-channel phase hodographs, we analyze branch evolution across integer, fractional, and transcendental incommensurable ($a = \pi, b = e$) exponents. It is demonstrated that logarithmic stratification eliminates branch cuts of fractional powers, ensuring global stability of iterative root refinement near branch points without auxiliary sheet-selection parameters. Finally, using the Sagdeev pseudopotential method in plasma physics, we confirm the continuous preservation of branch numbering in the Lambert $W$ limit and the physical invariance of the indices $n = 0$ (subsonic regime) and $n = -1$ (supersonic regime).
\end{abstract}


\tableofcontents
\section{Introduction}
\label{sec:intro}
This work is devoted to the study of multivalued solutions of the transcendental ultralogarithm equation
\begin{equation}
\label{eq:ulog_main}
e^{aw} = 1 + a\vert{}x\vert{} e^{i\arg(x)} e^{bw}, \quad a > b > 0,
\end{equation}
in the complex plane for real parameters $a$ and $b$ (including integer, fractional, and irrational values). Within the computational experiment, the evolution of the branches $w_n = \mathrm{ULog}_n(a;b;x)$ is examined along rays of constant argument $\arg(x) = \mathrm{const}$ as the modulus $\vert{}x\vert{}$ varies from $0$ to $x_{\max}$.

Inside the disk of convergence $\vert{}x\vert{} < R$, each branch with an integer index $n \in \mathbb{Z}$ is represented by a convergent master series (generalized Maclaurin expansion):
\begin{equation}
\label{eq:ulog_inner}
\mathrm{ULog}_n(a;b;x) = f^{0}_{n} + M\bigl(0;\, a;\, b;\, x e^{b f^{0}_{n}}\bigr), \quad \vert{}x\vert{} < R,
\end{equation}
where the master series operator is defined by the power expansion
\begin{equation}
\label{eq:master_series}
M(m;a;b;X) = m + X + \sum_{\ell=2}^{\infty} \left(\prod_{\gamma=1}^{\ell-1} (m - a\gamma + b\ell)\right) \frac{X^\ell}{\ell!},
\end{equation}
whose coefficients are \textbf{master numbers}---generalized factorial products with step $a$ and linear shift $b\ell$ introduced in~\cite{Berezin2026_arXiv} (see also standard notation for factorial powers in~\cite{Knuth1997}).

We refer to $f^{0}_{n}$ (from the Latin \textit{ferre}, ``to bear, to support'') as the \textbf{zero stationary reference value} (or \textbf{zero stationary anchor}) of the $n$-th branch: its real part vanishes identically ($\mathrm{Re}\,f^{0}_{n} = 0$) and is independent of $x$. In the ultralogarithm plane it sets the starting point on the imaginary axis at $x = 0$, in the ultraradical plane $y = e^w$ it determines the initial phase angle $\mathrm{Im}(f^{0}_{n}) = \frac{2\pi n}{a}$ on the unit circle, and it simultaneously fixes the integer Riemann sheet index $u = \bigl\lceil \frac{\mathrm{Im}(f^{0}_{n})}{2\pi} - \frac{1}{2} \bigr\rceil$ (where $\lceil \cdot \rceil$ denotes the standard ceiling function~\cite{Graham1994}).

Inside the disk of convergence, this reference value has the form
\begin{equation}
\label{eq:f1n}
f^{0}_{n} = \frac{2\pi i n}{a}, \quad n \in \mathbb{Z},
\end{equation}
and the radius of convergence $R$ is given by the formula:
\begin{equation}
\label{eq:radius_R}
R = \frac{\vert{}1 - a/b\vert{}^{b/a}}{\vert{}b - a\vert{}}.
\end{equation}
At the origin ($x = 0$), all branches satisfy the initial conditions $\mathrm{Re}(w_n(0)) = 0$ and $\mathrm{Im}(w_n(0)) = \frac{2\pi n}{a}$.

Outside the disk of convergence ($\vert{}x\vert{} > R$), the expansion~\eqref{eq:ulog_inner} diverges; however, the values of the same branches of the ultralogarithm can be computed using master series with transformed parameters and shifted initial points of two types:
\begin{align}
\label{eq:ulog_outer_h}
\mathrm{ULog}_n(a;b;x) &= f^{+}_{n} + M\left(0;\, b-a;\, -a;\, \frac{e^{-a f^{+}_{n}}}{a-b}\right), & f^{+}_{n} &= \frac{\ln\vert{}1/(ax)\vert{} + \bigl(\arg(1/(ax)) + 2\pi n_{+}\bigr)i}{b-a}, \\
\label{eq:ulog_outer_k}
\mathrm{ULog}_n(a;b;x) &= f^{-}_{n} + M\left(0;\, -b;\, a-b;\, \frac{e^{(a-b) f^{-}_{n}}}{bax}\right), & f^{-}_{n} &= \frac{\ln\vert{}ax\vert{} + \bigl(\arg(-ax) + 2\pi n_{-}\bigr)i}{-b},
\end{align}
where $n_{\pm} \in \mathbb{Z}$ are auxiliary integer indices of the outer expansions (within the anchor notation $f^{\pm}_{n}$ and phase indicators $L^{\pm}_{n}$, the subscript sign on $n$ is omitted for brevity, as the index type is uniquely specified by the superscript $+$ or $-$). Unlike the zero stationary anchor $f^{0}_{n}$, the quantities $f^{\pm}_{n}(x)$ act as \textbf{$x$-dependent moving outer anchors}. This explicit dependence on $x$ plays a dual role: on the one hand, as $\arg(x)$ varies, the imaginary parts of the anchors drift along the imaginary axis, so the correspondence between $n$ and $n_{\pm}$ is not permanently fixed and requires selecting the conjugate series; on the other hand, because all information about $x$ is already encoded in $f^{\pm}_{n}(x)$, we can determine \emph{a priori}~--- without computing the power series for all branches~--- which specific index $n_{+}$ or $n_{-}$ lands inside the sector of a target branch $n$. Since for $a > b > 0$ the radius of convergence~\eqref{eq:radius_R} strictly satisfies $R > 1/a$, the inequality $\ln\vert{}ax\vert{} > 0$ holds throughout the outer region $\vert{}x\vert{} > R$, ensuring that the real parts of the anchors never change sign: $\mathrm{Re}\,f^{+}_{n}(x) > 0$ and $\mathrm{Re}\,f^{-}_{n}(x) < 0$. Accordingly, we call $f^{+}_{n}(x)$ the \textbf{positive moving anchor} (corresponding to the exterior $\vert{}y\vert{} > 1$ of the unit disk in the ultraradical plane $y = e^w$) and $f^{-}_{n}(x)$ the \textbf{negative moving anchor} (corresponding to the interior $\vert{}y\vert{} < 1$). Note that near the boundary $\vert{}x\vert{} = R$, where the modulus of the transformed argument $X$ of the outer master series $M$ attains its maximum, the series correction $M$ is largest and can temporarily give $\mathrm{Re}\,w_n = \mathrm{Re}\,f^{\pm}_{n} + \mathrm{Re}\,M$ the sign opposite to that of its anchor; as $\vert{}x\vert{}$ grows, however, $X \to 0$, $M \to 0$, and the branch trajectory inevitably follows its anchor.

Indeed, setting $t = \frac{a-b}{b} > 0$ yields $aR = \frac{t+1}{t}\,t^{1/(t+1)}$; logarithmic differentiation shows that this function increases monotonically on $(0, e)$ and decreases on $(e, +\infty)$ with limits $\lim_{t \to 0^+} aR = \lim_{t \to +\infty} aR = 1$ and a unique maximum at $t = e$, establishing the sharp two-sided bound $1 < aR \le \frac{e+1}{e}\,e^{1/(e+1)} \approx 1.791$. Consequently, throughout the outer region $\vert{}x\vert{} > R$ one has $\vert{}ax\vert{} > aR > 1$, so $\ln\vert{}ax\vert{} > 0$ and the real parts of the anchors never change sign: $\mathrm{Re}\,f^{+}_{n}(x) > 0$ and $\mathrm{Re}\,f^{-}_{n}(x) < 0$.

The two-dimensional indexing $(0\pm, n)$ is purely technical in nature and is required for the proper substitution of parameters into the master series. For practical computations and Riemann surface analysis, a unified one-dimensional indexing of branches by an integer $n \in \mathbb{Z}$ is needed. This gives rise to the fundamental problem of uniquely determining which specific integer $n_{+}$ or $n_{-}$ corresponds to a given branch index $n$ outside the disk of convergence.

The geometric nature of this problem is transparently revealed in the complex $w$-plane. Each branch $w_n$ of the ultralogarithm is confined to its own horizontal phase strip (perpendicular to the imaginary axis), which we call the \textbf{sector} of that branch. Along a ray of constant argument $\arg(x) = \mathrm{const}$ with increasing modulus $\vert{}x\vert{}$, the evolution of the branch traces a continuous curve within this strip starting from the imaginary axis $\mathrm{Re}(w) = 0$. Trajectories belonging to different sectors never intersect or enter neighboring strips; at critical values of $\arg(x)$, they may only touch at isolated branch points on the shared boundaries of adjacent sectors.

Beyond the disk of convergence ($\vert{}x\vert{} > R$), the trajectory originating at $\mathrm{Re}(w) = 0$ admits a unique analytic continuation given by one of the outer master series~\eqref{eq:ulog_outer_h} or~\eqref{eq:ulog_outer_k}. The transformed power series generating this continuation is \emph{conjugate} to the initial series~\eqref{eq:ulog_inner}. Crucially, the parameters of this conjugate series can be determined \emph{a priori}~--- before evaluating $w_n$ itself~--- solely from the location of the reference value $f^{\pm}_{n}$ within the phase sector.

To construct a rigorous theory for reducing the two-dimensional index to a one-dimensional one, we introduce three interrelated concepts:
\begin{itemize}
\item \textbf{Conjugate series}~--- the power expansions~\eqref{eq:ulog_inner} and~\eqref{eq:ulog_outer_h}--\eqref{eq:ulog_outer_k} operating in tandem to represent the same analytic branch $w_n$ on opposite sides of the convergence boundary $\vert{}x\vert{} = R$ for a fixed $\arg(x)$ (the term \emph{conjugate} is used here in the sense of structural pairing, analogous to conjugate points in the calculus of variations, rather than conjugate trigonometric series in harmonic analysis);
\item \textbf{Sector-based conjugation criterion}~--- the rule for a priori identifying the expansion type ($\pm$) and the integer index $n_{+}$ or $n_{-}$ from the inclusion of the \emph{scaled phase coordinate} (reduced imaginary part) of the reference value $f^{\pm}_{n}$ inside the open scaled sector of the $n$-th branch;
\item \textbf{Boundary conjugation convention}~--- the rule for uniquely selecting the conjugate series in the degenerate case when the ray direction $\arg(x)$ places the reduced phase point directly on the shared boundary between two adjacent sectors.
\end{itemize}

\noindent\textbf{Theorem 1 (Gemini Sector-Based Conjugation Criterion).}
\label{th:sector_join}
Let $a > b > 0$ be real numbers. To align the periods of fractional and transformed exponents, multiply the imaginary parts of all reference values $f^{0\pm}_{n}$ ($n_{0\pm} \in \mathbb{Z}$) by the exponent $b$, assigning to each anchor a \textbf{scaled phase coordinate} $L^{0\pm}_{n} = \mathrm{Im}(b f^{0\pm}_{n})$ on the imaginary axis. Partition the imaginary axis into the corresponding \textbf{scaled phase sectors} $s_n = \left[L_{\min}^{(n)}, L_{\max}^{(n)}\right]$ of width $\frac{2\pi b}{a}$ (representing the original phase strips of width $\frac{2\pi}{a}$ in the $w$-plane) centered at the points $L^{0}_{n} = \frac{2\pi n b}{a}$:
\begin{equation}
\label{eq:sector_bounds}
L_{\min}^{(n)} = L^{0}_{n} - \frac{\pi b}{a}, \qquad L_{\max}^{(n)} = L^{0}_{n} + \frac{\pi b}{a}.
\end{equation}
Then:
\begin{enumerate}
\item \textbf{(Sector-based conjugation criterion).} Two points of the form $L^{+}_{n}$ or $L^{-}_{n}$ cannot simultaneously lie inside the open interval $\bigl(L_{\min}^{(n)}, L_{\max}^{(n)}\bigr)$. If a point $L^{+}_{n}$ (or $L^{-}_{n}$) lies inside the open sector $s_n$, then the corresponding integer $n_{+}$ (or $n_{-}$) uniquely defines the conjugate series for branch $n$.
\item \textbf{(Boundary conjugation convention).} If no such points lie inside the open sector $s_n$, while points $L^{+}_{n}$ and $L^{-}_{n}$ simultaneously coincide on the shared boundary between sectors $s_{n-1}$ and $s_n$, then the lower branch $n-1$ conjugates with the index $n_{-}$, and the upper branch $n$ conjugates with the index $n_{+}$.
\end{enumerate}

\begin{proof}
1. For the central points $L^{0}_{n}$, the step size is $\Delta_0 = \frac{2\pi b}{a}$; for $L^{+}_{n}$, it is $\Delta_{+} = \frac{2\pi b}{a-b}$; and for $L^{-}_{n}$, it is $\Delta_{-} = 2\pi$. For $a > b > 0$, the strict inequalities $\Delta_{+} > \Delta_0$ and $\Delta_{-} > \Delta_0$ hold; consequently, two points from the same family cannot fall into a single sector of width $\Delta_0$. Let us show that points from different families, $L^{+}_{n}$ and $L^{-}_{n}$, also cannot simultaneously lie inside an open sector. Set $\theta = \arg(x)$; then $\arg(1/(ax)) = -\theta + 2\pi p$ and $\arg(-ax) = \theta + \pi + 2\pi q$ ($p, q \in \mathbb{Z}$). Equating $L^{+}_{n} = L^{-}_{n} = L$ yields
\[
\frac{b}{b-a}(-\theta + 2\pi p + 2\pi n_{+}) = -(\theta + \pi + 2\pi q + 2\pi n_{-}) \implies a\theta = 2\pi b(p+n_{+}) + (b-a)\pi + 2\pi(b-a)(q+n_{-}),
\]
whence $L = -\frac{\pi b}{a}\bigl[2(p+n_{+}+q+n_{-}) + 1\bigr]$, meaning that coincidence is possible exclusively on the sector boundaries $L_{\max}^{(n-1)} = L_{\min}^{(n)}$. Since $L^{+}_{n}$ increases with $\theta$ while $L^{-}_{n}$ decreases, they swap sectors at the moment of boundary coincidence, so the total occupancy of open sectors is independent of $\theta$ and, by the frequency identity
\begin{equation}
\label{eq:freq_sum}
\frac{a}{2\pi b} = \frac{a-b}{2\pi b} + \frac{1}{2\pi},
\end{equation}
equals exactly one point ($L^{+}_{n}$ or $L^{-}_{n}$) per open sector.

2. The behavior on the sector boundaries is fixed by the standard convention for the principal value of the argument $\arg(z) \in (-\pi, \pi]$. This rule ensures complete consistency between the branch numbering of the ultralogarithm, the classical numbering of the $W_0$ and $W_{-1}$ branches of the Lambert $W$ function \cite{Corless1996, DLMF}, and the author's previous works~\cite{Berezin2025_arXiv}.
\end{proof}

\noindent\textbf{Corollary 1 (Logarithmic Analog of Vieta's Theorem and Branch Classification).}
Let the exponents $a > b > 0$ in Eq.~\eqref{eq:ulog_main} be integers (an analogous result holds for the general multinomial equation $e^{aw} = 1 + a\sum_{j} x_j e^{b_j w}$ with integers $a > b_j > 0$). Then:
\begin{enumerate}
\item The ultralogarithm possesses exactly $a$ \textbf{independent branches} $w_0$, $w_1$, \dots, $w_{a-1}$~--- matching the number of roots of the corresponding algebraic ultraradical $y_n = e^{w_n}$.
\item All remaining branches are \textbf{dependent} on them: for any $n, p \in \mathbb{Z}$,
\[
w_{n+ap} = w_n + 2\pi i p,
\]
because the exact same factor $e^{b f^{0}_{n+ap}} = e^{b f^{0}_{n}}$ enters the argument of the master series~\eqref{eq:master_series} for both $m=0$ and $m=1$. Hence, the selection of independent branches is conventional: any complete set of representatives of the residue classes modulo $a$ (such as $\{0,\dots,a-1\}$ or a zero-centered symmetric set) serves as a basis.
\item For the canonical basis $n \in \{0,\dots,a-1\}$, the \textbf{logarithmic Vieta's law} holds~--- a branch-sum invariant independent of $x$ (and of the coefficients $x_j$ in the multinomial case):
\begin{equation}
\label{eq:log_vieta}
\sum_{n=0}^{a-1} w_n = i\pi(a - 1).
\end{equation}
Under a change of basis, the sum is preserved modulo $2\pi i$.
\end{enumerate}

\section{Historical Background}
\label{sec:history}

The tradition of representing solutions to trinomial algebraic and transcendental equations directly as power series, bypassing attempts to express them through finite combinations of radicals, dates back to the works of Johann Lambert (1758)~\cite{Mezo2022} and Leonhard Euler (1779)~\cite{Corless1996}. In that era, the very concept of an analytic function as a unified multivalued object on a Riemann surface was still in its infancy: mathematicians worked primarily with local power expansions without addressing the global matching of all branches.

Over the subsequent two centuries, this line of research expanded substantially. Numerous methods for local series inversion emerged (including the Lagrange inversion formula~\cite{Gessel2016, SuryaWarnke2023}), alongside a vast array of approaches to analytic continuation beyond the disk of convergence\cite{Whittaker1927}~--- ranging from re-expansions along chains of disks and Padé approximants to integral representations and numerical contour methods. In modern works by D.~Belki\'c~\cite{Belkic2018, Belkic2019}, power expansions for the roots of trinomial equations and their logarithms were systematized in terms of Fox--Wright functions and Bell polynomials, while in the author's previous studies~\cite{Berezin2025_arXiv}, all classical expansions of Lambert, Euler, and Belki\'c were unified within the single master series formalism $M(m;a;b;X)$.

Nevertheless, a fundamental question remained open: how to analytically recover \emph{all} branches of a multivalued function outside its radius of convergence $\vert{}x\vert{} < R$ starting from a single-valued power series, by employing transformed power series with different parameters, while rigorously identifying which specific series continues each individual branch. The advent of modern computer algebra systems and large language models (AI) as tools for rapid generation and verification of algebraic hypotheses made it possible to overcome this combinatorial barrier: the sector-based conjugation criterion was discovered, transforming a collection of local ultralogarithm expansions into a unified, globally computable multivalued function.

In the author's previous works~\cite{Berezin2025_arXiv, Berezin2026_arXiv}, the primary focus was placed on studying the ultraradical $y_n$, whereas the ultralogarithm $w_n$ was treated mainly as an auxiliary function. At first glance, passing from the ultraradical to the ultralogarithm via the substitution $y_n = e^{w_n}$ might appear trivial; however, a global analysis of the branches reveals fundamental differences that make an independent investigation of the ultralogarithm essential:
\begin{enumerate}
\item \textbf{Infinite multivaluedness for integer parameters.} Even for integer exponents $a$ and $b$, the ultralogarithm possesses a countably infinite set of branches. An illustrative example is the case $a=2, b=1$: the ultraradical equation $y^2 = 1 + 2xy$ yields only two algebraic roots $y = x \pm \sqrt{x^2+1}$, whereas the corresponding ultralogarithm $w_n = \ln\bigl(x \pm \sqrt{x^2+1}\bigr) + 2\pi i n$ generates an infinite family of branches of the generalized inverse hyperbolic sine ($\operatorname{arsinh}$).
\item \textbf{Geometric stratification of branch trajectories.} Unlike the ultraradical $y_n$, whose projections onto the complex plane form numerous loops and self-intersections at different values of the argument modulus, the trajectories of the ultralogarithm branches $w_n$ are stratified into their own phase strips of width $\frac{2\pi}{a}$: along non-critical rays, the branches do not intersect at all, whereas along critical rays they can only touch at isolated branch points on the shared boundaries of adjacent strips at $\vert{}x\vert{} = R$, without crossing into neighboring strips.
\item \textbf{Computational stability and simplicity of expansions.} When refining roots of the ultraradical numerically via Newton's or Schröder's iterative methods near branch cuts and branch points, standard algorithms frequently lose stability or require explicitly supplying an auxiliary branch-verification parameter $u$ (to select the proper sheet of the fractional powers $y^a$ and $y^b$). For the ultralogarithm equation~\eqref{eq:ulog_main}, which involves only exponential functions, iterative refinement converges smoothly without additional branch-selection constraints. Moreover, as demonstrated in~\cite{Berezin2026_arXiv}, constructing the Taylor series of the ultralogarithm around an arbitrary point $x_0$ is an order of magnitude simpler than for the ultraradical.
\end{enumerate}
For these reasons, the present study is devoted entirely to analyzing the branches of the ultralogarithm, while systematically comparing their behavior in computational experiments against the plots of the well-established ultraradical.

\section{Branch Phase Hodographs}
\label{sec:experiments}

\subsection{Multi-Channel Branch Visualization Methodology}
\label{subsec:visualization}

Traditional real-valued function plots cannot capture the simultaneous evolution of the real and imaginary parts of a multivalued solution in the complex plane, whereas three-dimensional Riemann surfaces become unreadable when a large number of branches are superimposed. In this work, solutions are analyzed using \textbf{phase hodographs with multi-channel encoding}:
\begin{itemize}
\item \textbf{Planar coordinates}~--- the real and imaginary parts of the solution itself: $(\mathrm{Re}\,w, \mathrm{Im}\,w)$ for the ultralogarithm and $(\mathrm{Re}\,y, \mathrm{Im}\,y)$ for the ultraradical;
\item \textbf{Color coding and radius of convergence indication}~--- the trajectory segment inside the disk of convergence ($\vert{}x\vert{} < R$), computed via the basic expansion~\eqref{eq:ulog_inner}, is plotted in \textbf{blue}; crossing the boundary $\vert{}x\vert{} = R$ and switching to the conjugate series~\eqref{eq:ulog_outer_h} or~\eqref{eq:ulog_outer_k} is highlighted by a \textbf{thicker green segment}; further evolution of the branch in the outer region ($\vert{}x\vert{} > R$) as $\vert{}x\vert{}$ increases to $x_{\max}$ is represented by a \textbf{yellow-to-red gradient};
\item \textbf{Line thickness}~--- the branch index $n$: for unambiguous visual identification, each successive branch is plotted with an increasing line thickness.
\end{itemize}
For the ultraradical, a sparse dashed unit circle $\vert{}y\vert{} = 1$ is additionally displayed in the plots.

\textbf{Base points and geometric branching structure.} A key property of the phase hodographs is that at $\vert{}x\vert{} = 0$ (the start of the blue trajectory segment), all branches originate from fixed base points independent of the ray direction $\arg(x)$. These base points reflect the fundamental difference between the two projections of the master equation:

\begin{itemize}
\item \textbf{Ultralogarithm} ($m=0$): the base points $w_n(0) = \frac{2\pi i n}{a}$ lie \textbf{on the imaginary axis} with a uniform spacing of $\frac{2\pi}{a}$. This corresponds to an \textbf{additive branching structure}: the branches are ``unfolded'' along the imaginary axis and evolve within their own isolated phase strips.

\item \textbf{Ultraradical} ($m=1$): the base points $y_n(0) = e^{2\pi i n/a}$ lie \textbf{on the unit circle} $\vert{}y\vert{} = 1$ with an angular spacing of $\frac{2\pi}{a}$. This corresponds to a \textbf{multiplicative branching structure}: the branches are ``wrapped'' onto a circle.
\end{itemize}

The connection between the two representations is given by the exponential mapping $y_n = e^{w_n}$, which at the level of master series is expressed by the fundamental power series identity
\begin{equation}
\label{eq:master_exp_id}
M(1;a;b;X) = e^{M(0;a;b;X)}.
\end{equation}
The very existence of such an identity between two power series reveals a deep structural connection between this pair of functions: it ``wraps'' the imaginary axis of the ultralogarithm ($m=0$) onto the unit circle of the ultraradical ($m=1$). It is precisely this ``wrapping'' that masks the topological branching structure in the ultraradical and makes the ultralogarithm a more natural object for analyzing branch continuity and verifying the deterministic analytic continuation criterion.

In the phase hodographs, this difference is immediately apparent: the ultralogarithm trajectories remain localized within their phase strips near the levels $\mathrm{Im}(w) \approx \frac{2\pi n}{a}$ (touching only at isolated branch points when the ray hits a critical angle), whereas the ultraradical trajectories form loops and self-intersections upon crossing the boundary $\vert{}x\vert{} = R$.

\subsection{Benchmark Case of Integer Parameters: a = 2, b = 1}
\label{subsec:a2_b1}

It is natural to begin the experimental analysis with the simplest special case $a = 2$, $b = 1$, for which Eq.~\eqref{eq:ulog_main} takes the form $e^{2w} = 1 + 2x e^w$ and the radius of convergence of the basic series~\eqref{eq:radius_R} is $R = 1$. In this classical limit, the ultraradical ($m = 1$) describes the two algebraic branches of the quadratic equation $y^2 - 2xy - 1 = 0$, whereas the ultralogarithm ($m = 0$) yields the countably infinite family of branches of the inverse hyperbolic sine ($\operatorname{arsinh}$). Comparing the two on this benchmark example provides a clear verification of the sector-based conjugation criterion against known elementary functions.

\begin{figure}[htbp]
    \centering
    \begin{minipage}[b]{0.48\textwidth}
        \centering
        \includegraphics[width=\textwidth]{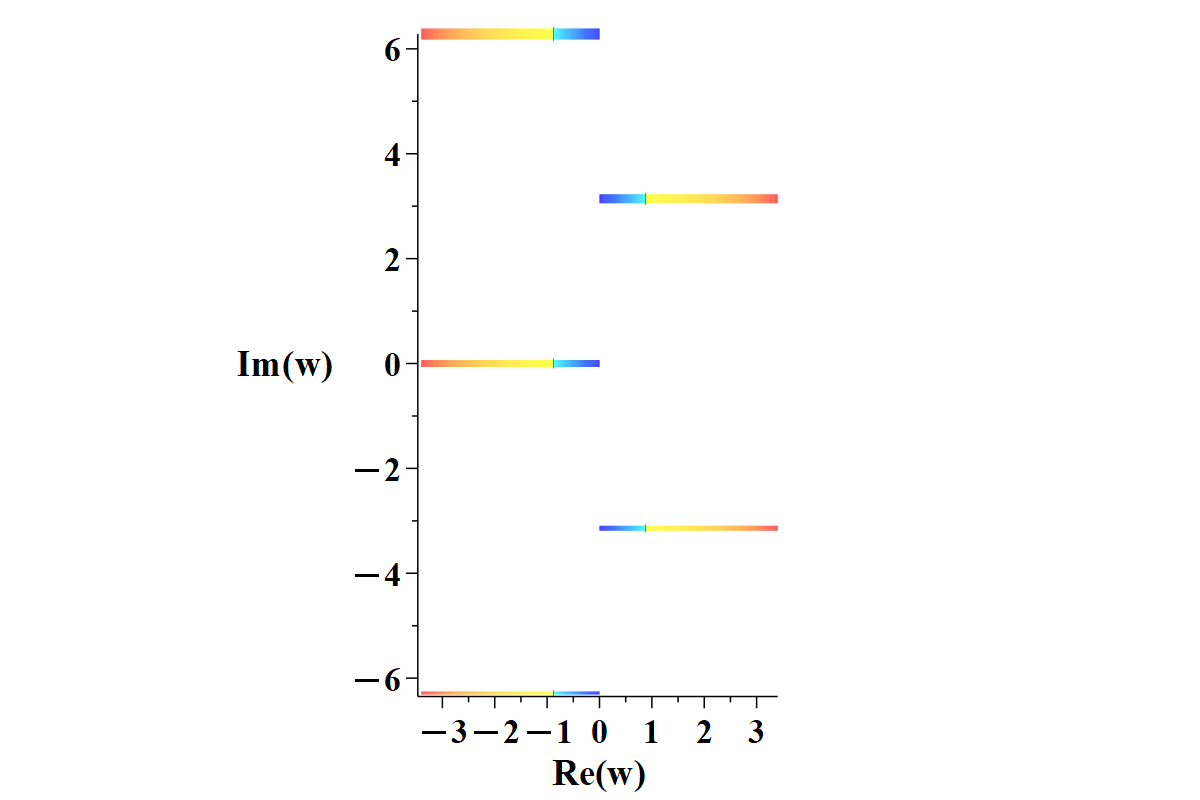}
        \textbf{(a)}
    \end{minipage}
    \hfill
    \begin{minipage}[b]{0.48\textwidth}
        \centering
        \includegraphics[width=\textwidth]{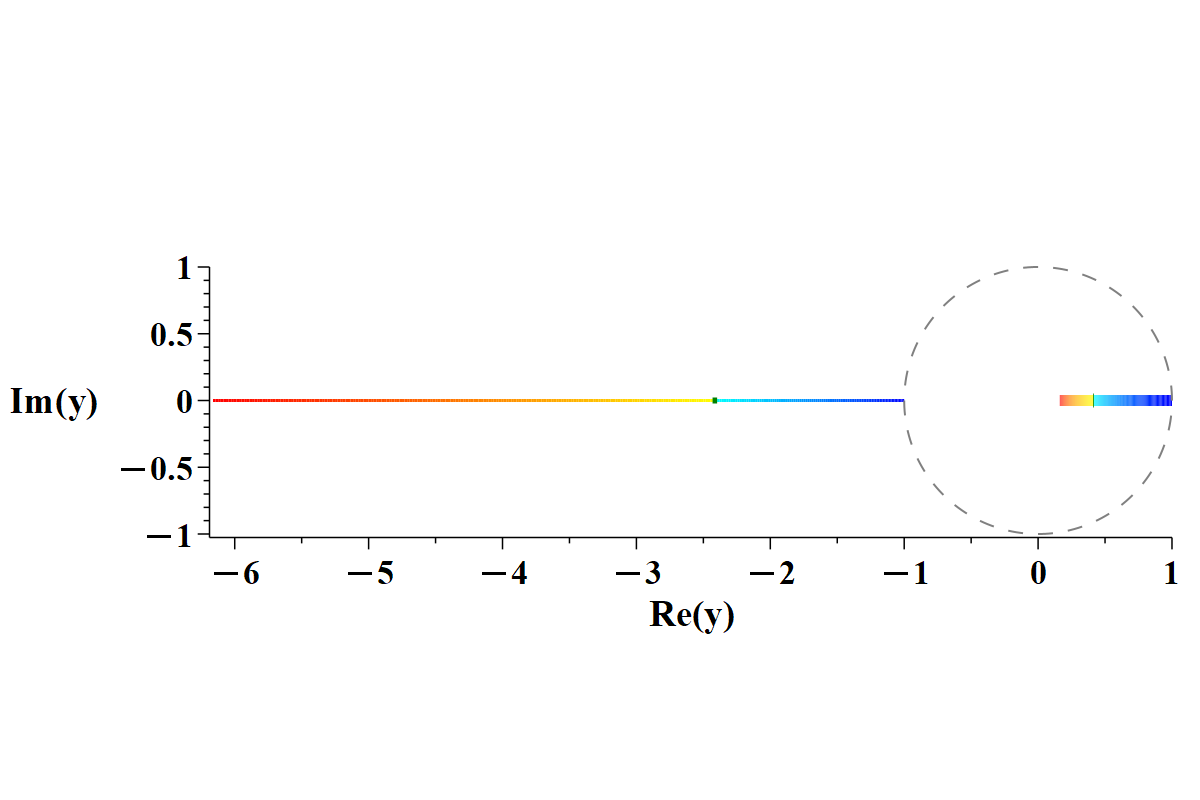}
        \textbf{(b)}
    \end{minipage}
    \caption{Phase hodographs for $a = 2$, $b = 1$ at angle $\arg(x) = -\pi$: 
    \textbf{(a)}~ultralogarithm ($m=0$, equation $e^{2w} = 1 + 2x e^{w}$), branches $n \in \{-2, \dots, 2\}$, $x_{\max} = 15$; 
    \textbf{(b)}~ultraradical ($m=1$, equation $y^{2} = 1 + 2xy$), branches $n \in \{-1, 0\}$, $x_{\max} = 3$.}
    \label{fig:2_1_minus_pi}
\end{figure}

\begin{figure}[htbp]
    \centering
    \begin{minipage}[b]{0.48\textwidth}
        \centering
        \includegraphics[width=\textwidth]{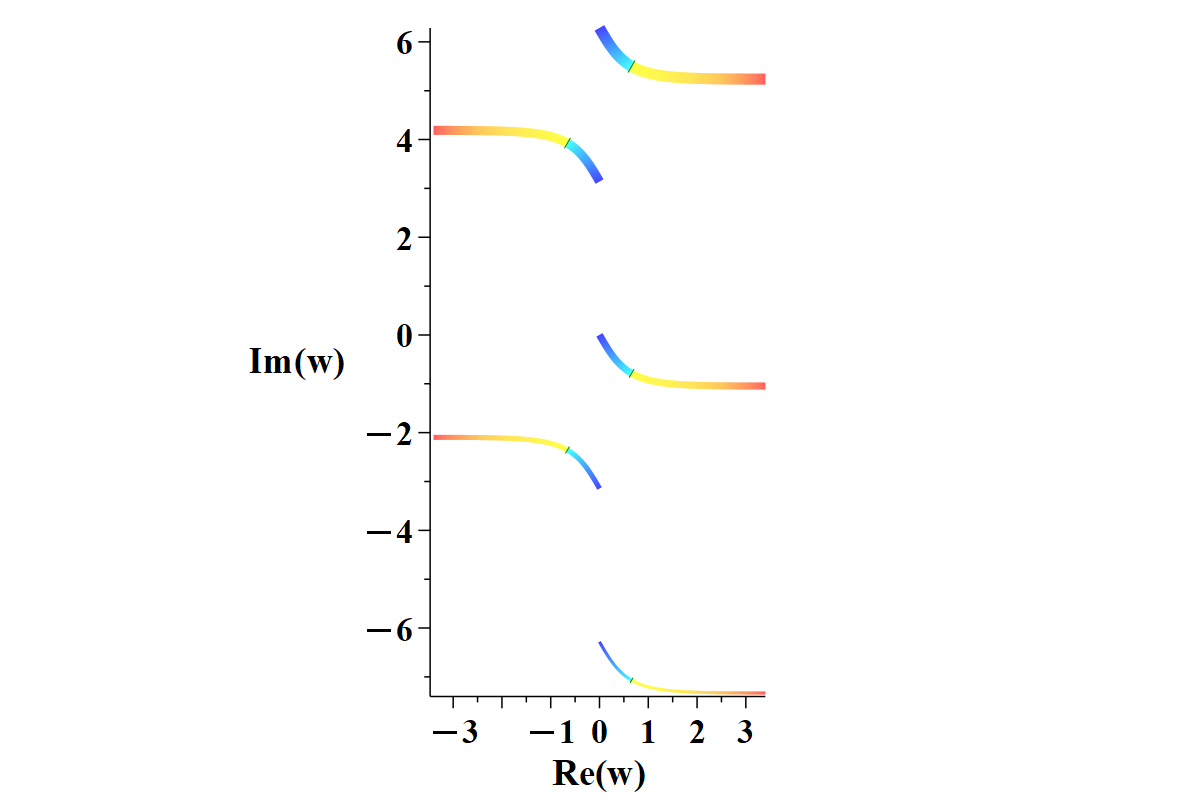}
        \textbf{(a)}
    \end{minipage}
    \hfill
    \begin{minipage}[b]{0.48\textwidth}
        \centering
        \includegraphics[width=\textwidth]{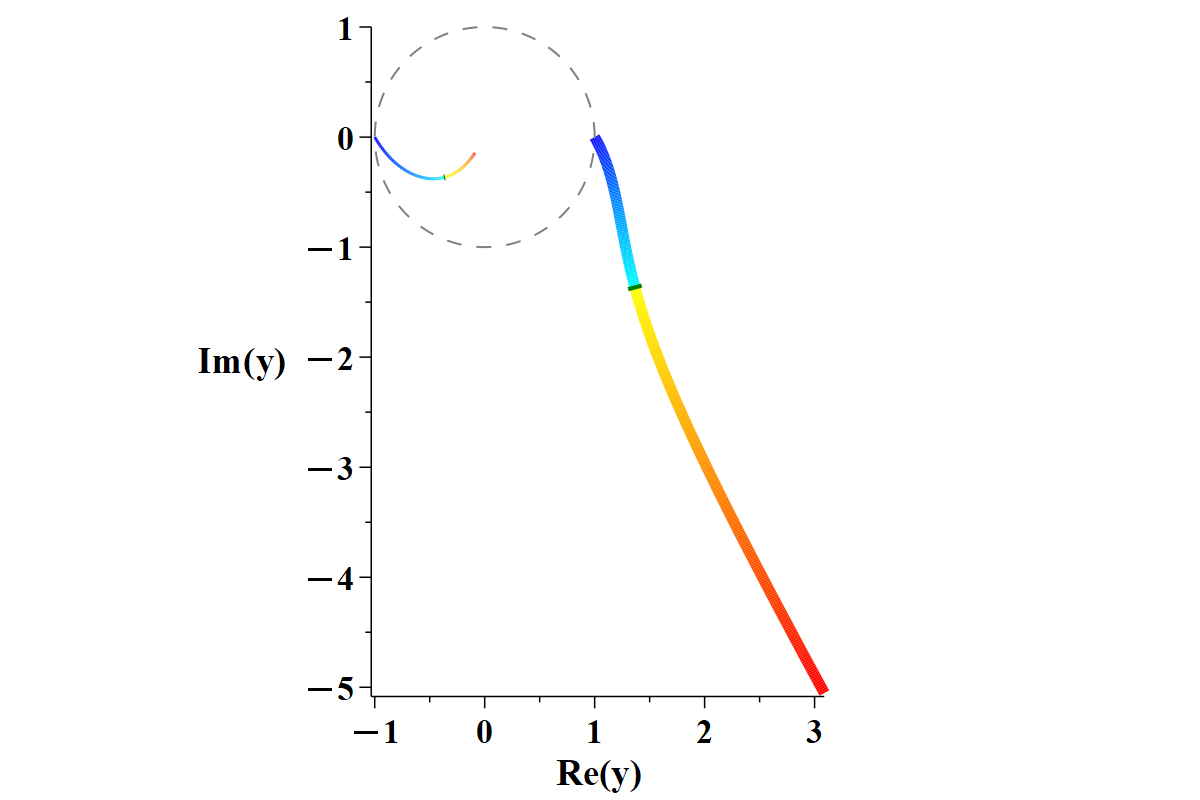}
        \textbf{(b)}
    \end{minipage}
    \caption{Phase hodographs for $a = 2$, $b = 1$ at angle $\arg(x) = -\pi/3$: 
    \textbf{(a)}~ultralogarithm ($m=0$, equation $e^{2w} = 1 + 2x e^{w}$), branches $n \in \{-2, \dots, 2\}$, $x_{\max} = 15$; 
    \textbf{(b)}~ultraradical ($m=1$, equation $y^{2} = 1 + 2xy$), branches $n \in \{-1, 0\}$, $x_{\max} = 3$.}
    \label{fig:2_1_minus_pi3}
\end{figure}

\begin{figure}[htbp]
    \centering
    \begin{minipage}[b]{0.48\textwidth}
        \centering
        \includegraphics[width=\textwidth]{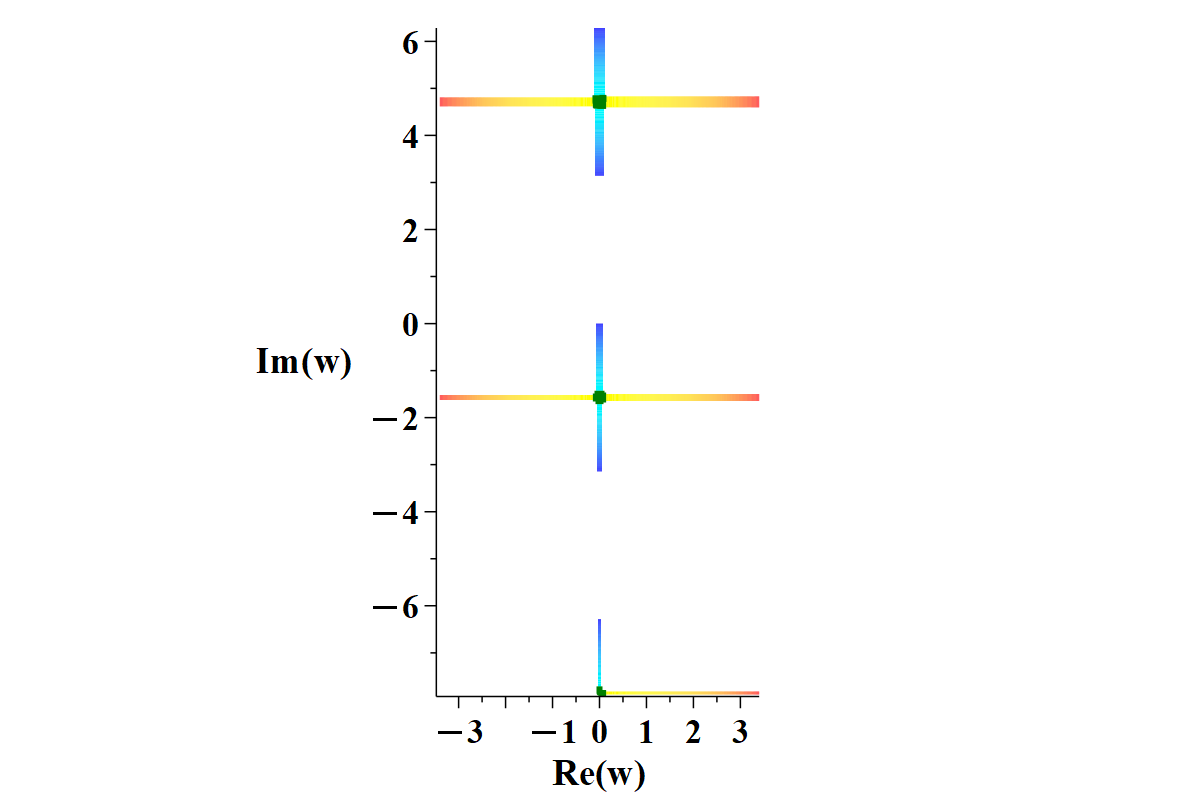}
        \textbf{(a)}
    \end{minipage}
    \hfill
    \begin{minipage}[b]{0.48\textwidth}
        \centering
        \includegraphics[width=\textwidth]{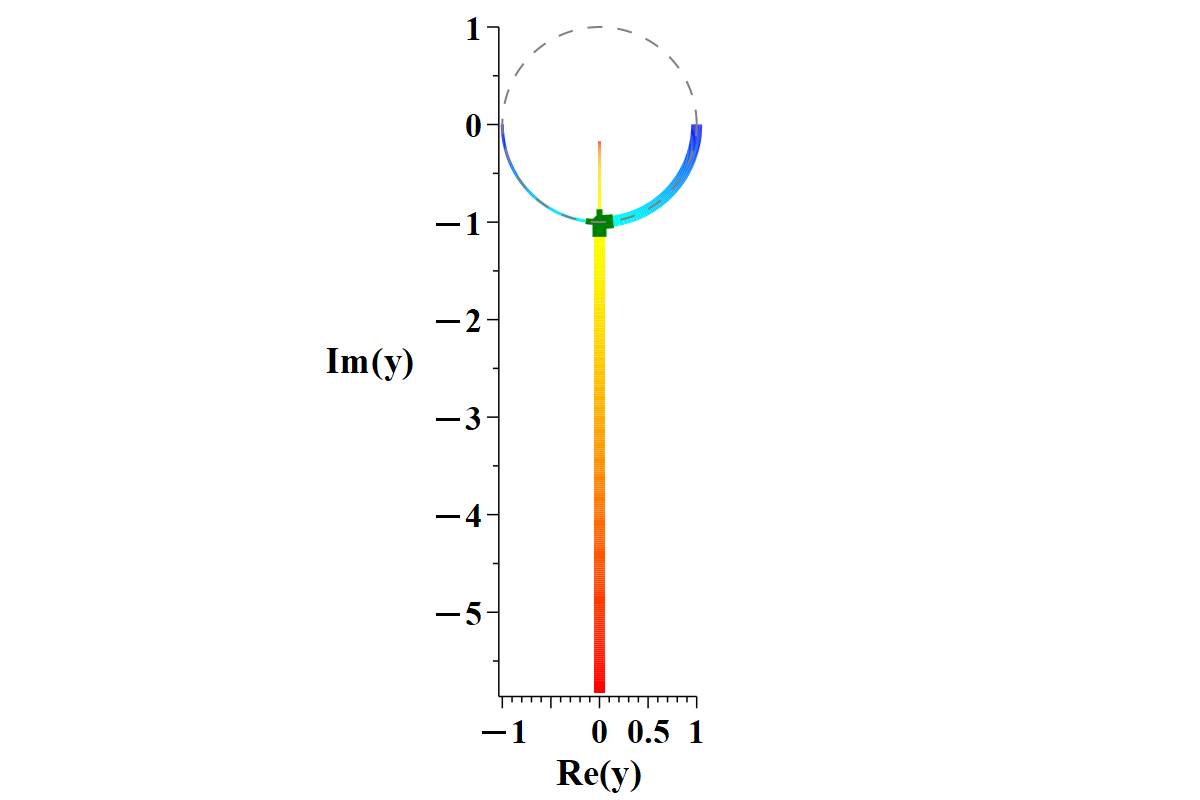}
        \textbf{(b)}
    \end{minipage}
    \caption{Phase hodographs for $a = 2$, $b = 1$ at angle $\arg(x) = -\pi/2$: 
    \textbf{(a)}~ultralogarithm ($m=0$, equation $e^{2w} = 1 + 2x e^{w}$), branches $n \in \{-2, \dots, 2\}$, $x_{\max} = 15$; 
    \textbf{(b)}~ultraradical ($m=1$, equation $y^{2} = 1 + 2xy$), branches $n \in \{-1, 0\}$, $x_{\max} = 3$.}
    \label{fig:2_1_minus_pi2}
\end{figure}

\begin{figure}[htbp]
    \centering
    \begin{minipage}[b]{0.48\textwidth}
        \centering
        \includegraphics[width=\textwidth]{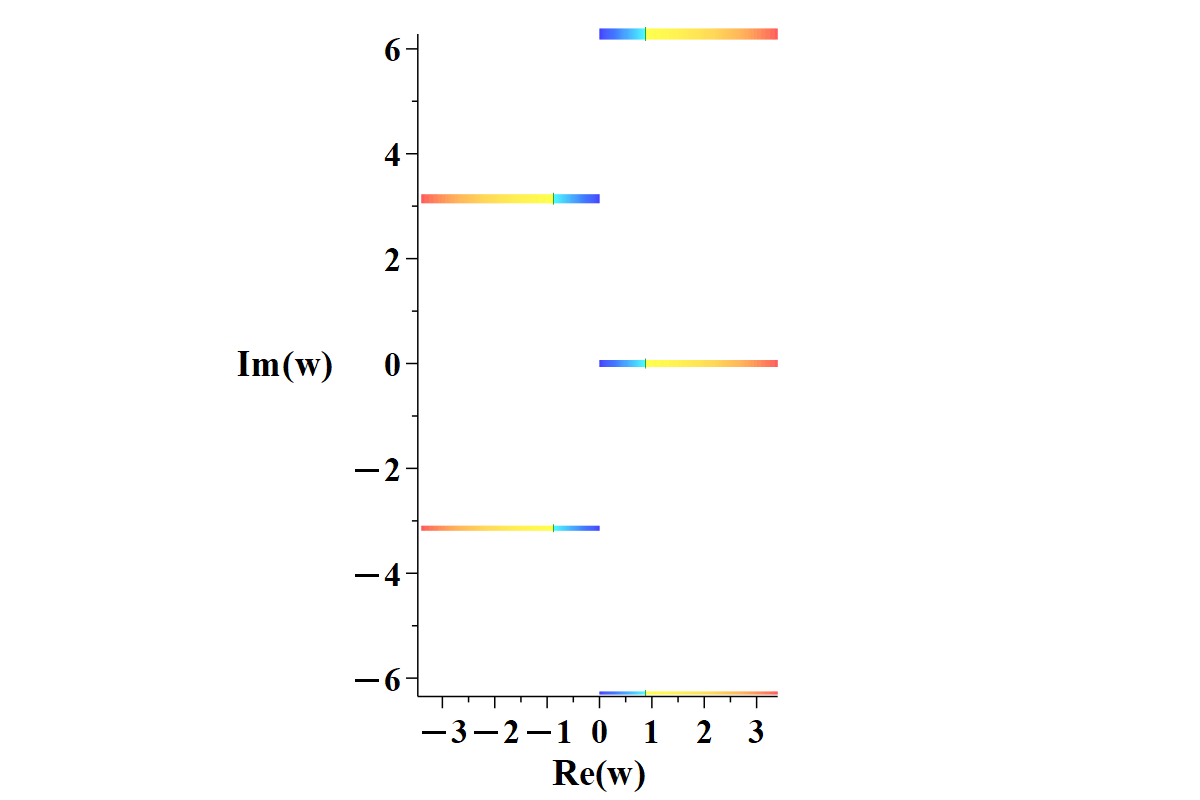}
        \textbf{(a)}
    \end{minipage}
    \hfill
    \begin{minipage}[b]{0.48\textwidth}
        \centering
        \includegraphics[width=\textwidth]{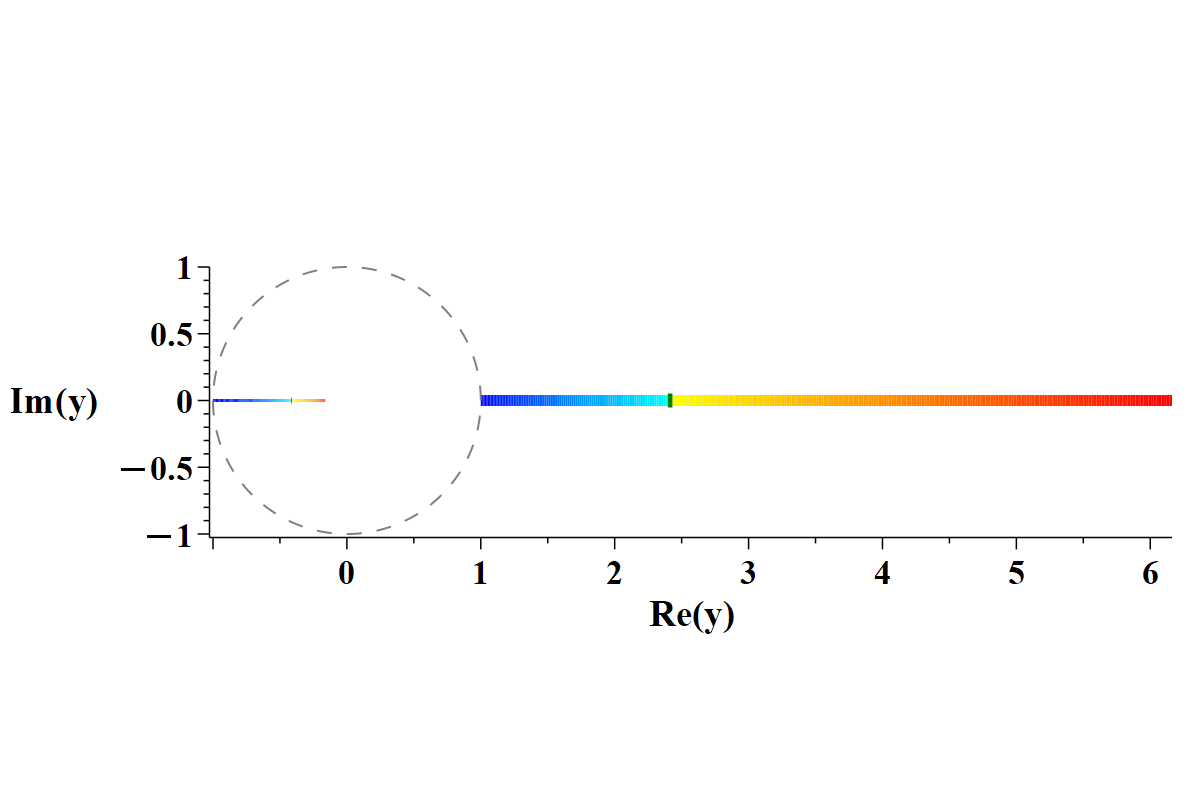}
        \textbf{(b)}
    \end{minipage}
    \caption{Phase hodographs for $a = 2$, $b = 1$ at angle $\arg(x) = 0$: 
    \textbf{(a)}~ultralogarithm ($m=0$, equation $e^{2w} = 1 + 2x e^{w}$), branches $n \in \{-2, \dots, 2\}$, $x_{\max} = 15$; 
    \textbf{(b)}~ultraradical ($m=1$, equation $y^{2} = 1 + 2xy$), branches $n \in \{-1, 0\}$, $x_{\max} = 3$.}
    \label{fig:2_1_zero}
\end{figure}

\begin{figure}[htbp]
    \centering
    \begin{minipage}[b]{0.48\textwidth}
        \centering
        \includegraphics[width=\textwidth]{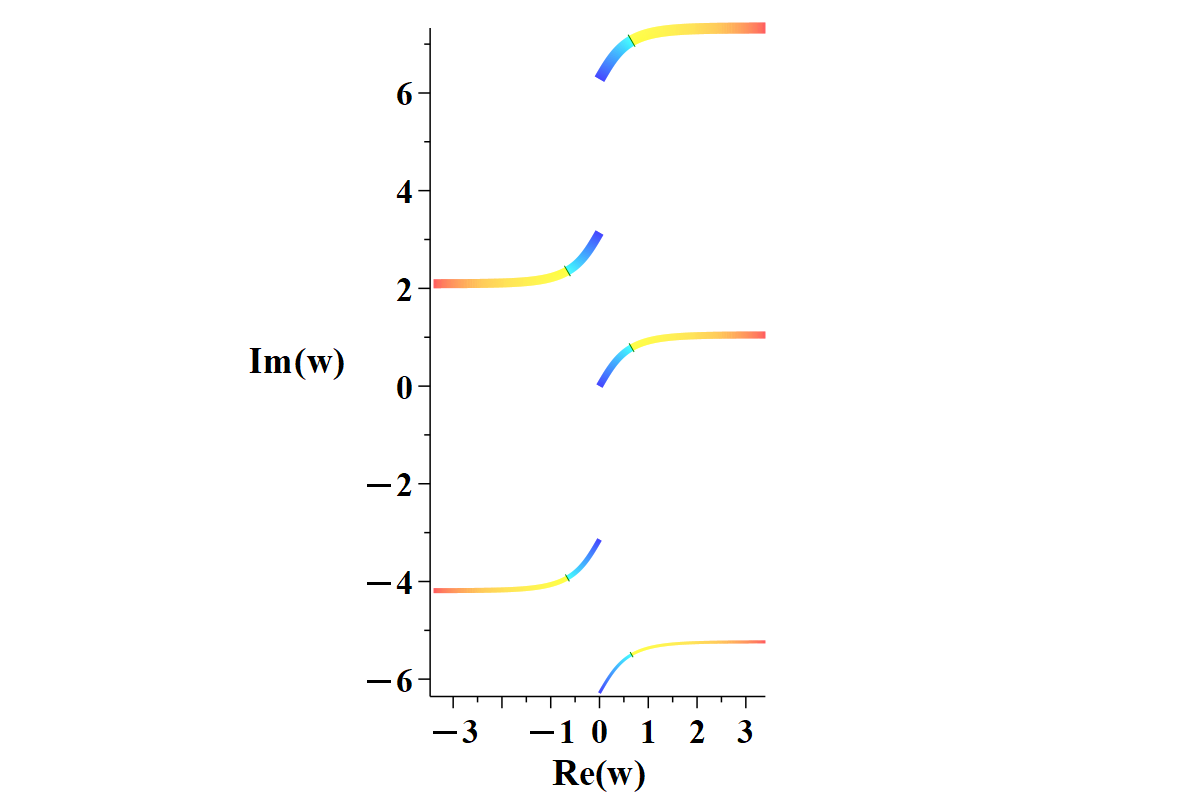}
        \textbf{(a)}
    \end{minipage}
    \hfill
    \begin{minipage}[b]{0.48\textwidth}
        \centering
        \includegraphics[width=\textwidth]{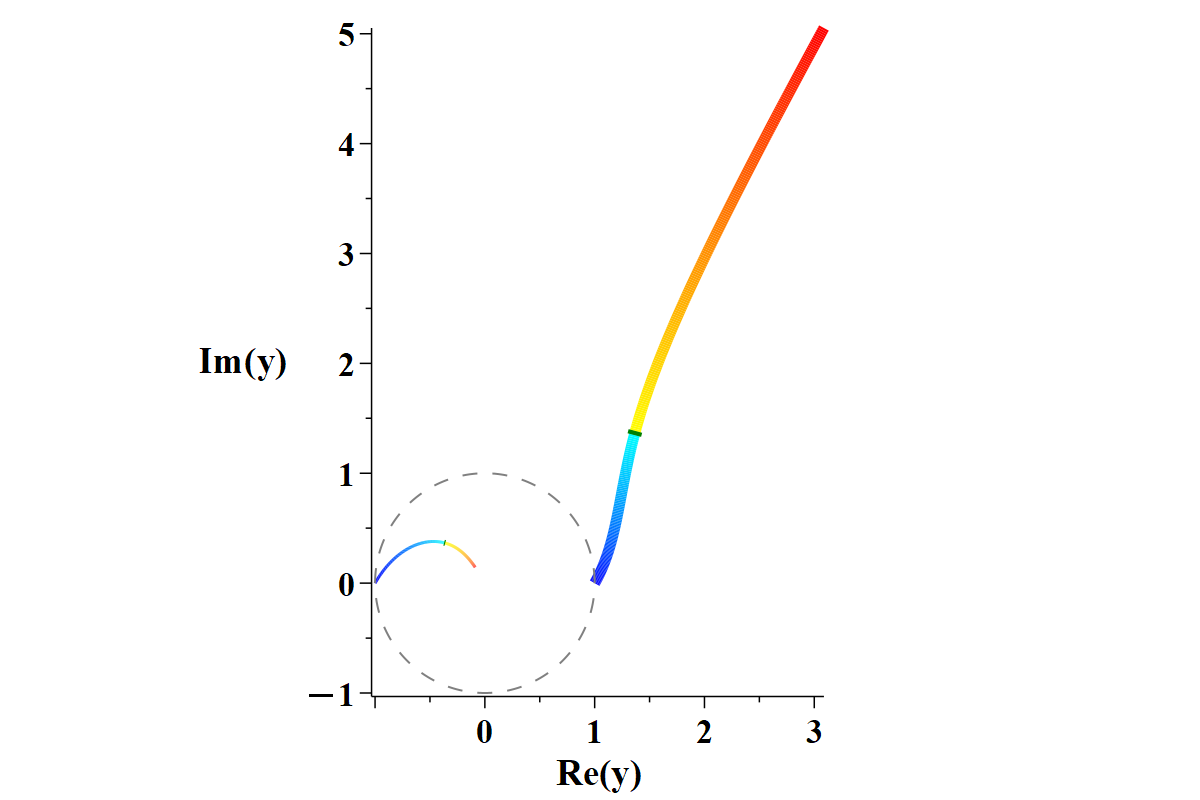}
        \textbf{(b)}
    \end{minipage}
    \caption{Phase hodographs for $a = 2$, $b = 1$ at angle $\arg(x) = \pi/3$: 
    \textbf{(a)}~ultralogarithm ($m=0$, equation $e^{2w} = 1 + 2x e^{w}$), branches $n \in \{-2, \dots, 2\}$, $x_{\max} = 15$; 
    \textbf{(b)}~ultraradical ($m=1$, equation $y^{2} = 1 + 2xy$), branches $n \in \{-1, 0\}$, $x_{\max} = 3$.}
    \label{fig:2_1_plus_pi3}
\end{figure}

Figures~\ref{fig:2_1_minus_pi}--\ref{fig:2_1_plus_pi3} show the paired phase hodographs of the ultralogarithm ($n \in \{-2, \dots, 2\}$, $x_{\max} = 15$) and the ultraradical ($n \in \{-1, 0\}$, $x_{\max} = 3$) for five characteristic ray directions $\arg(x) \in \{-\pi, -\pi/3, -\pi/2, 0, \pi/3\}$. The resulting plots reveal a striking dependence of the straightness of the trajectories on the angle $\arg(x)$:
\begin{itemize}
\item \textbf{Angles that are multiples of $\pi / a = \pi / 2$ (separation of real and imaginary parts).} At even multiples $\arg(x) = -\pi$ and $\arg(x) = 0$ (Figs.~\ref{fig:2_1_minus_pi}a and~\ref{fig:2_1_zero}a), only the real part $\mathrm{Re}(w_n)$ depends on the modulus $\vert{}x\vert{}$, while the imaginary part remains strictly constant ($\mathrm{Im}(w_n) = \pi n$), forming horizontal straight lines. At the odd multiple $\arg(x) = -\pi/2$ (Fig.~\ref{fig:2_1_minus_pi2}a), the roles swap: inside the disk of convergence ($\vert{}x\vert{} < 1$), the real part remains constant ($\mathrm{Re}(w_n) = 0$) and only the imaginary part varies (vertical segments), but after adjacent branches touch at isolated branch points on the boundary $\vert{}x\vert{} = R = 1$ (green markers), the trajectories turn by $90^\circ$~--- thereafter, the imaginary part freezes at $\mathrm{Im}(w_n) = \frac{\pi}{2}(2n \pm 1)$ while only $\mathrm{Re}(w_n)$ grows.
\item \textbf{Intermediate angles.} Away from these multiples, such as at $\arg(x) = -\pi/3$ and $\arg(x) = \pi/3$ (Figs.~\ref{fig:2_1_minus_pi3} and~\ref{fig:2_1_plus_pi3}), both components of the solution vary simultaneously with $\vert{}x\vert{}$, causing the trajectories to bend into mirror-symmetric smooth arcs while the ultralogarithm branches remain strictly isolated within their respective phase strips.
\end{itemize}

\subsection{Quintic Equation with Integer Parameters: a = 5, b = 1}
\label{subsec:a5_b1}

The next stage of the computational experiment examines the quintic equation unsolvable in radicals ($a = 5$, $b = 1$), for which the governing relation~\eqref{eq:ulog_main} takes the form
\begin{equation}
\label{eq:a5_b1}
e^{5w} = 1 + 5x e^{w},
\end{equation}
and the corresponding ultraradical $y = e^w$ satisfies the reduced Bring--Jerrard equation $y^5 - 5xy - 1 = 0$. The radius of convergence of the basic series~\eqref{eq:radius_R} is $R = 4^{-4/5} \approx 0.329877$. Since the branch point satisfies $y^5 = -1/4$ and $x = y^4$, the root modulus at the branch point is $\vert{}y\vert{} = 4^{-1/5} \approx 0.758 < 1$ (corresponding to $\mathrm{Re}(w) = -\frac{1}{5}\ln 4 \approx -0.277$), while the critical rays in the $x$-plane are determined by even multiples of the angle $\frac{\pi}{5}$: $\arg(x) \in \{0, \pm 2\pi/5, \pm 4\pi/5\}$. By Corollary~1, the equation possesses $5$ independent branches satisfying the sum invariant $\sum_{n=-2}^{2} w_n \equiv 4\pi i \pmod{2\pi i}$.

\begin{figure}[htbp]
    \centering
    \begin{minipage}[b]{0.48\textwidth}
        \centering
        \includegraphics[width=\textwidth]{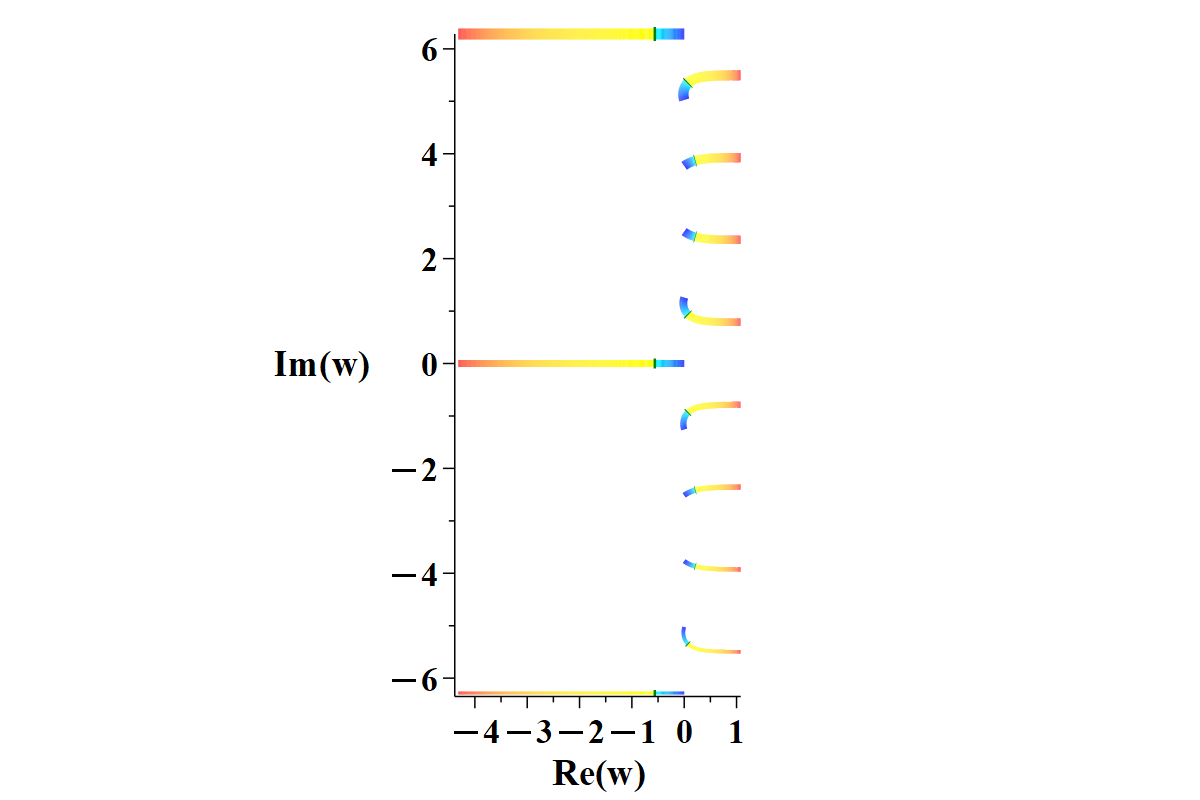}
        \textbf{(a)}
    \end{minipage}
    \hfill
    \begin{minipage}[b]{0.48\textwidth}
        \centering
        \includegraphics[width=\textwidth]{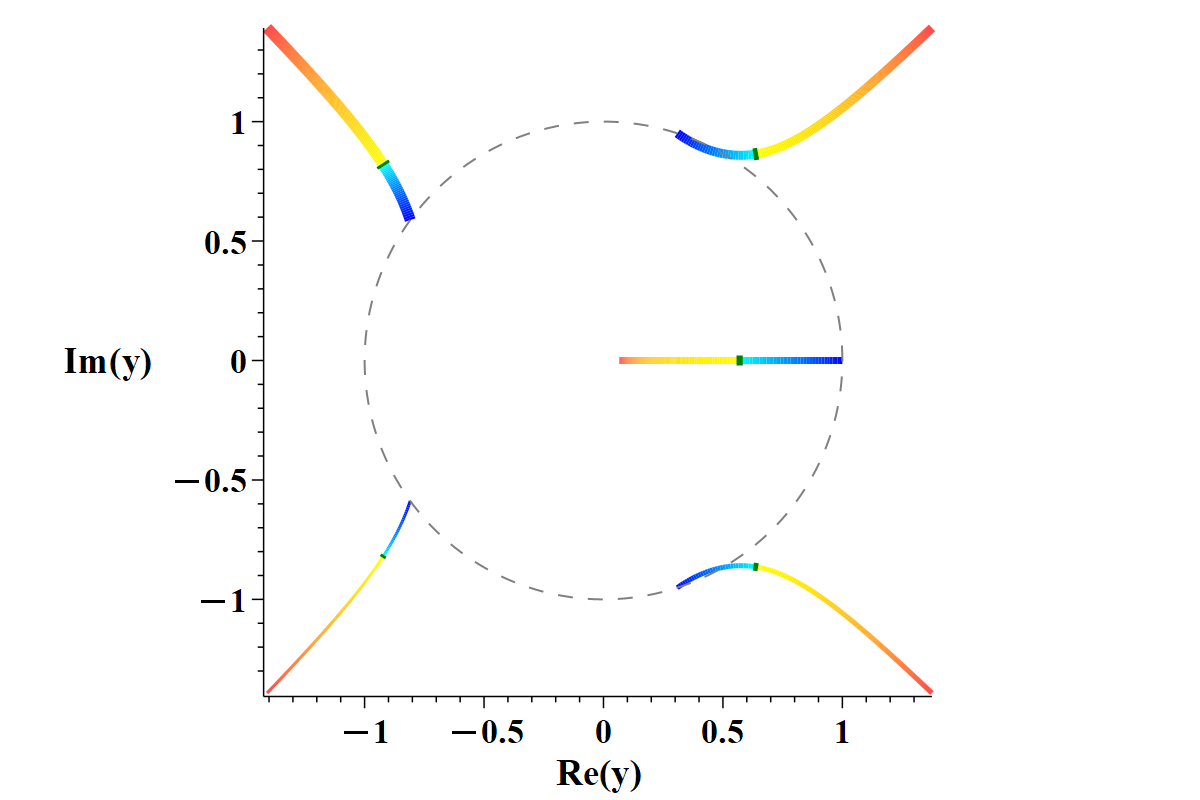}
        \textbf{(b)}
    \end{minipage}
    \caption{Phase hodographs for $a = 5$, $b = 1$ at angle $\arg(x) = -\pi$: 
    \textbf{(a)}~ultralogarithm ($m=0$, equation $e^{5w} = 1 + 5x e^{w}$), branches $n \in \{-5, \dots, 5\}$, $x_{\max} = 15$; 
    \textbf{(b)}~ultraradical ($m=1$, equation $y^{5} = 1 + 5xy$), branches $n \in \{-2, \dots, 2\}$, $x_{\max} = 3$.}
    \label{fig:5_1_minus_pi}
\end{figure}

\begin{figure}[htbp]
    \centering
    \begin{minipage}[b]{0.48\textwidth}
        \centering
        \includegraphics[width=\textwidth]{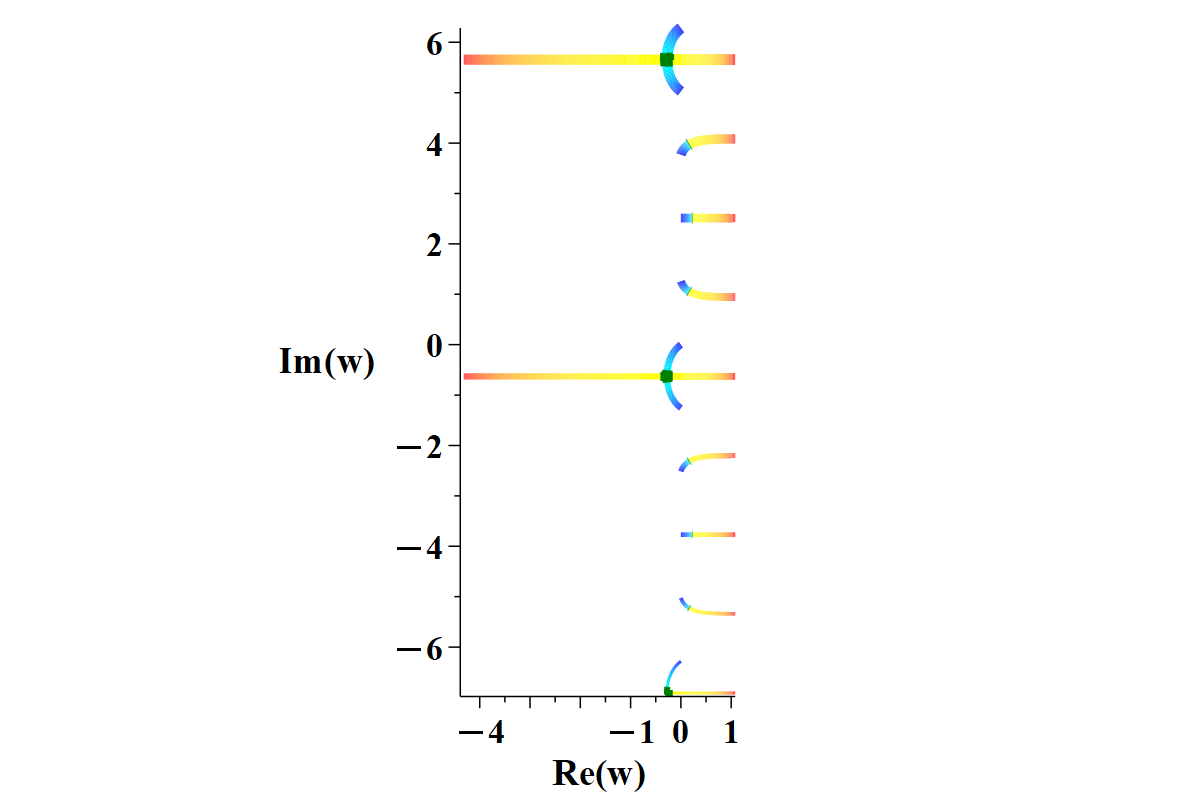}
        \textbf{(a)}
    \end{minipage}
    \hfill
    \begin{minipage}[b]{0.48\textwidth}
        \centering
        \includegraphics[width=\textwidth]{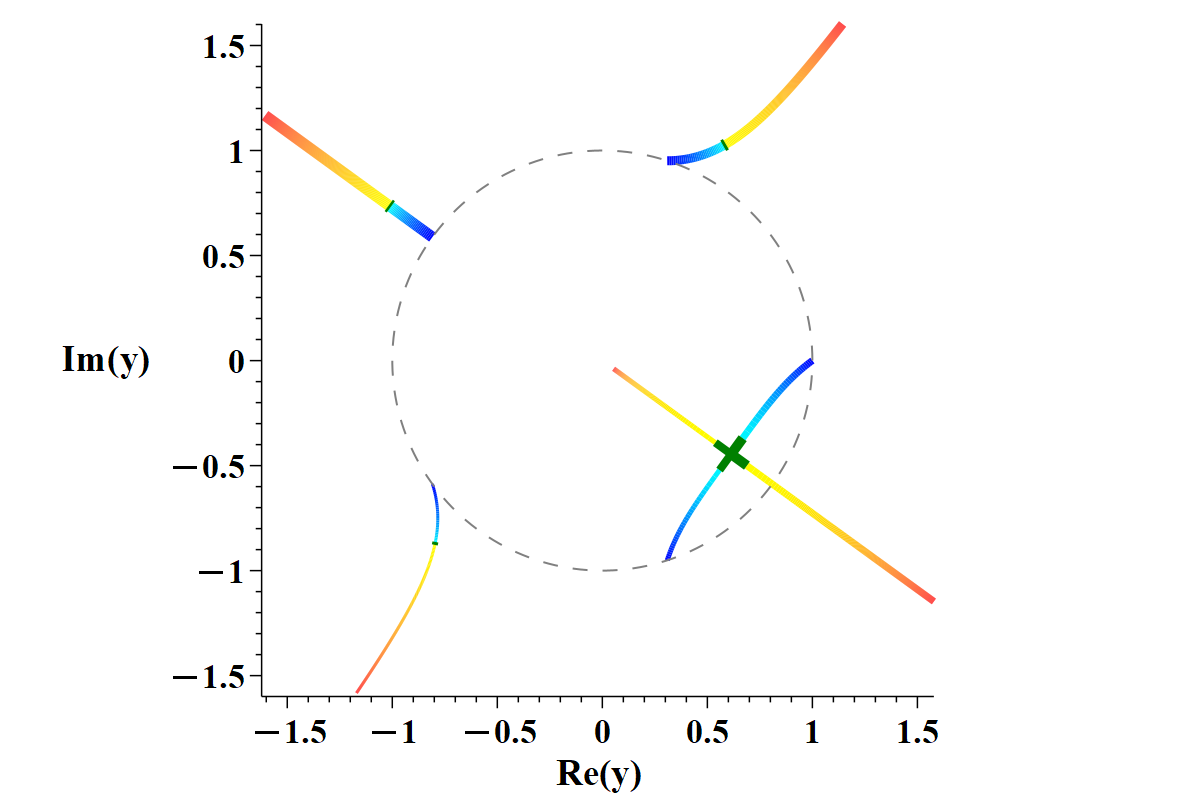}
        \textbf{(b)}
    \end{minipage}
    \caption{Phase hodographs for $a = 5$, $b = 1$ at angle $\arg(x) = -4\pi/5$: 
    \textbf{(a)}~ultralogarithm ($m=0$, equation $e^{5w} = 1 + 5x e^{w}$), branches $n \in \{-5, \dots, 5\}$, $x_{\max} = 15$; 
    \textbf{(b)}~ultraradical ($m=1$, equation $y^{5} = 1 + 5xy$), branches $n \in \{-2, \dots, 2\}$, $x_{\max} = 3$.}
    \label{fig:5_1_minus_4pi5}
\end{figure}

\begin{figure}[htbp]
    \centering
    \begin{minipage}[b]{0.48\textwidth}
        \centering
        \includegraphics[width=\textwidth]{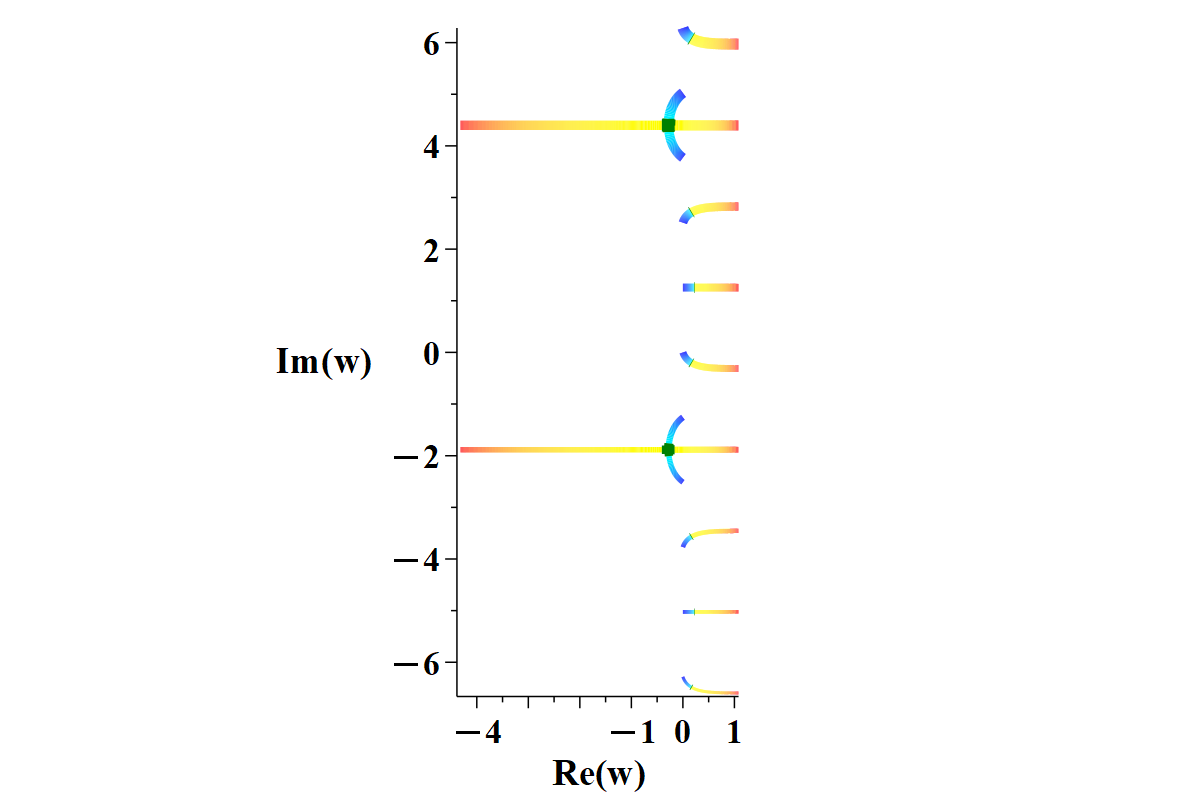}
        \textbf{(a)}
    \end{minipage}
    \hfill
    \begin{minipage}[b]{0.48\textwidth}
        \centering
        \includegraphics[width=\textwidth]{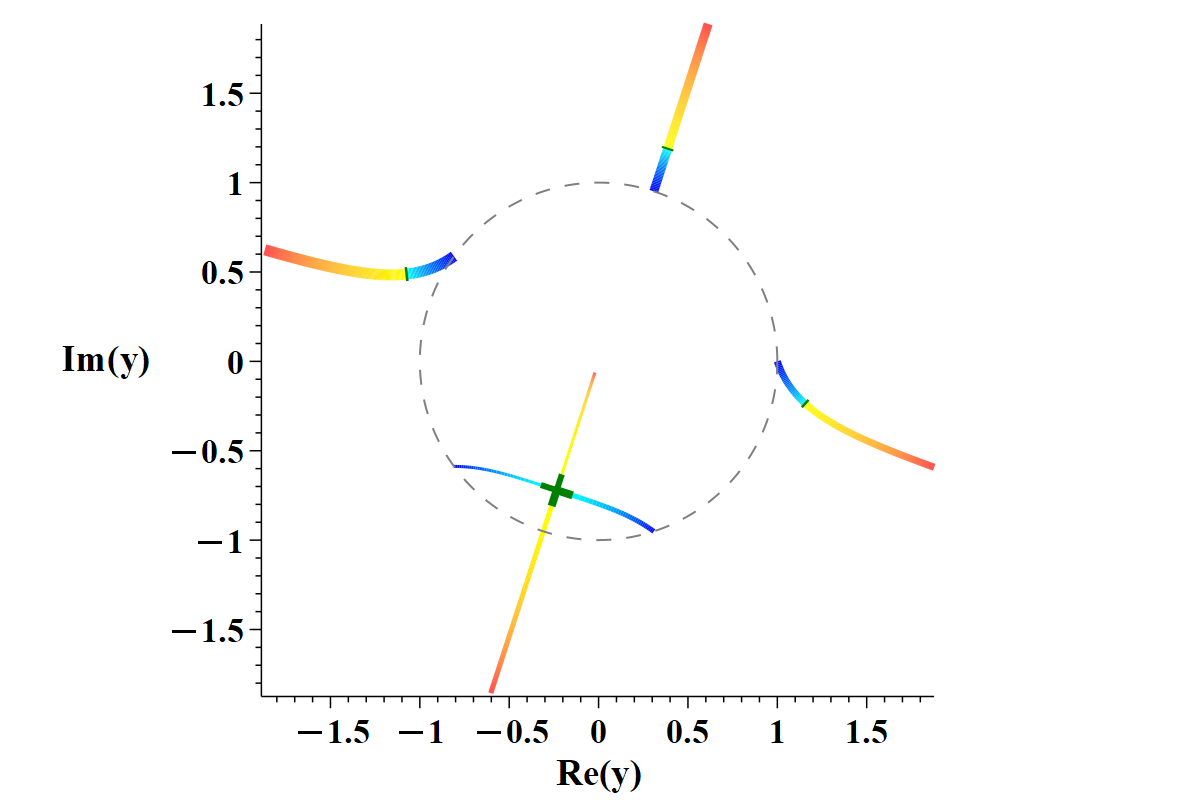}
        \textbf{(b)}
    \end{minipage}
    \caption{Phase hodographs for $a = 5$, $b = 1$ at angle $\arg(x) = -2\pi/5$: 
    \textbf{(a)}~ultralogarithm ($m=0$, equation $e^{5w} = 1 + 5x e^{w}$), branches $n \in \{-5, \dots, 5\}$, $x_{\max} = 15$; 
    \textbf{(b)}~ultraradical ($m=1$, equation $y^{5} = 1 + 5xy$), branches $n \in \{-2, \dots, 2\}$, $x_{\max} = 3$.}
    \label{fig:5_1_minus_2pi5}
\end{figure}

\begin{figure}[htbp]
    \centering
    \begin{minipage}[b]{0.48\textwidth}
        \centering
        \includegraphics[width=\textwidth]{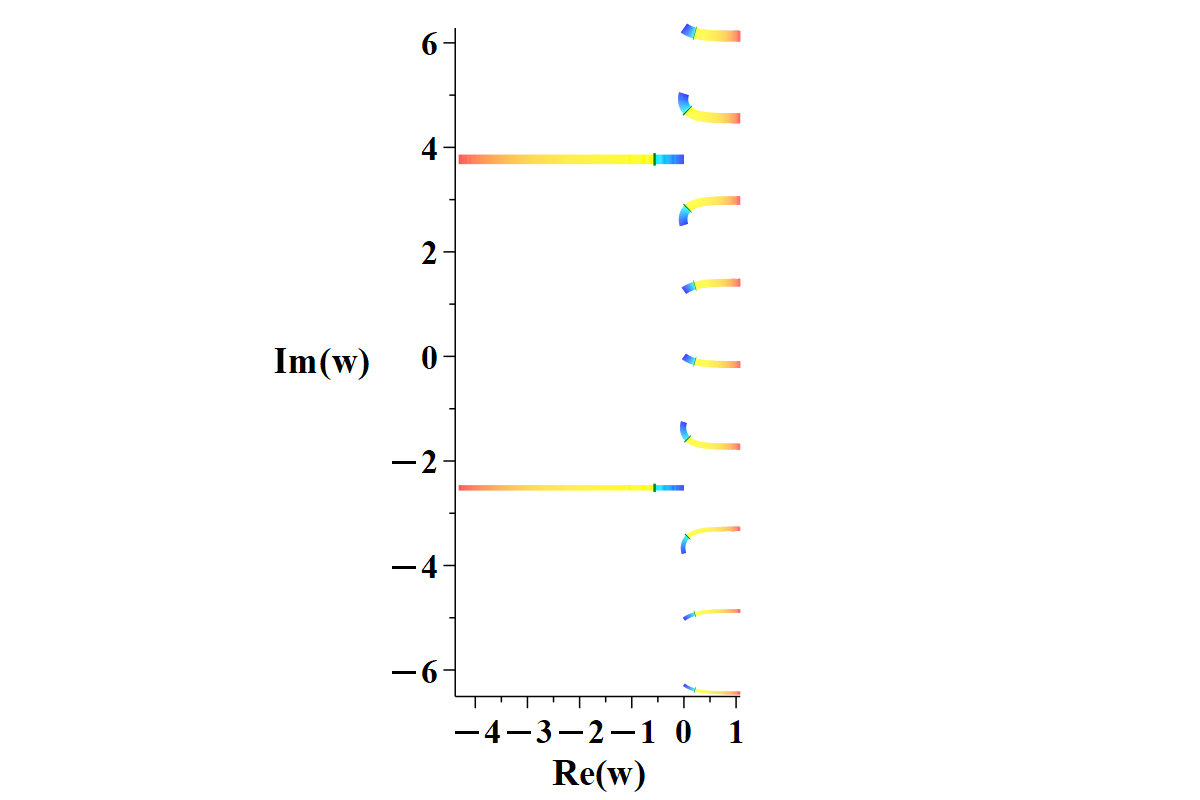}
        \textbf{(a)}
    \end{minipage}
    \hfill
    \begin{minipage}[b]{0.48\textwidth}
        \centering
        \includegraphics[width=\textwidth]{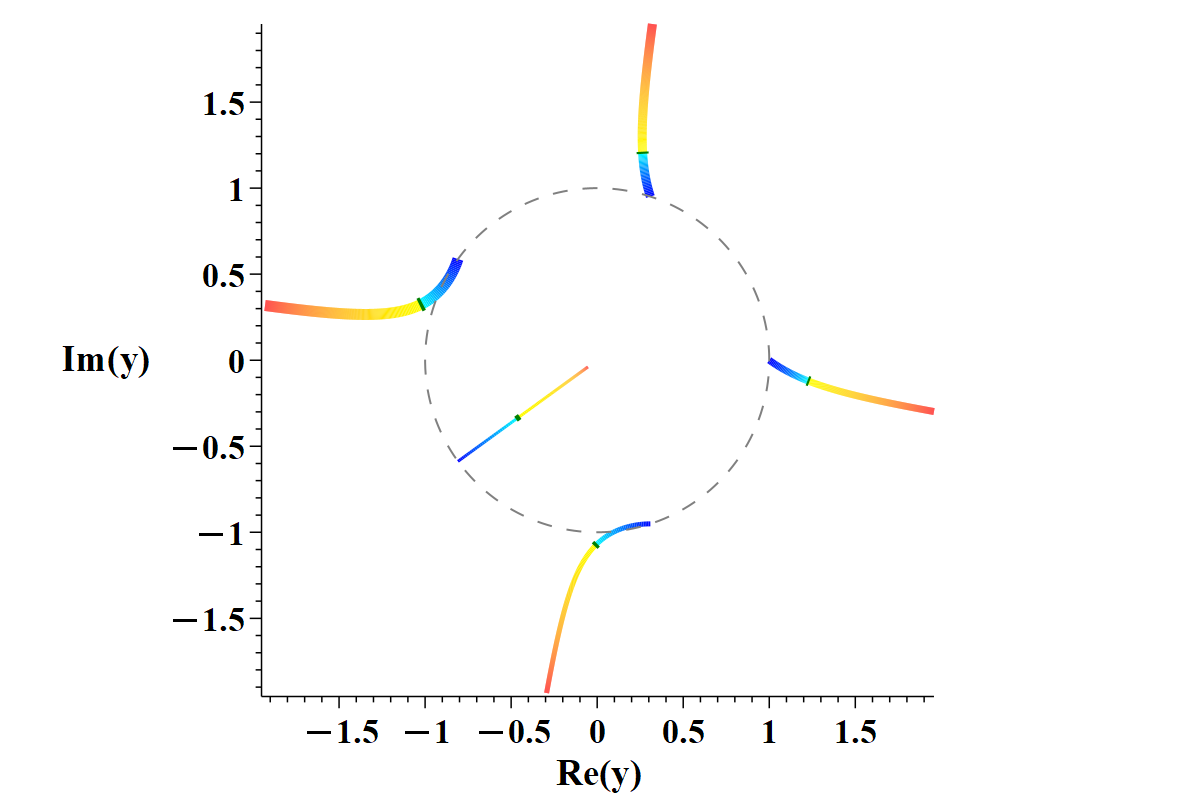}
        \textbf{(b)}
    \end{minipage}
    \caption{Phase hodographs for $a = 5$, $b = 1$ at angle $\arg(x) = -\pi/5$: 
    \textbf{(a)}~ultralogarithm ($m=0$, equation $e^{5w} = 1 + 5x e^{w}$), branches $n \in \{-5, \dots, 5\}$, $x_{\max} = 15$; 
    \textbf{(b)}~ultraradical ($m=1$, equation $y^{5} = 1 + 5xy$), branches $n \in \{-2, \dots, 2\}$, $x_{\max} = 3$.}
    \label{fig:5_1_minus_pi5}
\end{figure}

\begin{figure}[htbp]
    \centering
    \begin{minipage}[b]{0.48\textwidth}
        \centering
        \includegraphics[width=\textwidth]{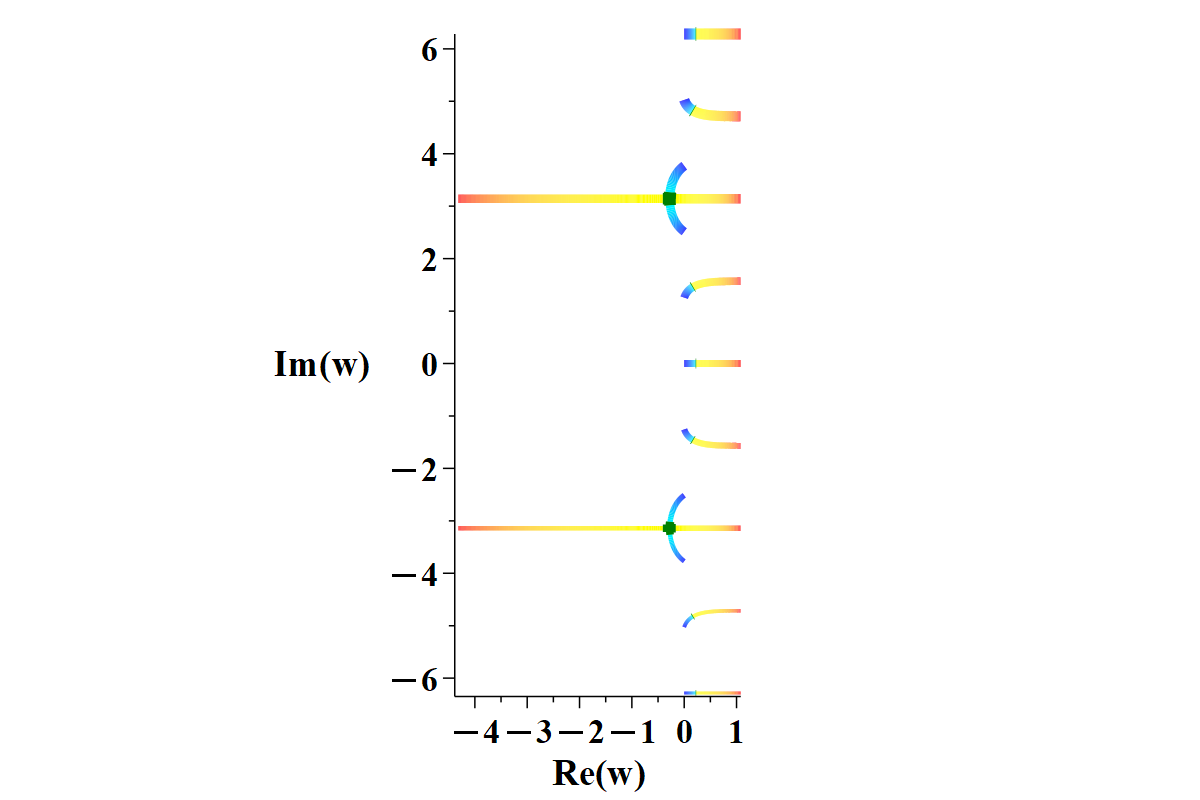}
        \textbf{(a)}
    \end{minipage}
    \hfill
    \begin{minipage}[b]{0.48\textwidth}
        \centering
        \includegraphics[width=\textwidth]{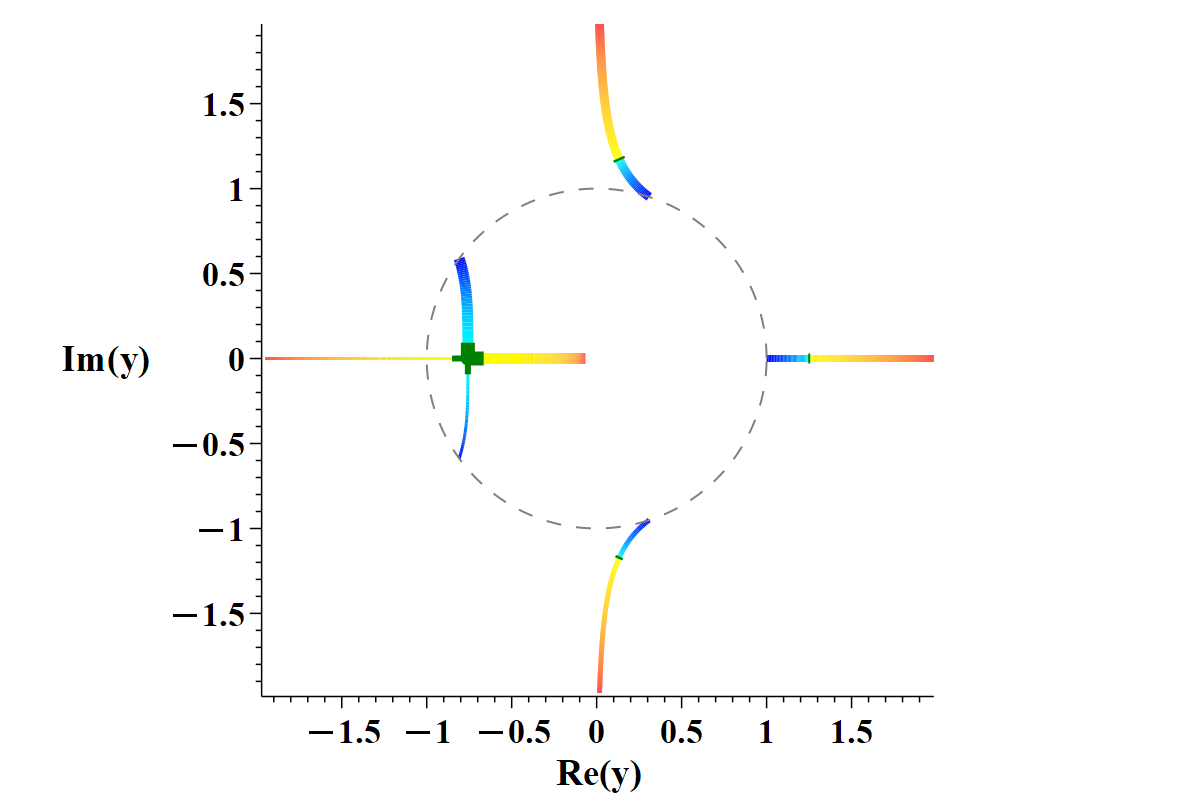}
        \textbf{(b)}
    \end{minipage}
    \caption{Phase hodographs for $a = 5$, $b = 1$ at angle $\arg(x) = 0$: 
    \textbf{(a)}~ultralogarithm ($m=0$, equation $e^{5w} = 1 + 5x e^{w}$), branches $n \in \{-5, \dots, 5\}$, $x_{\max} = 15$; 
    \textbf{(b)}~ultraradical ($m=1$, equation $y^{5} = 1 + 5xy$), branches $n \in \{-2, \dots, 2\}$, $x_{\max} = 3$.}
    \label{fig:5_1_zero}
\end{figure}

\begin{figure}[htbp]
    \centering
    \begin{minipage}[b]{0.48\textwidth}
        \centering
        \includegraphics[width=\textwidth]{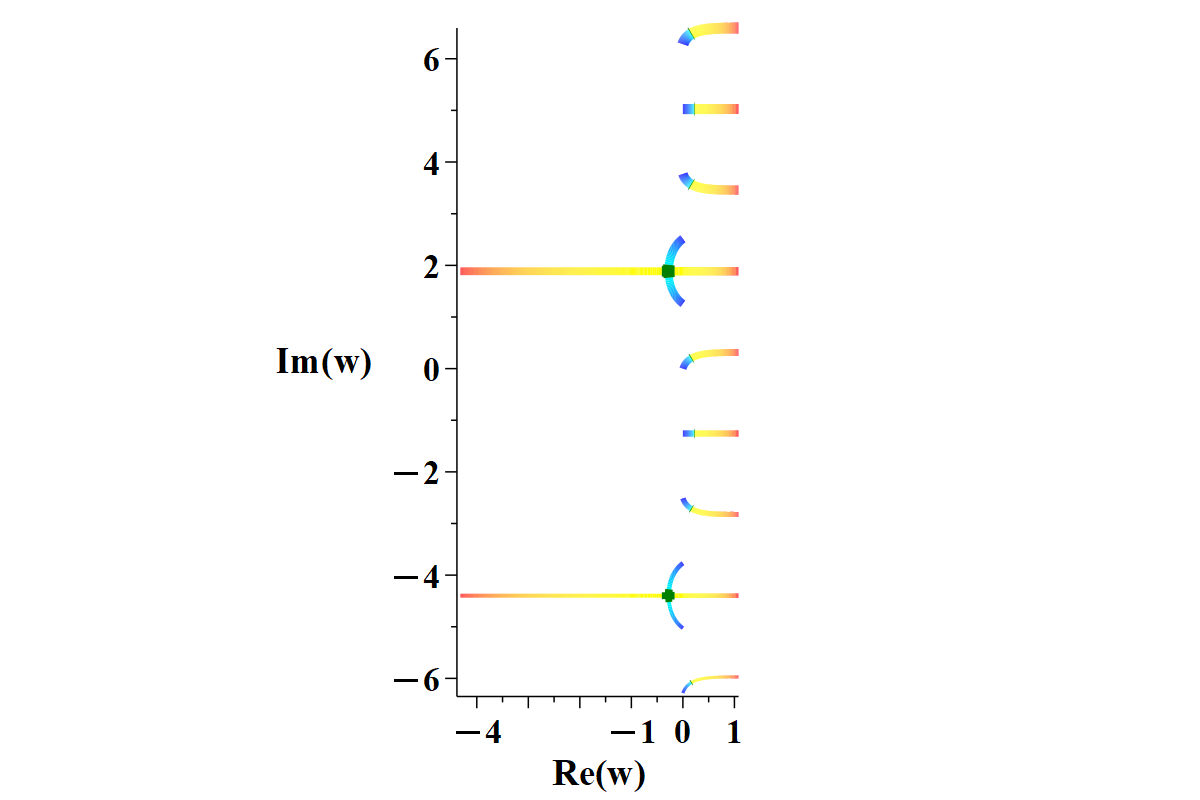}
        \textbf{(a)}
    \end{minipage}
    \hfill
    \begin{minipage}[b]{0.48\textwidth}
        \centering
        \includegraphics[width=\textwidth]{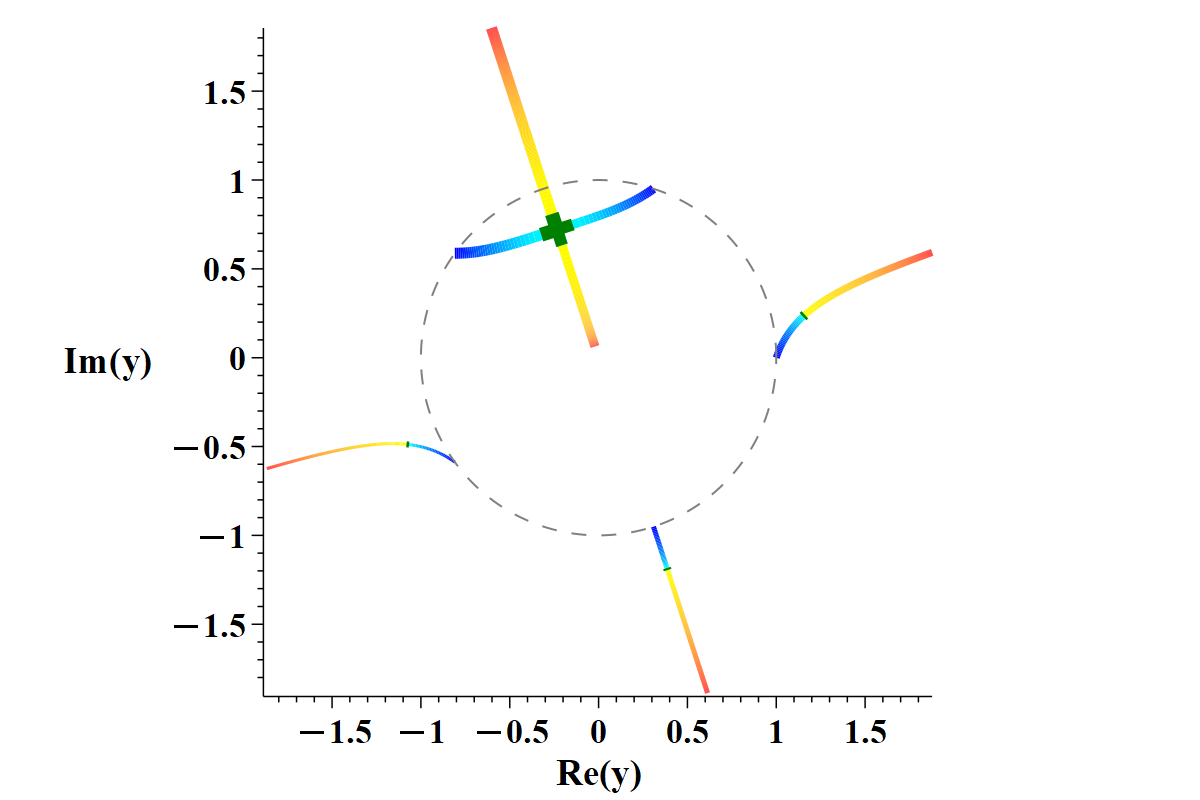}
        \textbf{(b)}
    \end{minipage}
    \caption{Phase hodographs for $a = 5$, $b = 1$ at angle $\arg(x) = 2\pi/5$: 
    \textbf{(a)}~ultralogarithm ($m=0$, equation $e^{5w} = 1 + 5x e^{w}$), branches $n \in \{-5, \dots, 5\}$, $x_{\max} = 15$; 
    \textbf{(b)}~ultraradical ($m=1$, equation $y^{5} = 1 + 5xy$), branches $n \in \{-2, \dots, 2\}$, $x_{\max} = 3$.}
    \label{fig:5_1_plus_2pi5}
\end{figure}

\begin{figure}[htbp]
    \centering
    \begin{minipage}[b]{0.48\textwidth}
        \centering
        \includegraphics[width=\textwidth]{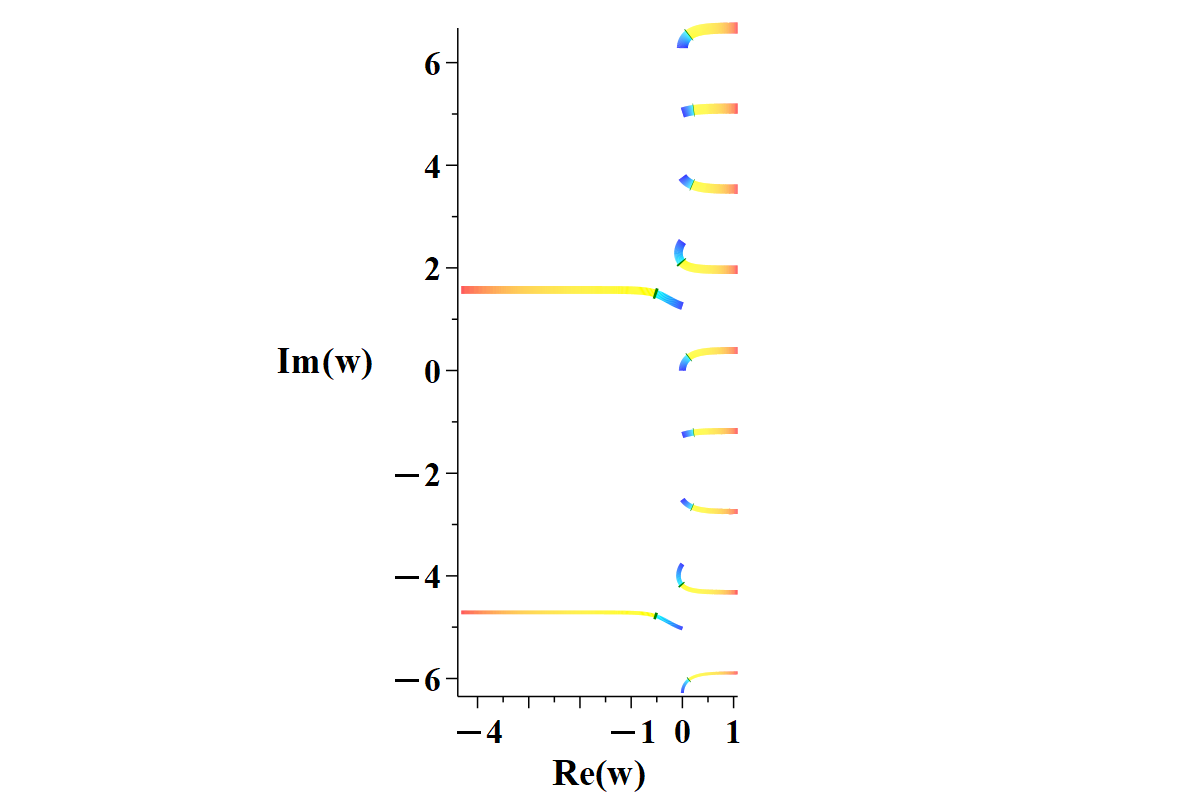}
        \textbf{(a)}
    \end{minipage}
    \hfill
    \begin{minipage}[b]{0.48\textwidth}
        \centering
        \includegraphics[width=\textwidth]{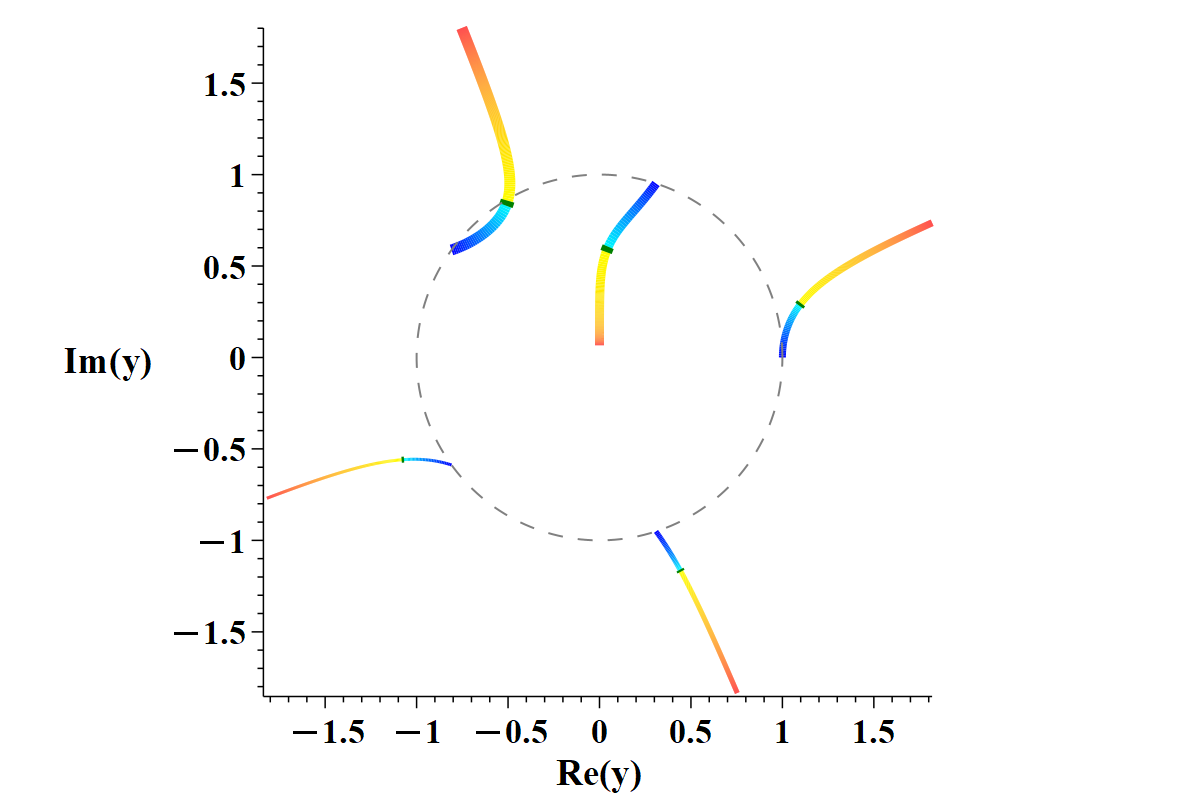}
        \textbf{(b)}
    \end{minipage}
    \caption{Phase hodographs for $a = 5$, $b = 1$ at angle $\arg(x) = \pi/2$: 
    \textbf{(a)}~ultralogarithm ($m=0$, equation $e^{5w} = 1 + 5x e^{w}$), branches $n \in \{-5, \dots, 5\}$, $x_{\max} = 15$; 
    \textbf{(b)}~ultraradical ($m=1$, equation $y^{5} = 1 + 5xy$), branches $n \in \{-2, \dots, 2\}$, $x_{\max} = 3$.}
    \label{fig:5_1_plus_pi2}
\end{figure}

To trace the periodic junction structure of independent and dependent branches across period boundaries, Figs.~\ref{fig:5_1_minus_pi}--\ref{fig:5_1_plus_pi2} display the paired phase hodographs for an extended set of $11$ ultralogarithm branches ($n \in \{-5, \dots, 5\}$, $x_{\max} = 15$) and $5$ independent ultraradical branches ($n \in \{-2, \dots, 2\}$, $x_{\max} = 3$) along seven characteristic rays $\arg(x) \in \{-\pi, -4\pi/5, -2\pi/5, -\pi/5, 0, 2\pi/5, \pi/2\}$. Analyzing all seven pairs reveals three general patterns:
\begin{itemize}
\item \textbf{Periodic lattice of branch points along critical rays.} Along all four critical rays $\arg(x) = -4\pi/5$, $-2\pi/5$, $0$, and $+2\pi/5$ (Figs.~\ref{fig:5_1_minus_4pi5}, \ref{fig:5_1_minus_2pi5}, \ref{fig:5_1_zero}, and~\ref{fig:5_1_plus_2pi5}), two branches in the ultraradical plane (panels \textbf{b}) collide \emph{inside} the unit circle at $\vert{}y\vert{} \approx 0.758$ (green cross). On the ultralogarithm hodographs (panels \textbf{a}), this collision repeats periodically with step $\Delta n = a = 5$ (vertical spacing $2\pi$): at $\arg(x) = -4\pi/5$, the pairs $(n=-1, 0)$ and $(n=4, 5)$ touch at levels $-\pi/5 + 2\pi p$; at $\arg(x) = -2\pi/5$, the pairs $(n=-2, -1)$ and $(n=3, 4)$ touch at $-3\pi/5 + 2\pi p$; at $\arg(x) = 0$, the pairs $(n=-3, -2)$ and $(n=2, 3)$ meet at the period boundaries $\pm \pi + 2\pi p$; and at $\arg(x) = +2\pi/5$, the pairs $(n=-4, -3)$ and $(n=1, 2)$ touch at $+3\pi/5 + 2\pi p$.
\item \textbf{Asymptotic ``$1$-left vs. $(a-1)$-right'' splitting.} Within each period of five branches as $\vert{}x\vert{} \to \infty$, strictly \emph{one} branch recedes into the left half-plane ($\mathrm{Re}(w) \to -\infty$, corresponding to the root $y \sim -(ax)^{-1} \to 0$ inside the unit circle and governed by the conjugate series~\eqref{eq:ulog_outer_k} with index $n_{-}$), whereas the remaining $a - 1 = 4$ branches recede into the right half-plane ($\mathrm{Re}(w) \to +\infty$, corresponding to the roots $y \sim (ax)^{1/4} \to \infty$ and the series~\eqref{eq:ulog_outer_h} with indices $n_{+}$). Along critical rays, the branch point acts as the bifurcation point between these two asymptotic families.
\item \textbf{Smoothness of trajectories along non-critical rays.} For $\arg(x) = -\pi$, $-\pi/5$, and $\pi/2$ (Figs.~\ref{fig:5_1_minus_pi}, \ref{fig:5_1_minus_pi5}, and~\ref{fig:5_1_plus_pi2}), no branch encounters a branch point: all $11$ ultralogarithm trajectories form smooth, non-intersecting curves confined within their respective phase strips of width $\frac{2\pi}{5}$.
\end{itemize}

\subsection{Septic Equation with Exponent: a = 7, b = 3}
\label{subsec:a7_b3}

To test the sector-based conjugation criterion for a non-trivial intermediate exponent $b > 1$, we consider the septic equation with coprime integer parameters $a = 7$ and $b = 3$:
\begin{equation}
\label{eq:a7_b3}
e^{7w} = 1 + 7x e^{3w},
\end{equation}
which corresponds to the algebraic ultraradical $y^7 - 7xy^3 - 1 = 0$. Here, the radius of convergence of the basic series~\eqref{eq:radius_R} is $R = \frac{(4/3)^{3/7}}{4} \approx 0.282811$. In the complex $w$-plane itself, the actual width of each branch's phase strip (the spacing between the base levels $\mathrm{Im}\,f^{0}_{n} = \frac{2\pi n}{7}$) equals $\frac{2\pi}{7}$, so every seven branches fit into a layer of height $2\pi \approx 6.28$; upon passing to the scaled phase indicators $L^{0}_{n}, L^{\pm}_{n} = \mathrm{Im}(3 f^{0,\pm}_{n})$, where all reference values are multiplied by $b = 3$ to align the periods of fractional and transformed exponents, the scaled width of the sector $s_n$ equals $\frac{6\pi}{7}$. The branch-point conditions $y^7 = -3/4$ and $x = -1/(4y^3)$ yield the root modulus $\vert{}y\vert{} = (3/4)^{1/7} \approx 0.960 < 1$ ($\mathrm{Re}(w) \approx -0.041$) and lead to $x^7 = (27 \cdot 4^4)^{-1} > 0$, showing that the critical rays in the $x$-plane are determined by even multiples of the angle $\frac{\pi}{7}$: $\arg(x) \in \{0, \pm 2\pi/7, \pm 4\pi/7, \pm 6\pi/7\}$. By Corollary~1, the equation possesses $7$ independent branches ($n \in \{-3, \dots, 3\}$) satisfying the sum invariant $\sum_{n=-3}^{3} w_n \equiv 6\pi i \pmod{2\pi i}$.

\begin{figure}[htbp]
    \centering
    \begin{minipage}[b]{0.48\textwidth}
        \centering
        \includegraphics[width=\textwidth]{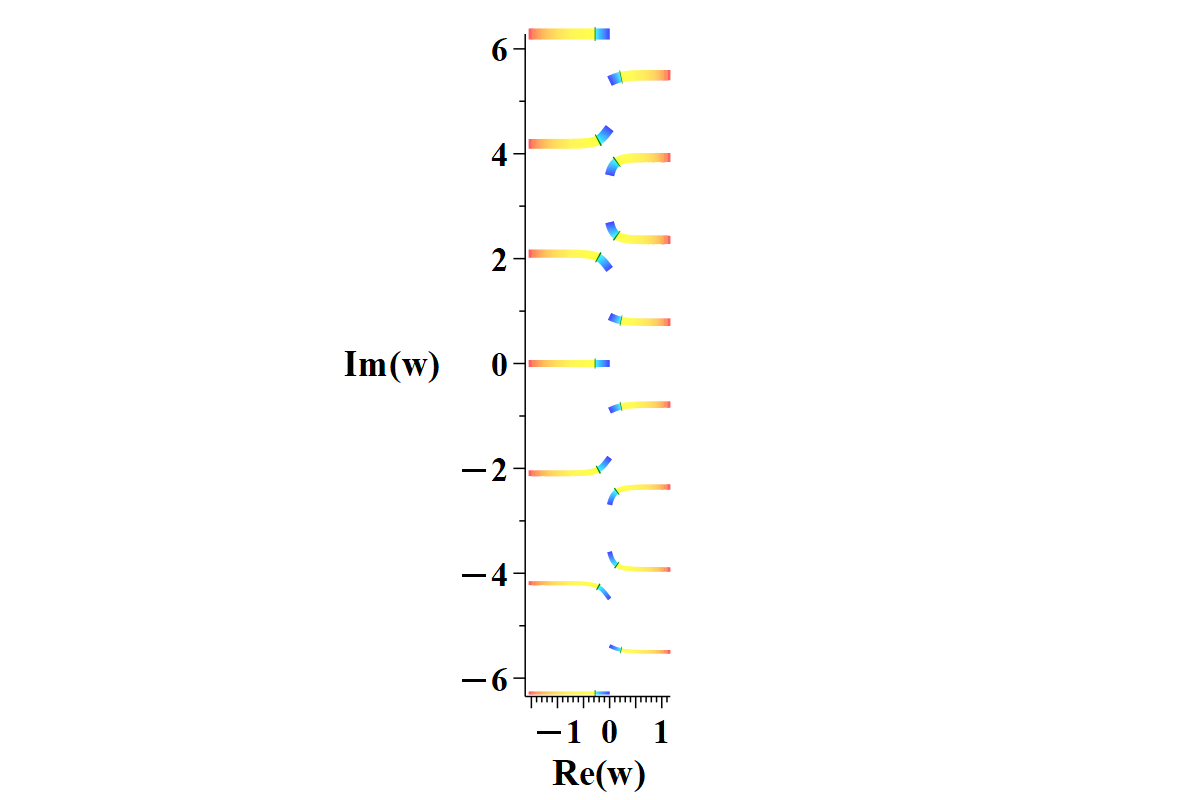}
        \textbf{(a)}
    \end{minipage}
    \hfill
    \begin{minipage}[b]{0.48\textwidth}
        \centering
        \includegraphics[width=\textwidth]{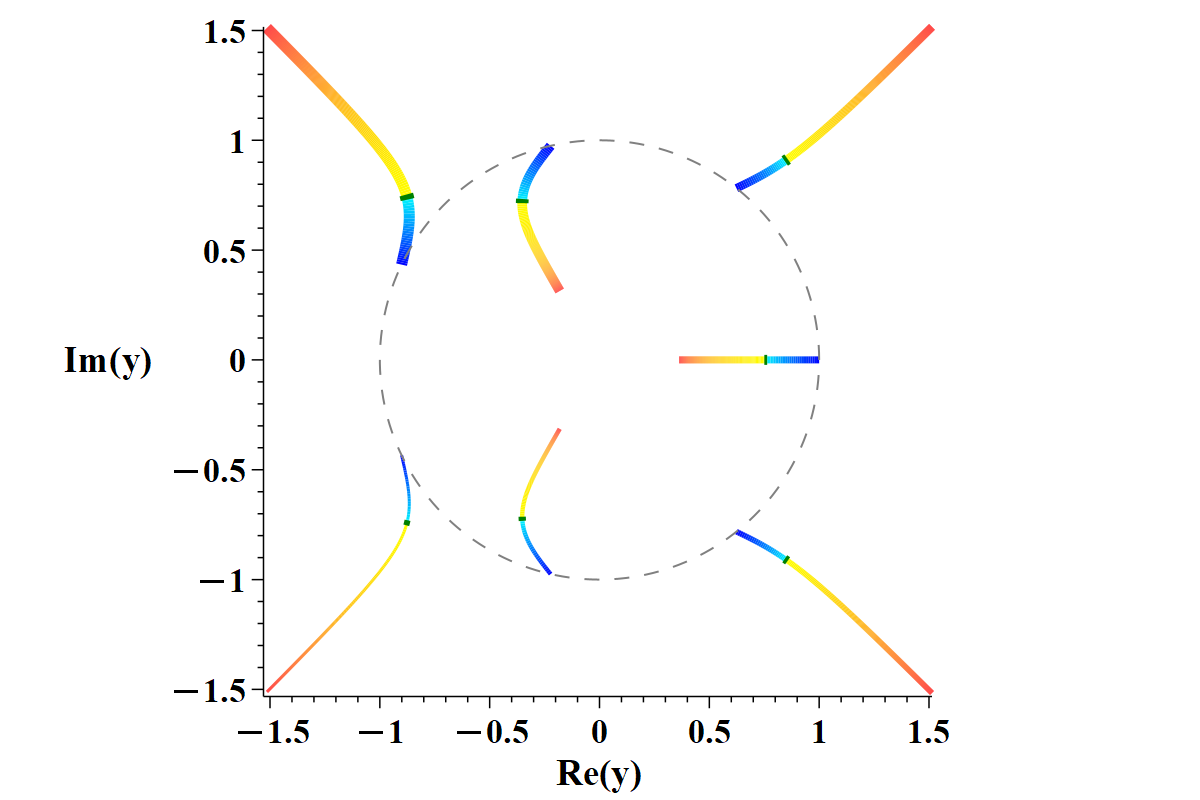}
        \textbf{(b)}
    \end{minipage}
    \caption{Phase hodographs for $a = 7$, $b = 3$ at angle $\arg(x) = -\pi$: 
    \textbf{(a)}~ultralogarithm ($m=0$, equation $e^{7w} = 1 + 7x e^{3w}$), branches $n \in \{-5, \dots, 5\}$, $x_{\max} = 15$; 
    \textbf{(b)}~ultraradical ($m=1$, equation $y^{7} = 1 + 7xy^{3}$), branches $n \in \{-3, \dots, 3\}$, $x_{\max} = 3$.}
    \label{fig:7_3_minus_pi}
\end{figure}

\begin{figure}[htbp]
    \centering
    \begin{minipage}[b]{0.48\textwidth}
        \centering
        \includegraphics[width=\textwidth]{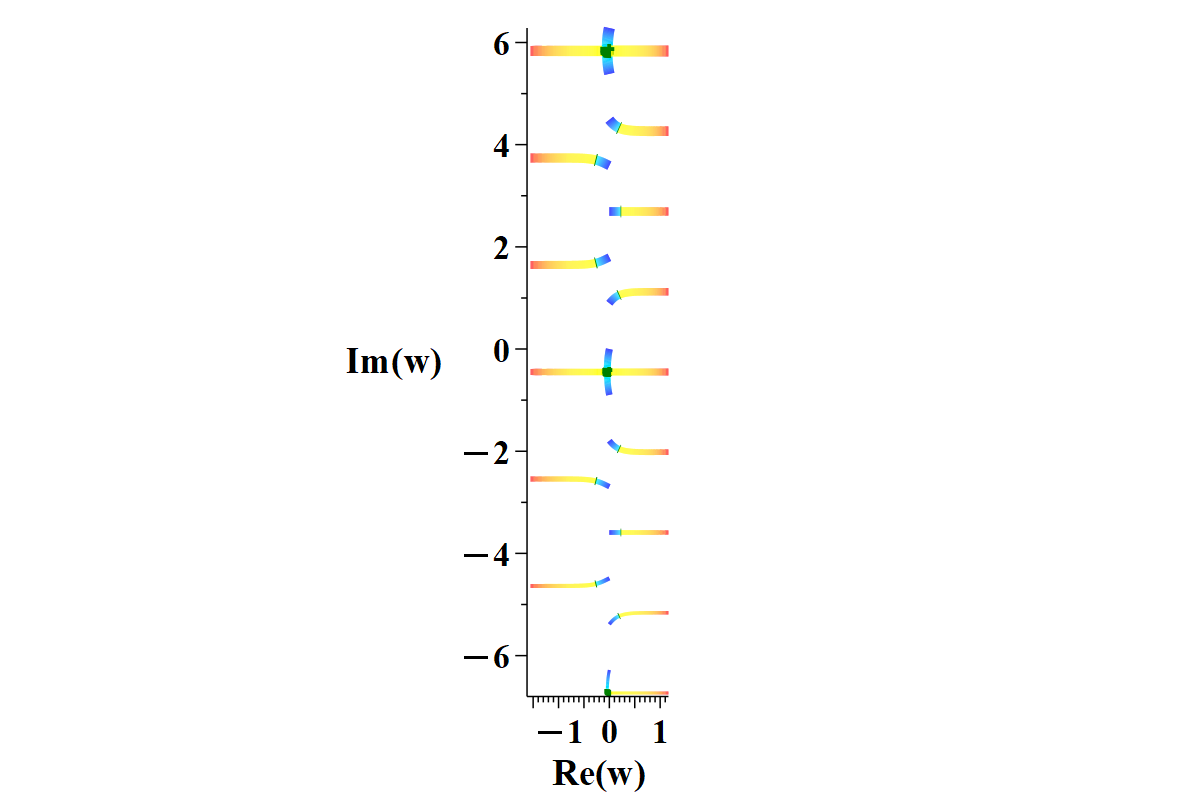}
        \textbf{(a)}
    \end{minipage}
    \hfill
    \begin{minipage}[b]{0.48\textwidth}
        \centering
        \includegraphics[width=\textwidth]{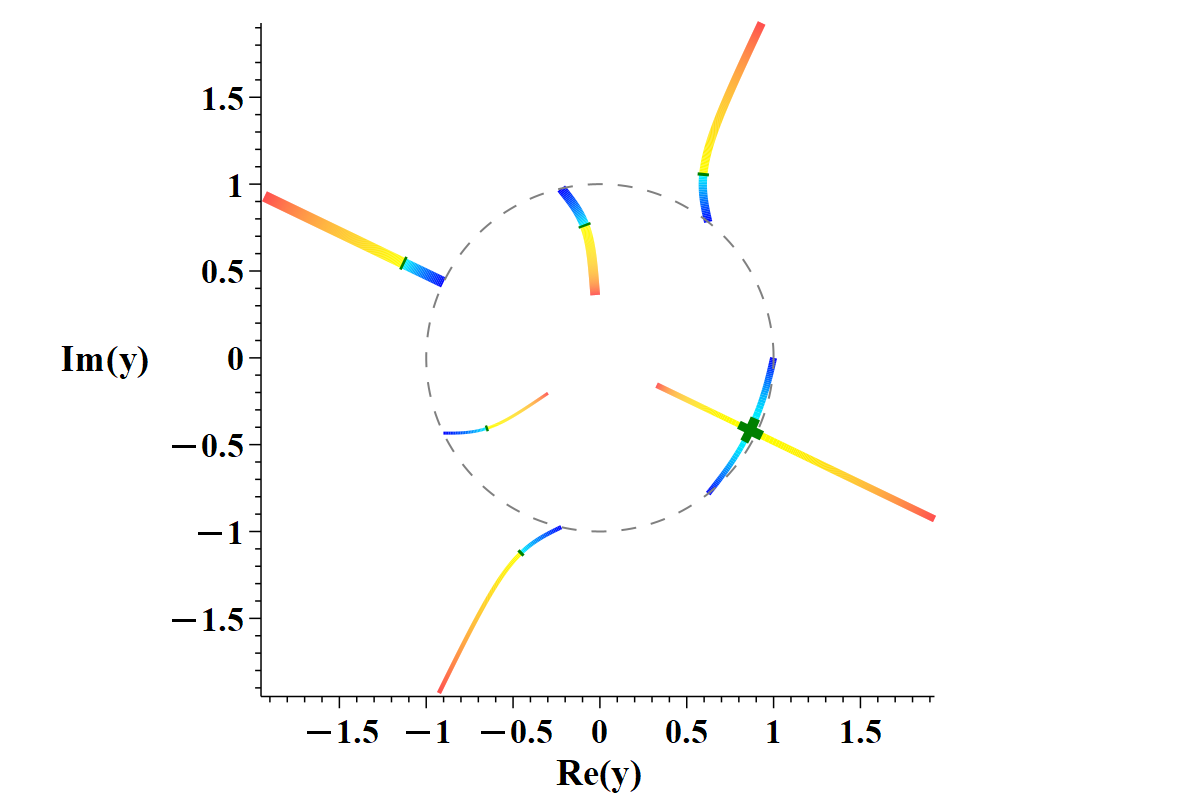}
        \textbf{(b)}
    \end{minipage}
    \caption{Phase hodographs for $a = 7$, $b = 3$ at angle $\arg(x) = -4\pi/7$: 
    \textbf{(a)}~ultralogarithm ($m=0$, equation $e^{7w} = 1 + 7x e^{3w}$), branches $n \in \{-5, \dots, 5\}$, $x_{\max} = 15$; 
    \textbf{(b)}~ultraradical ($m=1$, equation $y^{7} = 1 + 7xy^{3}$), branches $n \in \{-3, \dots, 3\}$, $x_{\max} = 3$.}
    \label{fig:7_3_minus_4pi7}
\end{figure}

\begin{figure}[htbp]
    \centering
    \begin{minipage}[b]{0.48\textwidth}
        \centering
        \includegraphics[width=\textwidth]{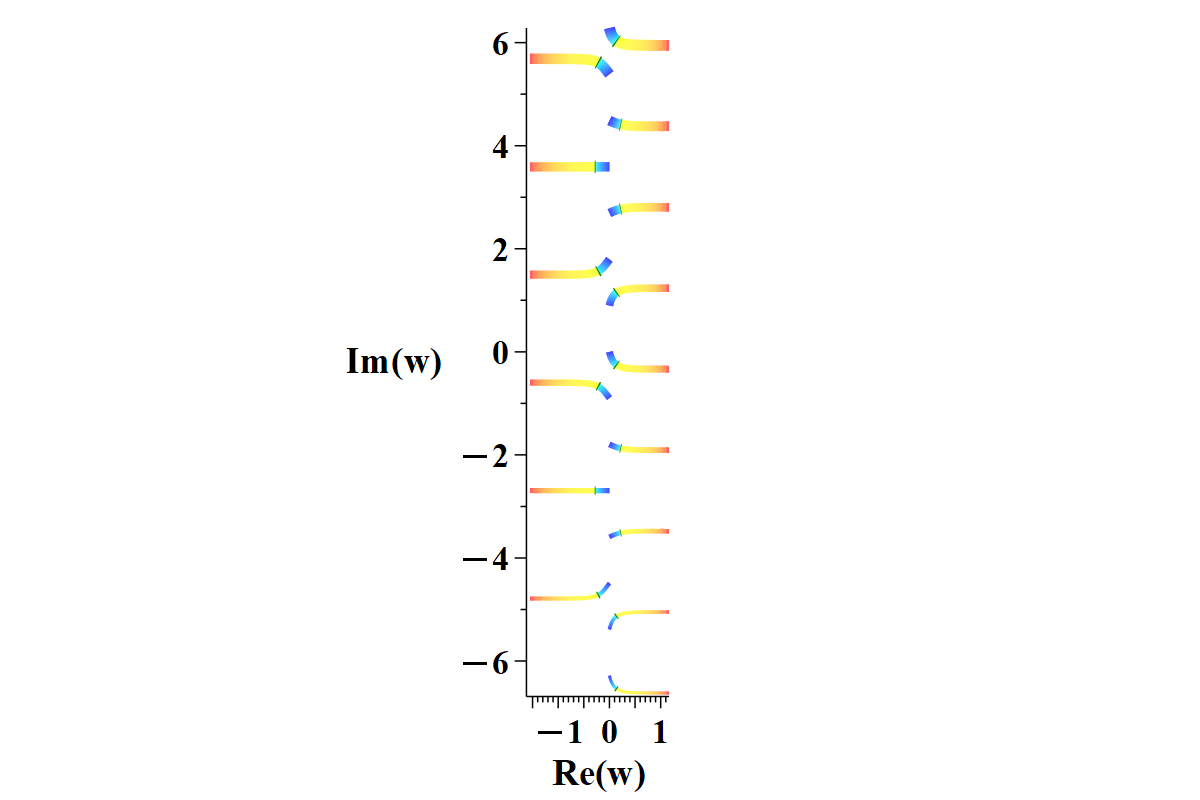}
        \textbf{(a)}
    \end{minipage}
    \hfill
    \begin{minipage}[b]{0.48\textwidth}
        \centering
        \includegraphics[width=\textwidth]{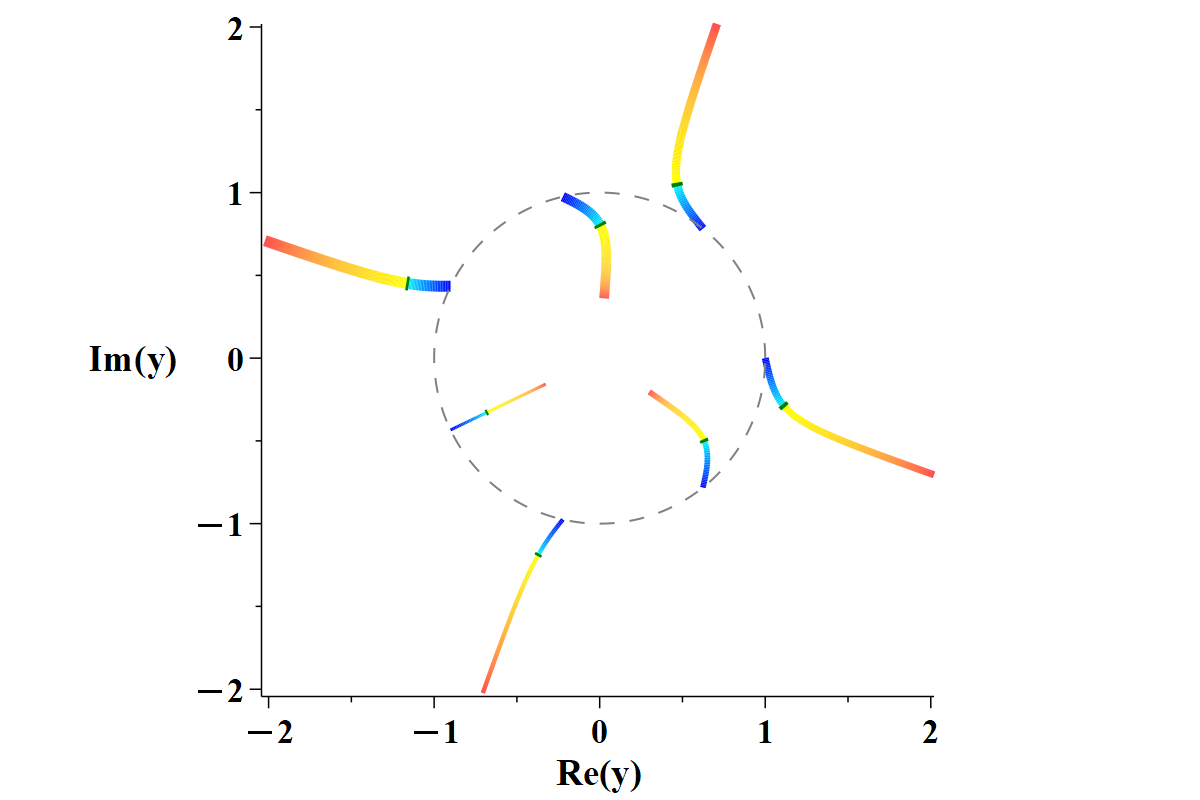}
        \textbf{(b)}
    \end{minipage}
    \caption{Phase hodographs for $a = 7$, $b = 3$ at angle $\arg(x) = -3\pi/7$: 
    \textbf{(a)}~ultralogarithm ($m=0$, equation $e^{7w} = 1 + 7x e^{3w}$), branches $n \in \{-5, \dots, 5\}$, $x_{\max} = 15$; 
    \textbf{(b)}~ultraradical ($m=1$, equation $y^{7} = 1 + 7xy^{3}$), branches $n \in \{-3, \dots, 3\}$, $x_{\max} = 3$.}
    \label{fig:7_3_minus_3pi7}
\end{figure}

\begin{figure}[htbp]
    \centering
    \begin{minipage}[b]{0.48\textwidth}
        \centering
        \includegraphics[width=\textwidth]{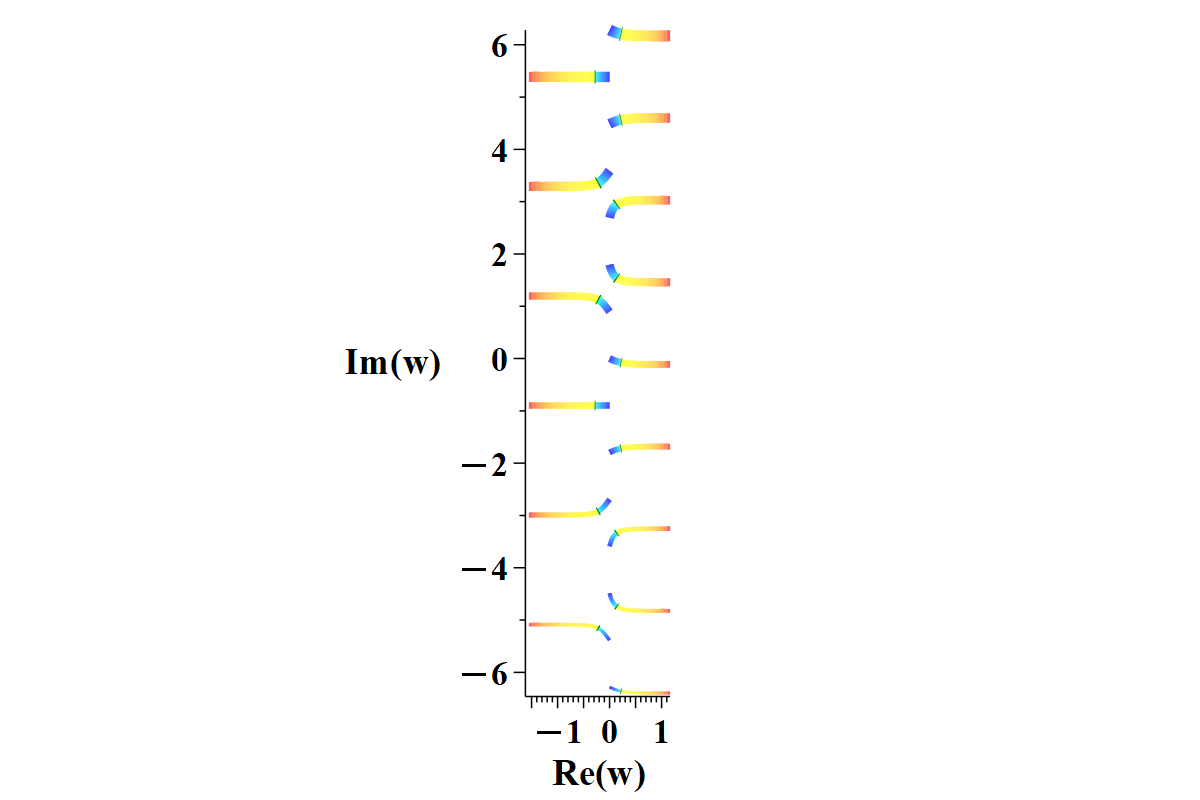}
        \textbf{(a)}
    \end{minipage}
    \hfill
    \begin{minipage}[b]{0.48\textwidth}
        \centering
        \includegraphics[width=\textwidth]{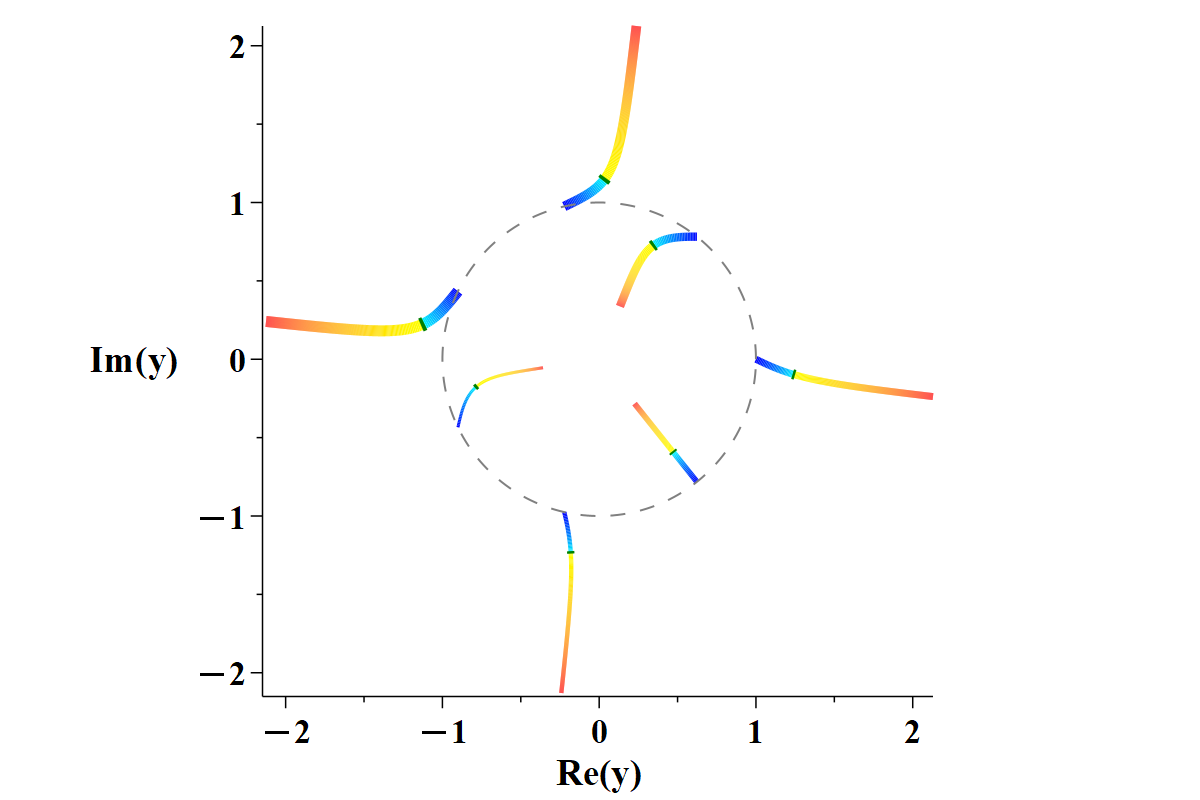}
        \textbf{(b)}
    \end{minipage}
    \caption{Phase hodographs for $a = 7$, $b = 3$ at angle $\arg(x) = -\pi/7$: 
    \textbf{(a)}~ultralogarithm ($m=0$, equation $e^{7w} = 1 + 7x e^{3w}$), branches $n \in \{-5, \dots, 5\}$, $x_{\max} = 15$; 
    \textbf{(b)}~ultraradical ($m=1$, equation $y^{7} = 1 + 7xy^{3}$), branches $n \in \{-3, \dots, 3\}$, $x_{\max} = 3$.}
    \label{fig:7_3_minus_pi7}
\end{figure}

\begin{figure}[htbp]
    \centering
    \begin{minipage}[b]{0.48\textwidth}
        \centering
        \includegraphics[width=\textwidth]{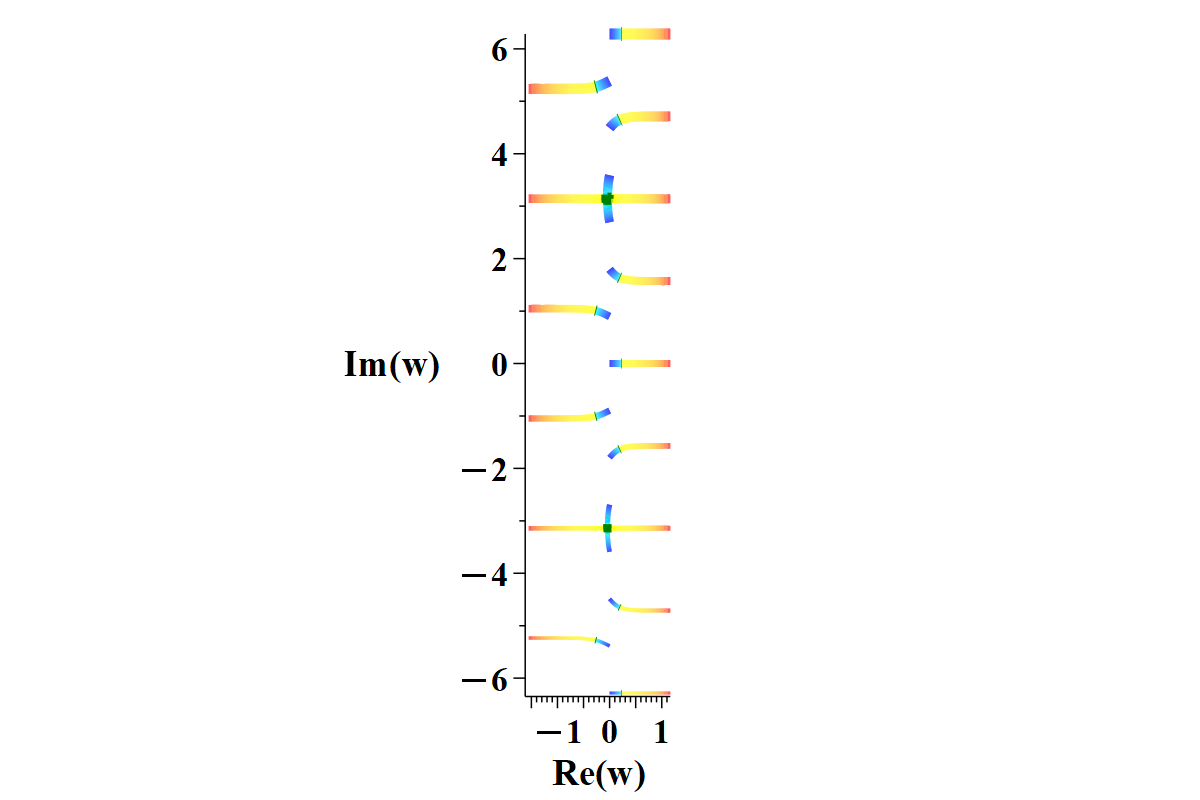}
        \textbf{(a)}
    \end{minipage}
    \hfill
    \begin{minipage}[b]{0.48\textwidth}
        \centering
        \includegraphics[width=\textwidth]{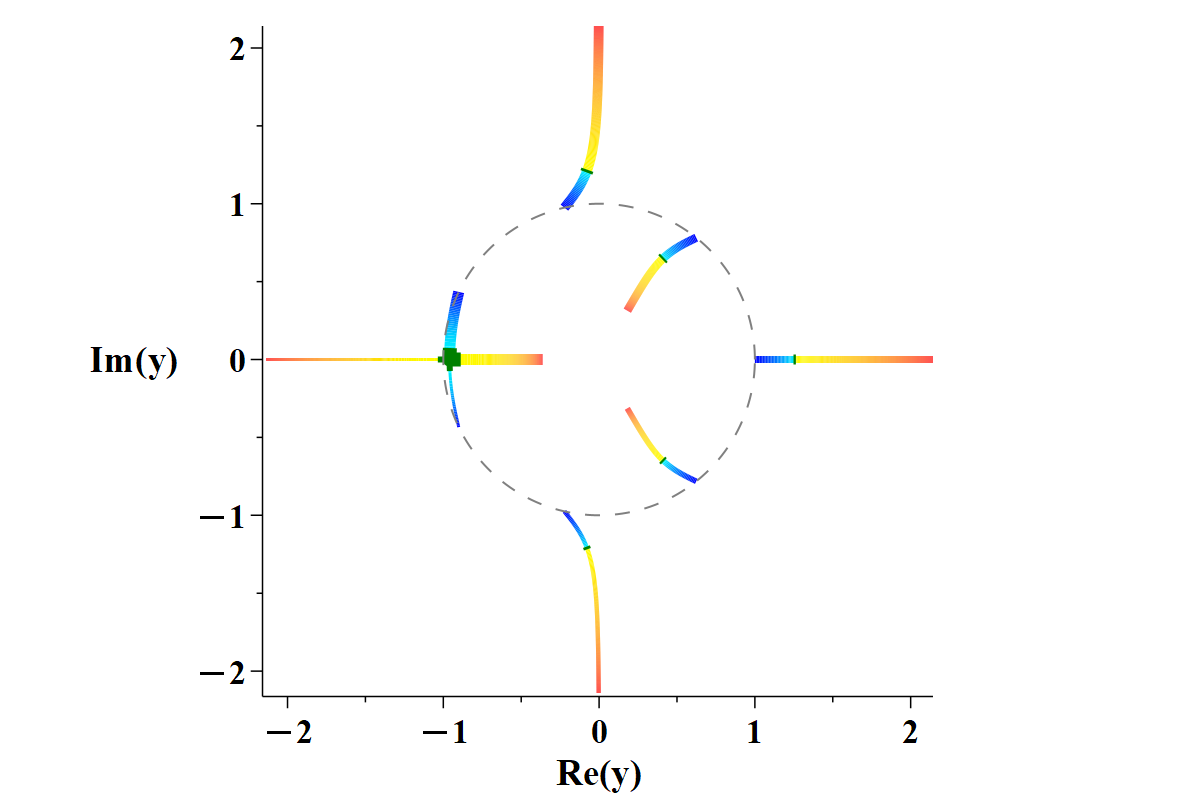}
        \textbf{(b)}
    \end{minipage}
    \caption{Phase hodographs for $a = 7$, $b = 3$ at angle $\arg(x) = 0$: 
    \textbf{(a)}~ultralogarithm ($m=0$, equation $e^{7w} = 1 + 7x e^{3w}$), branches $n \in \{-5, \dots, 5\}$, $x_{\max} = 15$; 
    \textbf{(b)}~ultraradical ($m=1$, equation $y^{7} = 1 + 7xy^{3}$), branches $n \in \{-3, \dots, 3\}$, $x_{\max} = 3$.}
    \label{fig:7_3_zero}
\end{figure}

\begin{figure}[htbp]
    \centering
    \begin{minipage}[b]{0.48\textwidth}
        \centering
        \includegraphics[width=\textwidth]{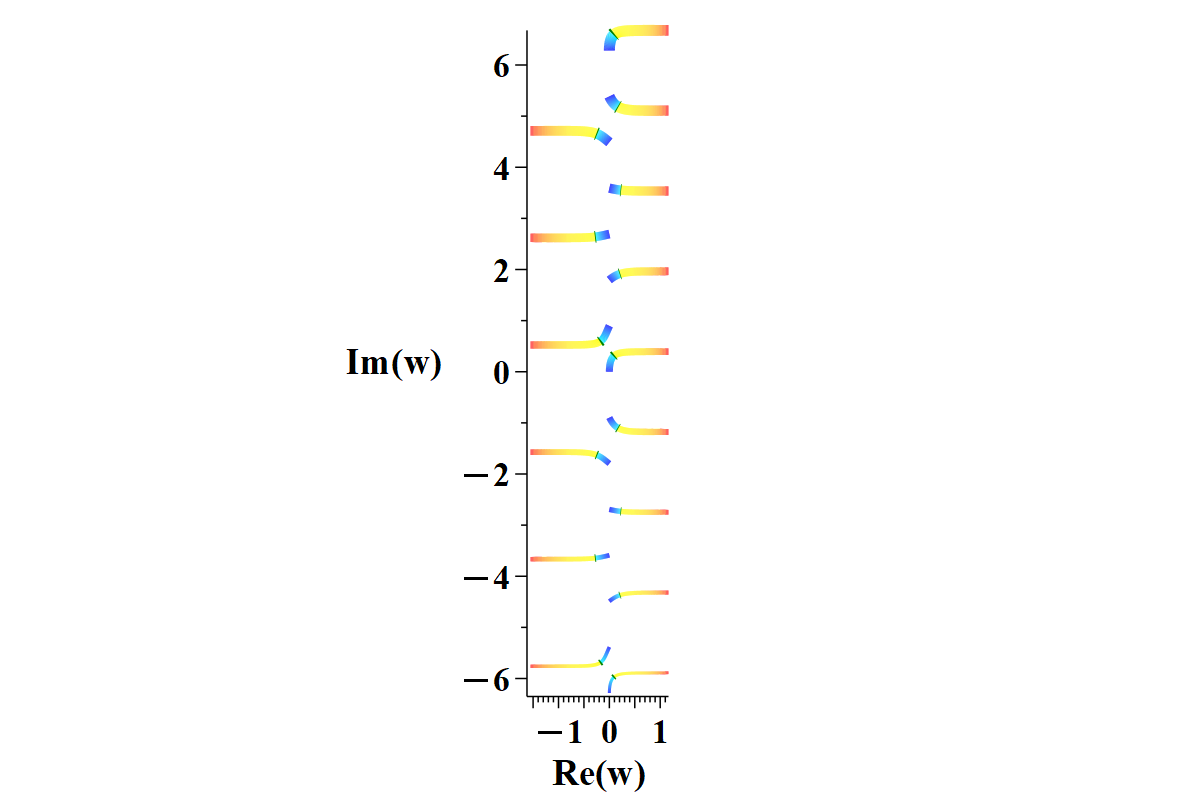}
        \textbf{(a)}
    \end{minipage}
    \hfill
    \begin{minipage}[b]{0.48\textwidth}
        \centering
        \includegraphics[width=\textwidth]{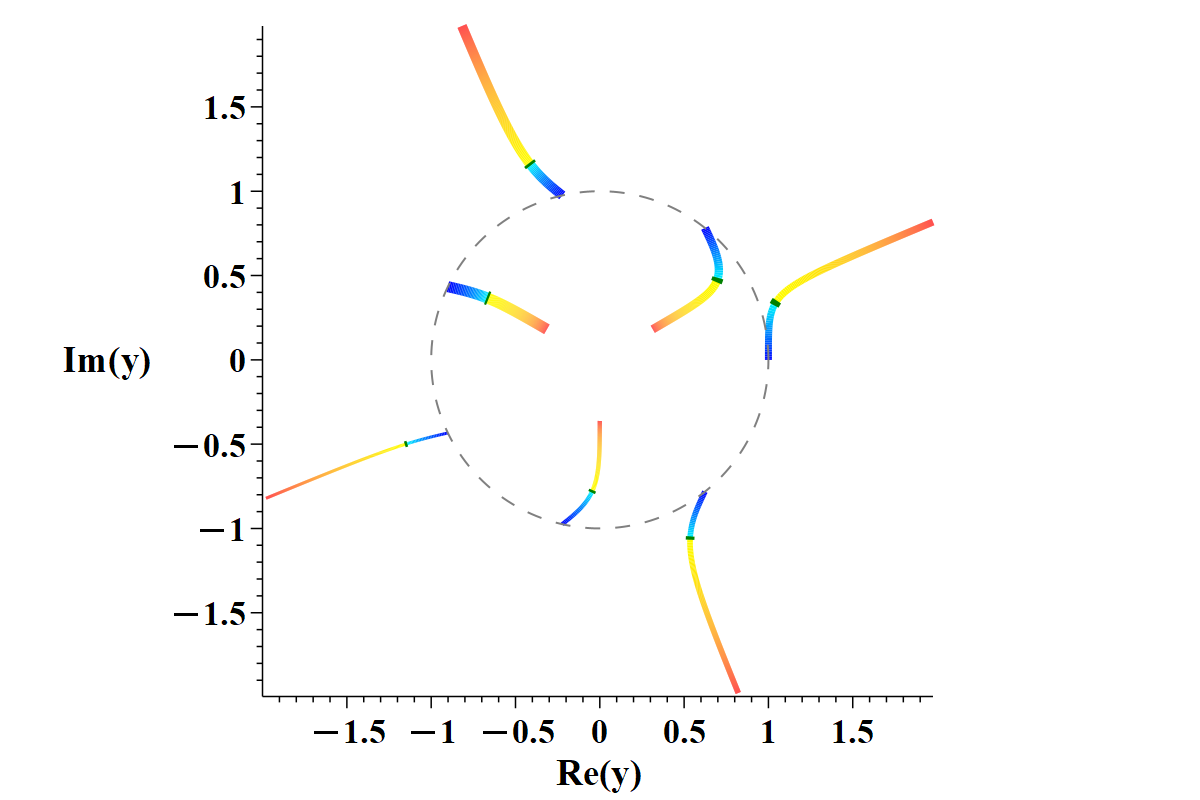}
        \textbf{(b)}
    \end{minipage}
    \caption{Phase hodographs for $a = 7$, $b = 3$ at angle $\arg(x) = \pi/2$: 
    \textbf{(a)}~ultralogarithm ($m=0$, equation $e^{7w} = 1 + 7x e^{3w}$), branches $n \in \{-5, \dots, 5\}$, $x_{\max} = 15$; 
    \textbf{(b)}~ultraradical ($m=1$, equation $y^{7} = 1 + 7xy^{3}$), branches $n \in \{-3, \dots, 3\}$, $x_{\max} = 3$.}
    \label{fig:7_3_plus_pi2}
\end{figure}

\begin{figure}[htbp]
    \centering
    \begin{minipage}[b]{0.48\textwidth}
        \centering
        \includegraphics[width=\textwidth]{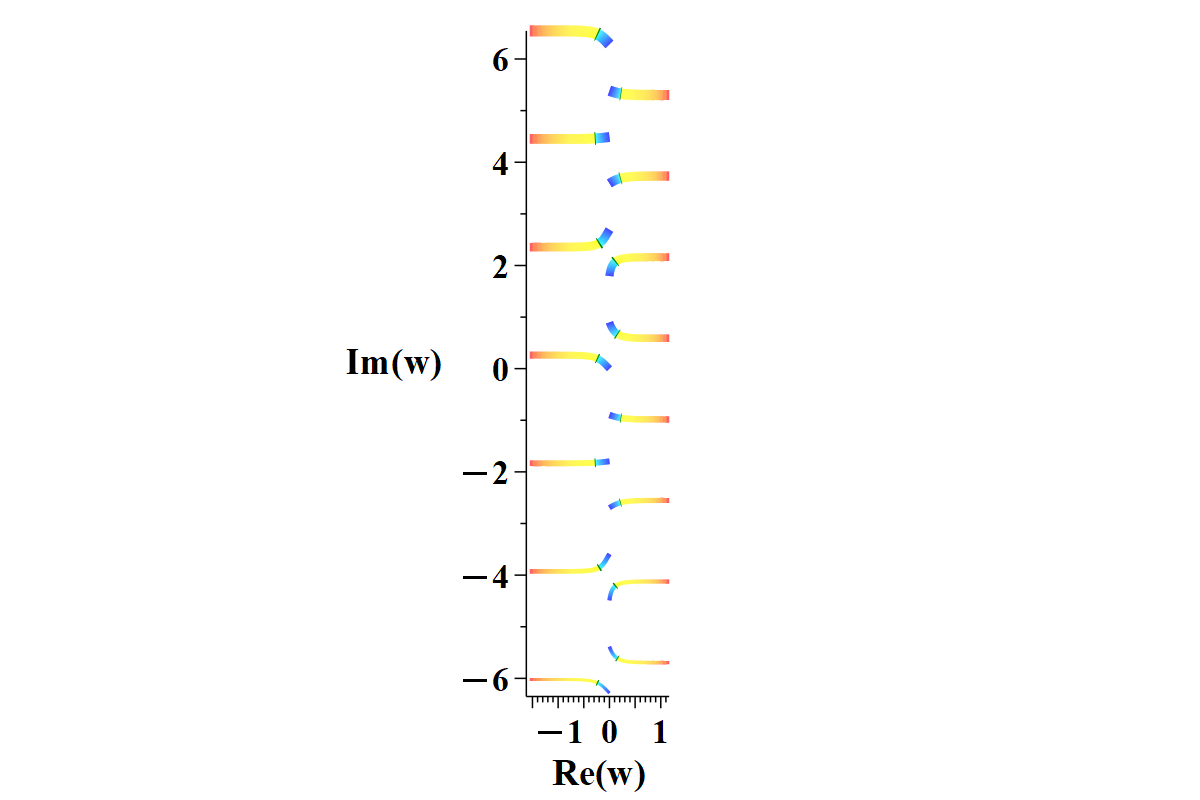}
        \textbf{(a)}
    \end{minipage}
    \hfill
    \begin{minipage}[b]{0.48\textwidth}
        \centering
        \includegraphics[width=\textwidth]{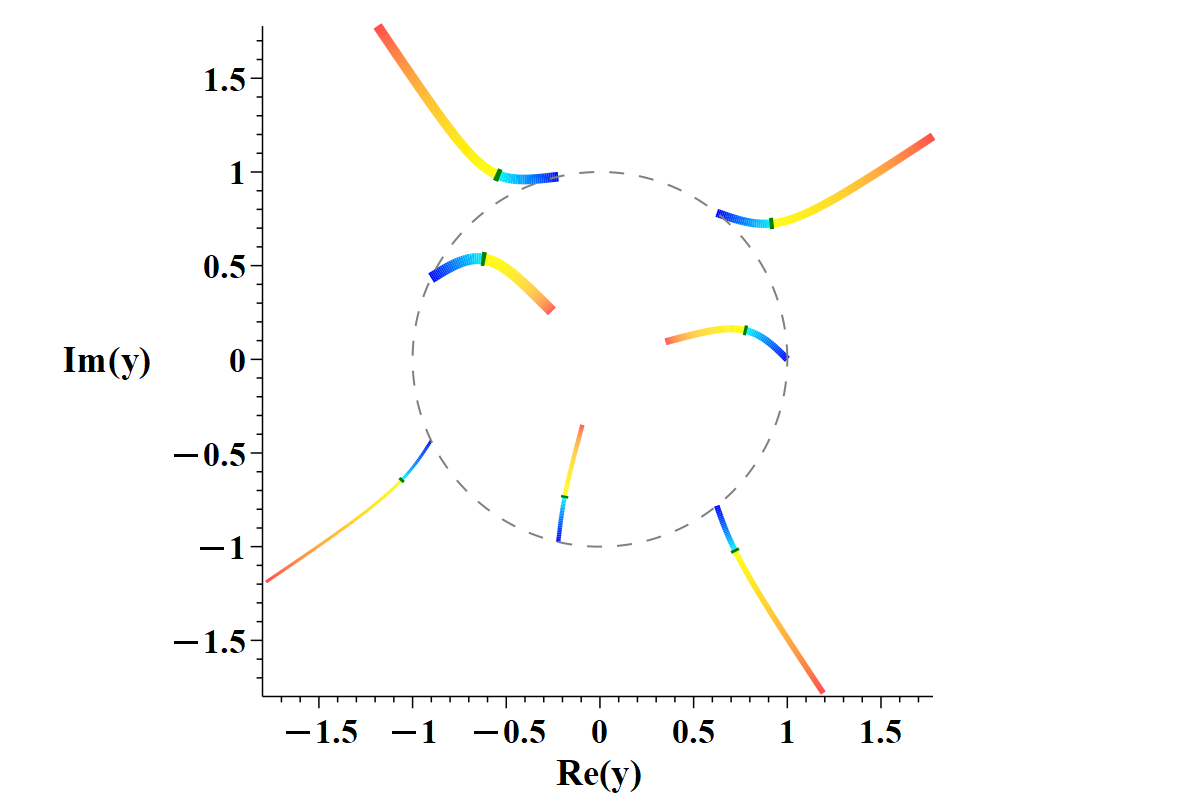}
        \textbf{(b)}
    \end{minipage}
    \caption{Phase hodographs for $a = 7$, $b = 3$ at angle $\arg(x) = 3\pi/4$: 
    \textbf{(a)}~ultralogarithm ($m=0$, equation $e^{7w} = 1 + 7x e^{3w}$), branches $n \in \{-5, \dots, 5\}$, $x_{\max} = 15$; 
    \textbf{(b)}~ultraradical ($m=1$, equation $y^{7} = 1 + 7xy^{3}$), branches $n \in \{-3, \dots, 3\}$, $x_{\max} = 3$.}
    \label{fig:7_3_plus_3pi4}
\end{figure}

Figures~\ref{fig:7_3_minus_pi}--\ref{fig:7_3_plus_3pi4} compare the paired phase hodographs of $11$ ultralogarithm branches ($n \in \{-5, \dots, 5\}$, $x_{\max} = 15$) and all $7$ independent ultraradical branches ($n \in \{-3, \dots, 3\}$, $x_{\max} = 3$) along seven rays $\arg(x) \in \{-\pi, -4\pi/7, -3\pi/7, -\pi/7, 0, \pi/2, 3\pi/4\}$. The resulting portraits clearly illustrate three key properties:
\begin{itemize}
\item \textbf{Asymmetric ``$b = 3$-left vs. $(a - b) = 4$-right'' splitting.} Unlike the cases with $b = 1$, where only a single branch per period receded into the left half-plane, for $b = 3$ within each period of seven ultralogarithm branches (panels \textbf{a}) exactly \emph{three} branches recede to the left ($\mathrm{Re}(w) \to -\infty$, conjugate series~\eqref{eq:ulog_outer_k} with index $n_{-}$), whereas \emph{four} branches recede to the right ($\mathrm{Re}(w) \to +\infty$, conjugate series~\eqref{eq:ulog_outer_h} with indices $n_{+}$).
\item \textbf{Visualization of Newton's polygon law.} On the ultraradical phase hodographs (panels \textbf{b}), one can directly observe the asymptotic rearrangement of roots governed by the Newton polygon of $y^7 - 7xy^3 - 1 = 0$: starting at $\vert{}x\vert{} = 0$ from the vertices of a \emph{regular heptagon} on the unit circle, the trajectories split into two regular polygons as $\vert{}x\vert{}$ grows. Four outer branches ($y \sim (7x)^{1/4}$) align toward the vertices of an \emph{outer square}, while three inner branches ($y \sim -(7x)^{-1/3}$) curve inward inside the unit circle and converge toward the origin along the rays of an \emph{equilateral triangle}, forming characteristic symmetric arc configurations (particularly prominent in Figs.~\ref{fig:7_3_plus_pi2}b and~\ref{fig:7_3_plus_3pi4}b). Moreover, the conjugate series formulas~\eqref{eq:ulog_outer_h} and~\eqref{eq:ulog_outer_k} imply counter-rotation of these two polygons as the ray angle $\arg(x)$ varies: the outer square rotates counterclockwise at an angular rate of $+\frac{1}{a-b} = +\frac{1}{4}$, whereas the inner triangle rotates clockwise at a rate of $-\frac{1}{b} = -\frac{1}{3}$.
\item \textbf{Critical angle lattice and continuation stability.} In general, the branch-point conditions $y^a = -b/(a-b)$ and $x = -1/\bigl((a-b)y^b\bigr)$ yield an explicit formula for the arguments of all critical rays:
\begin{equation}
\label{eq:crit_angles}
\arg(x_c) = \frac{(a-b)(\pi + 2\pi j)}{a}, \quad j \in \mathbb{Z},
\end{equation}
which, for odd $a$, directly ties the parity of the critical multiples of $\pi/a$ to the parity of the difference $(a - b)$. For $a = 7, b = 3$ ($a - b = 4$ is even), the critical rays are even multiples of $\pi/7$: along $\arg(x) = -4\pi/7$ (Fig.~\ref{fig:7_3_minus_4pi7}, where branches $n=-2$ and $n=-1$ touch at $\mathrm{Im}(w) = -3\pi/7$) and $\arg(x) = 0$ (Fig.~\ref{fig:7_3_zero}, branch point at $y = -(3/4)^{1/7}$ on the negative real semi-axis), the green markers pinpoint the crossing of the branch point at $\vert{}x\vert{} = R \approx 0.283$. Along all five non-critical rays ($\arg(x) = -\pi, -3\pi/7, -\pi/7, \pi/2, 3\pi/4$), the ultralogarithm branches remain strictly smooth and isolated within their respective phase strips of width $\frac{2\pi}{7}$.
\end{itemize}

\subsection{Fractional Exponents and Computational Stability: a = 3, b = 1/3}
\label{subsec:a3_b13}

When passing to fractional (rational) exponents in the ultraradical equation $y^a = 1 + a x y^b$, direct numerical root refinement via Newton's and Schröder's iterative methods encounters a fundamental topological obstacle: the multivaluedness of fractional powers induces a system of branch cuts in the complex $y$-plane. In the immediate vicinity of a branch point (on the boundary of the disk of convergence $\vert{}x\vert{} = R$), where the derivative vanishes and the convergence rate of the basic power series slows down, neither increasing the number of retained series terms to several thousand nor raising the working precision guarantees that iterations remain on the prescribed sheet of the Riemann surface. Upon crossing a branch cut, the algorithm abruptly jumps onto a neighboring branch.

\begin{figure}[htbp]
    \centering
    \begin{minipage}[b]{0.32\textwidth}
        \centering
        \includegraphics[width=\textwidth]{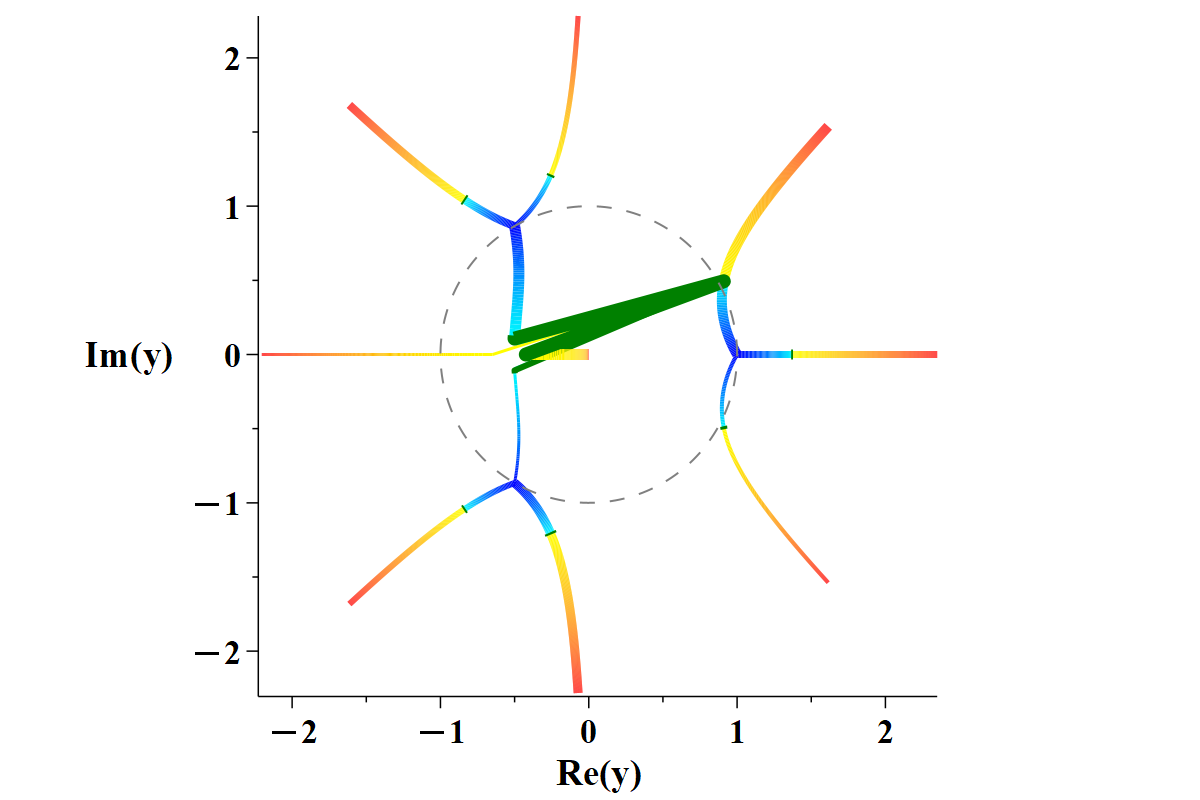}
        \textbf{(a)}
    \end{minipage}
    \hfill
    \begin{minipage}[b]{0.32\textwidth}
        \centering
        \includegraphics[width=\textwidth]{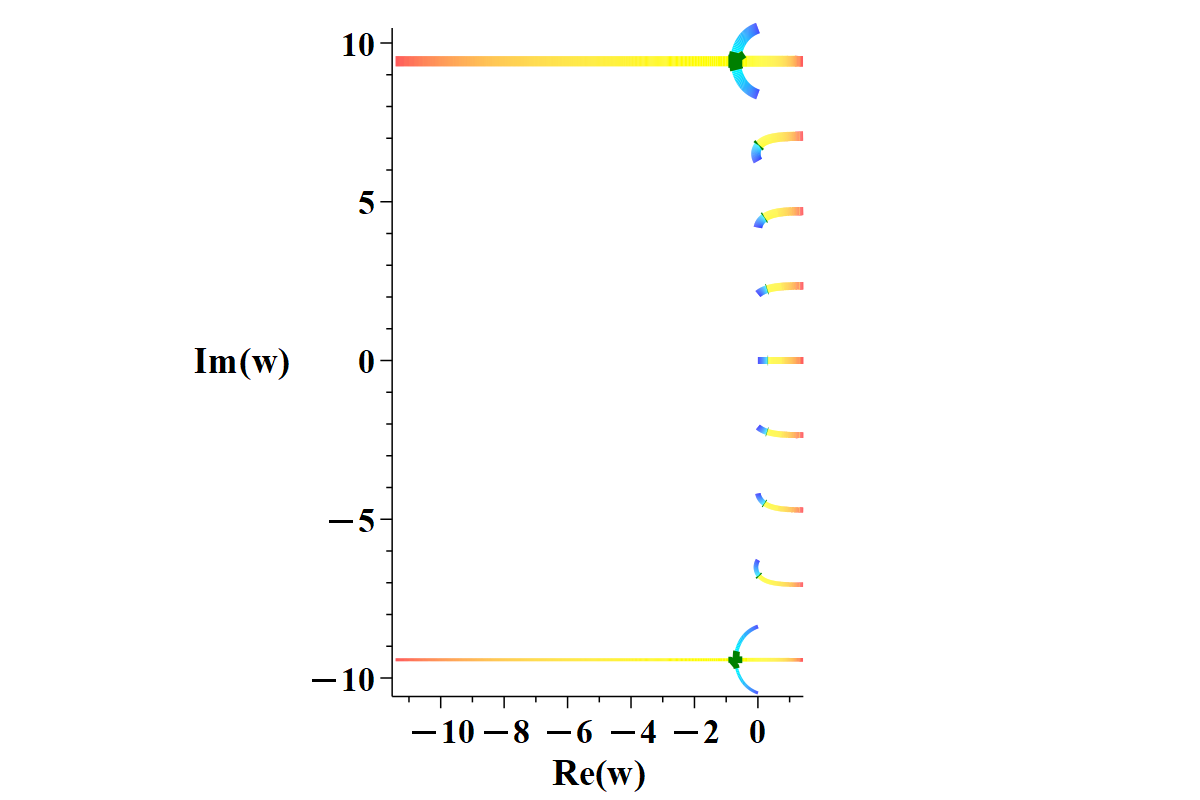}
        \textbf{(b)}
    \end{minipage}
    \hfill
    \begin{minipage}[b]{0.32\textwidth}
        \centering
        \includegraphics[width=\textwidth]{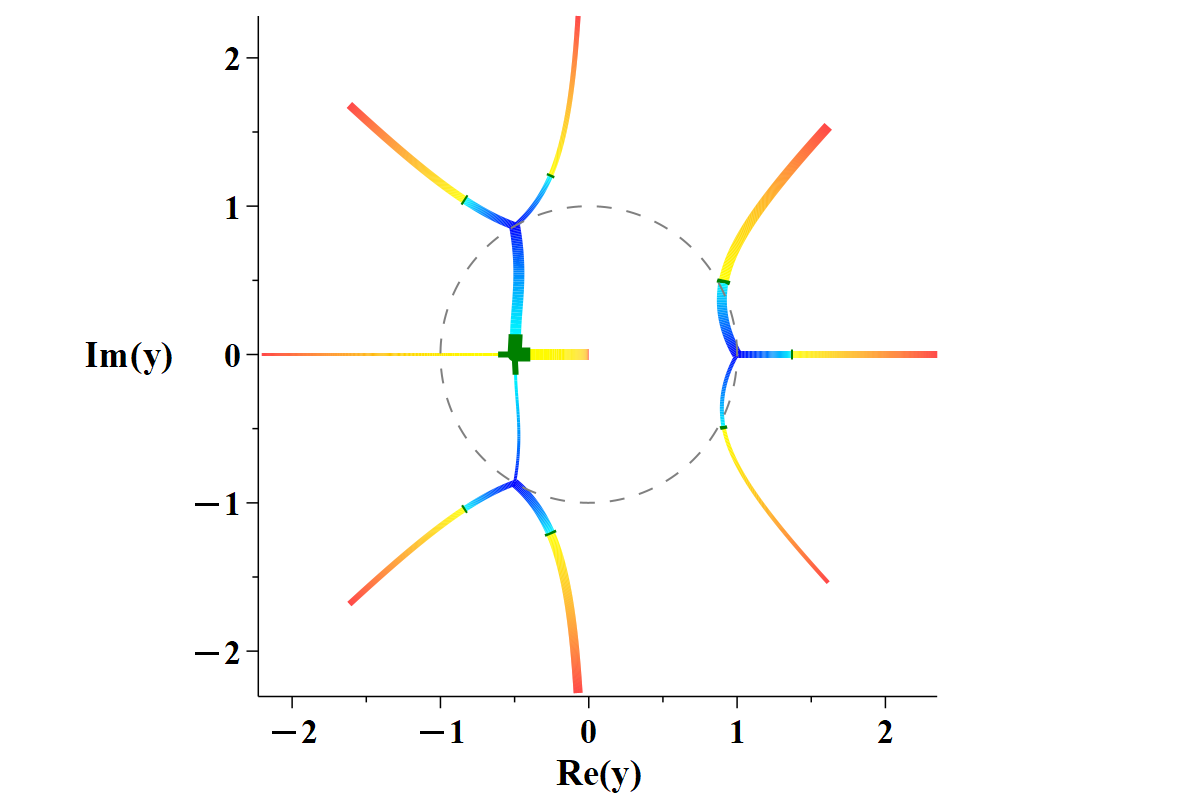}
        \textbf{(c)}
    \end{minipage}
    \caption{Demonstration of computational stability for fractional exponents $a = 3$, $b = 1/3$ along the critical ray $\arg(x) = 0$ ($x_{\max} = 3$ for $y$, $x_{\max} = 15$ for $w$): 
    \textbf{(a)}~direct ultraradical computation ($m = 1$, $n \in \{-4, \dots, 4\}$) exhibiting branch jumping near the branch point; 
    \textbf{(b)}~stratified ultralogarithm trajectories ($m = 0$, $n \in \{-5, \dots, 5\}$) free of branch cuts; 
    \textbf{(c)}~ultraradical recovered by exponentiating the refined ultralogarithm ($m = 2$, $y_n = e^{w_n}$, $n \in \{-4, \dots, 4\}$).}
    \label{fig:3_13_comparison}
\end{figure}

Consider the illustrative example with reciprocal exponents $a = 3$ and $b = 1/3$:
\begin{equation}
\label{eq:a3_b13}
e^{3w} = 1 + 3x e^{w/3}, \qquad y^3 = 1 + 3x y^{1/3},
\end{equation}
for which the radius of convergence of the basic series~\eqref{eq:radius_R} is $R = \frac{8^{1/9}}{8/3} = \frac{3 \cdot 2^{1/3}}{8} \approx 0.47247$, and the critical rays include the positive real semi-axis $\arg(x) = 0$. Substituting $y = t^3$ reduces the algebraic part of the problem to a degree-$9$ equation $t^9 - 3xt - 1 = 0$, which possesses $9$ ultraradical branches ($n \in \{-4, \dots, 4\}$). However, the starting points $y_n(0) = e^{2\pi i n / 3}$ form only $3$ base nodes on the unit circle: due to the fractional exponent $b = 1/3$, \textbf{three distinct branches} (corresponding to indices $n$ congruent modulo $3$) emerge from each base node as $\vert{}x\vert{}$ grows, whereas on the ultralogarithm hodograph (Fig.~\ref{fig:3_13_comparison}b) all of these branches are unfolded into separate horizontal strips.

Figure~\ref{fig:3_13_comparison} compares three phase hodographs along the critical ray $\arg(x) = 0$:
\begin{itemize}
\item \textbf{Direct ultraradical refinement ($m = 1$, Fig.~\ref{fig:3_13_comparison}a).} When applying the hybrid Newton--Schröder algorithm directly to $y^3 - 3x y^{1/3} - 1 = 0$ with explicit sheet-index tracking $u$, iterations lose stability near the branch point $y = -(1/8)^{1/3} = -1/2$ ($\vert{}x\vert{} = R$) due to the proximity of a branch cut: the trajectory undergoes an unphysical jump to another branch, appearing on the hodograph as a prominent green ``scar'' cutting across the interior of the unit disk.
\item \textbf{Uniformization in the ultralogarithm plane ($m = 0$, Fig.~\ref{fig:3_13_comparison}b).} In the logarithmic variable $w = \ln y$, Eq.~\eqref{eq:a3_b13} involves solely exponential functions and possesses neither finite-order branch points nor branch cuts. All branches $w_n$ remain strictly stratified within horizontal phase strips of width $\frac{2\pi}{3}$; on the boundary $\vert{}x\vert{} = R$, adjacent trajectories merely touch smoothly at the isolated branch point $\mathrm{Re}(w) = -\ln 2 \approx -0.693$, $\mathrm{Im}(w) = \pm 3\pi$, allowing Schröder iterations to continue each branch effortlessly along its own corridor.
\item \textbf{Computing the ultraradical via the ultralogarithm ($m = 2$, Fig.~\ref{fig:3_13_comparison}c).} Performing root refinement in the stable $w$-space ($m = 0$) and returning to the $y$-plane only at the final step via $y_n = e^{w_n}$ completely eliminates branch jumping. The resulting hodograph exhibits a clean green cross at the branch point $y = -0.5$ without any numerical artifacts and without requiring an auxiliary sheet-selection parameter $u$.
\end{itemize}

Thus, the ultralogarithm acts as a natural \emph{uniformizing map} for multivalued algebraic and transcendental equations with fractional exponents: shifting iterative refinement into the $w$-plane transforms a topologically delicate branch-cut navigation problem into a smooth computational procedure.

\subsection{Transcendental Incommensurable Exponents: a = \texorpdfstring{$\pi$}{pi}, b = e}
\label{subsec:api_be}

The final stage of the computational experiment addresses the transcendental case with incommensurable irrational exponents $a = \pi$ and $b = e$:
\begin{equation}
\label{eq:api_be}
e^{\pi w} = 1 + \pi x e^{ew},
\end{equation}
which corresponds to the ultraradical $y^{\pi} = 1 + \pi x y^{e}$. Here, the radius of convergence of the basic series~\eqref{eq:radius_R} is $R = \frac{\vert{}1 - \pi/e\vert{}^{e/\pi}}{\pi - e} \approx 0.47281$. Unlike integer-exponent equations, an incommensurable ratio $b/a = e/\pi$ destroys periodicity in the index $n$ with period $2\pi i$: the notion of a finite set of ``independent'' roots and the algebraic Vieta's law~\eqref{eq:log_vieta} no longer apply, and the entire countably infinite family of branches $w_n$ ($n \in \mathbb{Z}$) becomes genuinely independent.

\begin{figure}[htbp]
    \centering
    \begin{minipage}[b]{0.48\textwidth}
        \centering
        \includegraphics[width=\textwidth]{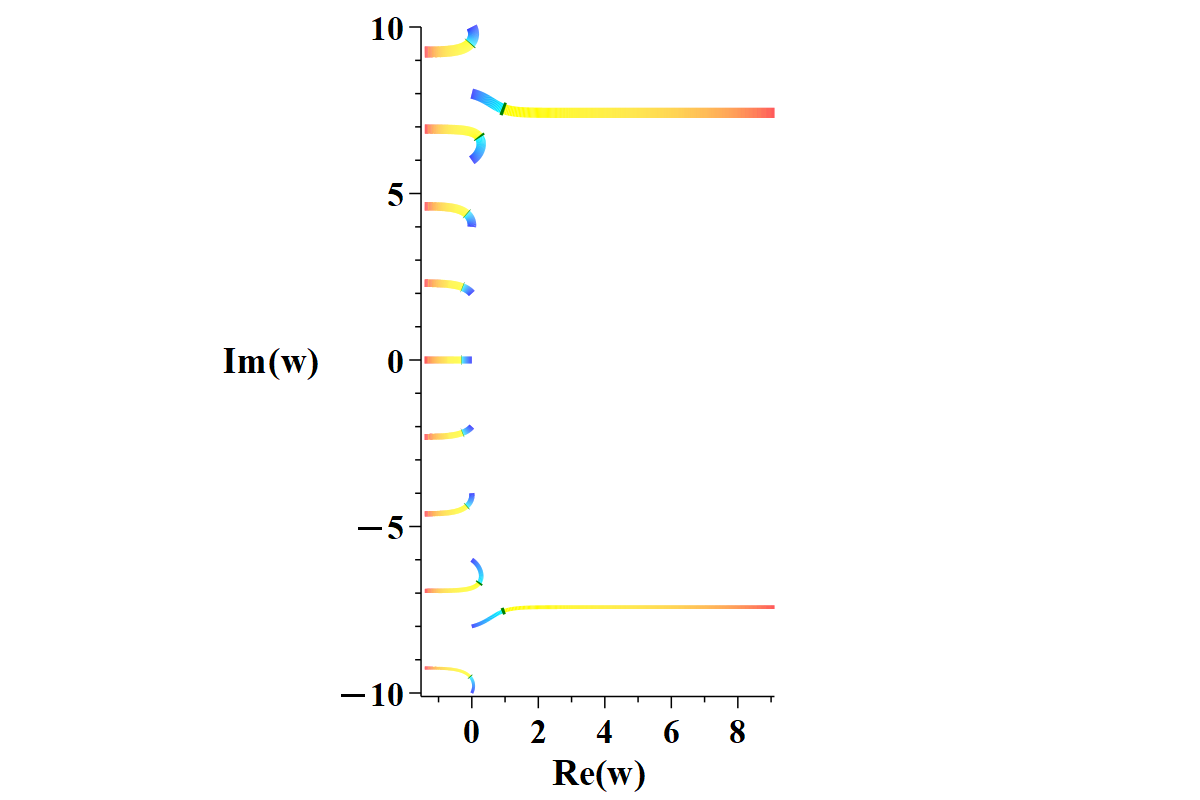}
        \textbf{(a)}
    \end{minipage}
    \hfill
    \begin{minipage}[b]{0.48\textwidth}
        \centering
        \includegraphics[width=\textwidth]{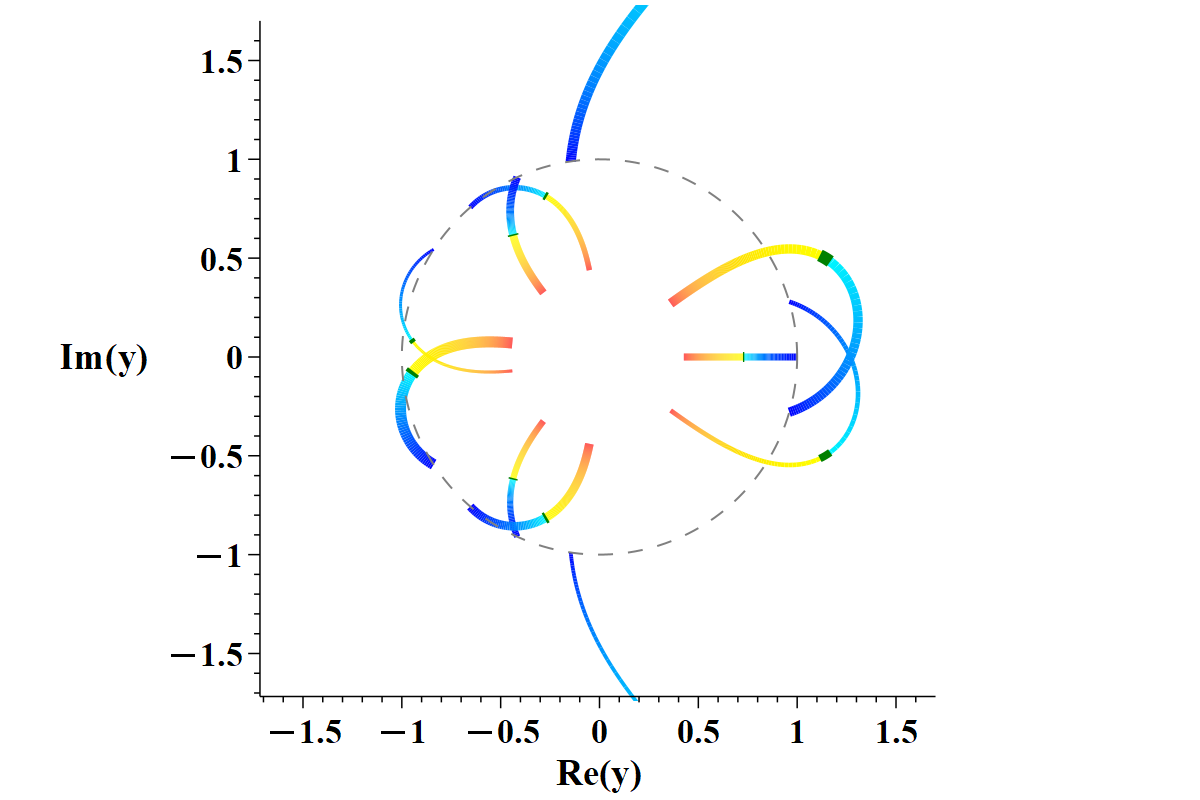}
        \textbf{(b)}
    \end{minipage}
    \caption{Phase hodographs for $a = \pi$, $b = e$ at angle $\arg(x) = -\pi$: 
    \textbf{(a)}~ultralogarithm ($m=0$, equation $e^{\pi w} = 1 + \pi x e^{ew}$), branches $n \in \{-5, \dots, 5\}$, $x_{\max} = 15$; 
    \textbf{(b)}~ultraradical ($m=1$, equation $y^{\pi} = 1 + \pi x y^{e}$), branches $n \in \{-5, \dots, 5\}$, $x_{\max} = 3$.}
    \label{fig:pi_e_minus_pi}
\end{figure}

\begin{figure}[htbp]
    \centering
    \begin{minipage}[b]{0.48\textwidth}
        \centering
        \includegraphics[width=\textwidth]{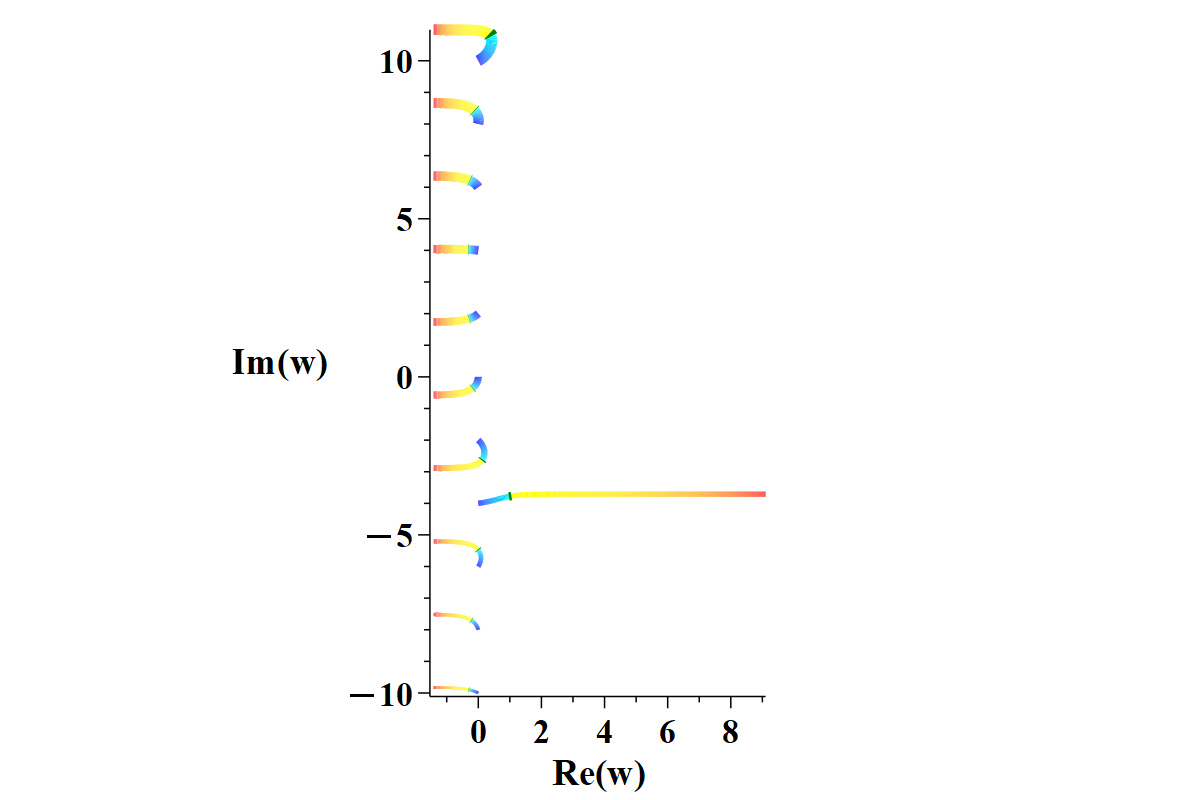}
        \textbf{(a)}
    \end{minipage}
    \hfill
    \begin{minipage}[b]{0.48\textwidth}
        \centering
        \includegraphics[width=\textwidth]{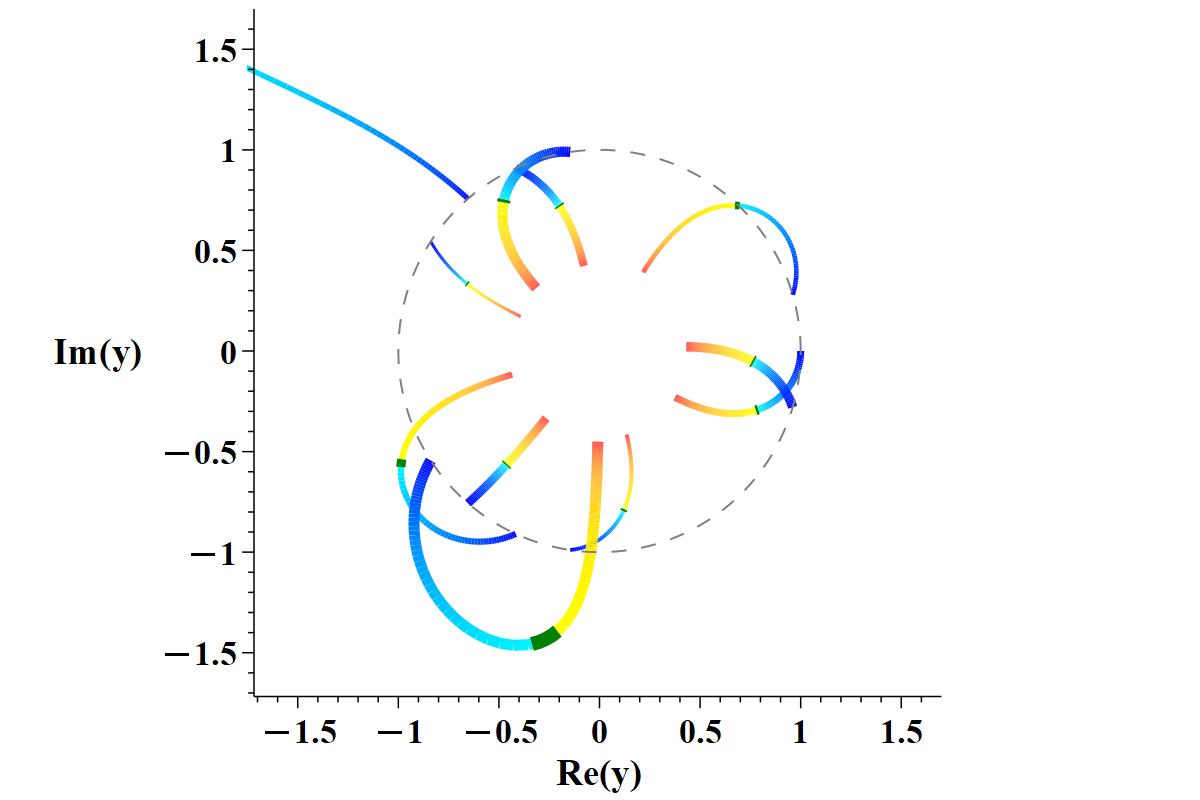}
        \textbf{(b)}
    \end{minipage}
    \caption{Phase hodographs for $a = \pi$, $b = e$ at angle $\arg(x) = -\pi/2$: 
    \textbf{(a)}~ultralogarithm ($m=0$, equation $e^{\pi w} = 1 + \pi x e^{ew}$), branches $n \in \{-5, \dots, 5\}$, $x_{\max} = 15$; 
    \textbf{(b)}~ultraradical ($m=1$, equation $y^{\pi} = 1 + \pi x y^{e}$), branches $n \in \{-5, \dots, 5\}$, $x_{\max} = 3$.}
    \label{fig:pi_e_minus_pi2}
\end{figure}

\begin{figure}[htbp]
    \centering
    \begin{minipage}[b]{0.48\textwidth}
        \centering
        \includegraphics[width=\textwidth]{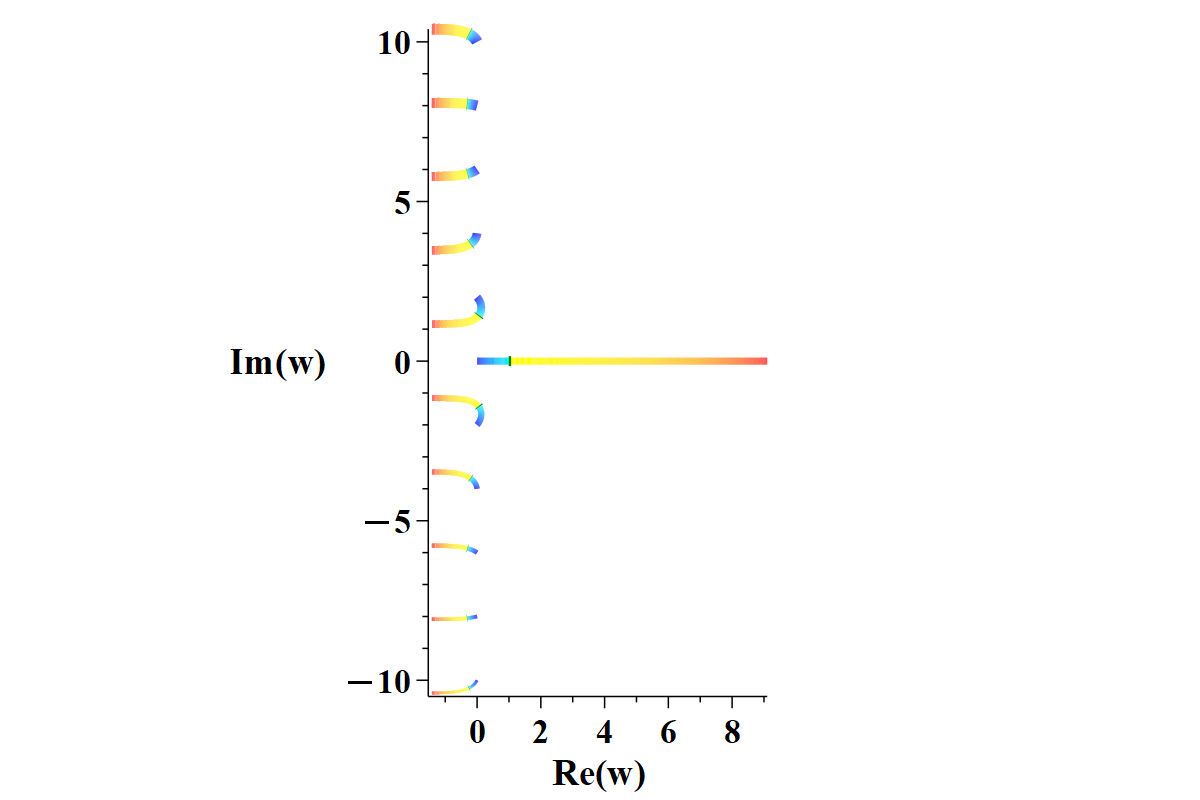}
        \textbf{(a)}
    \end{minipage}
    \hfill
    \begin{minipage}[b]{0.48\textwidth}
        \centering
        \includegraphics[width=\textwidth]{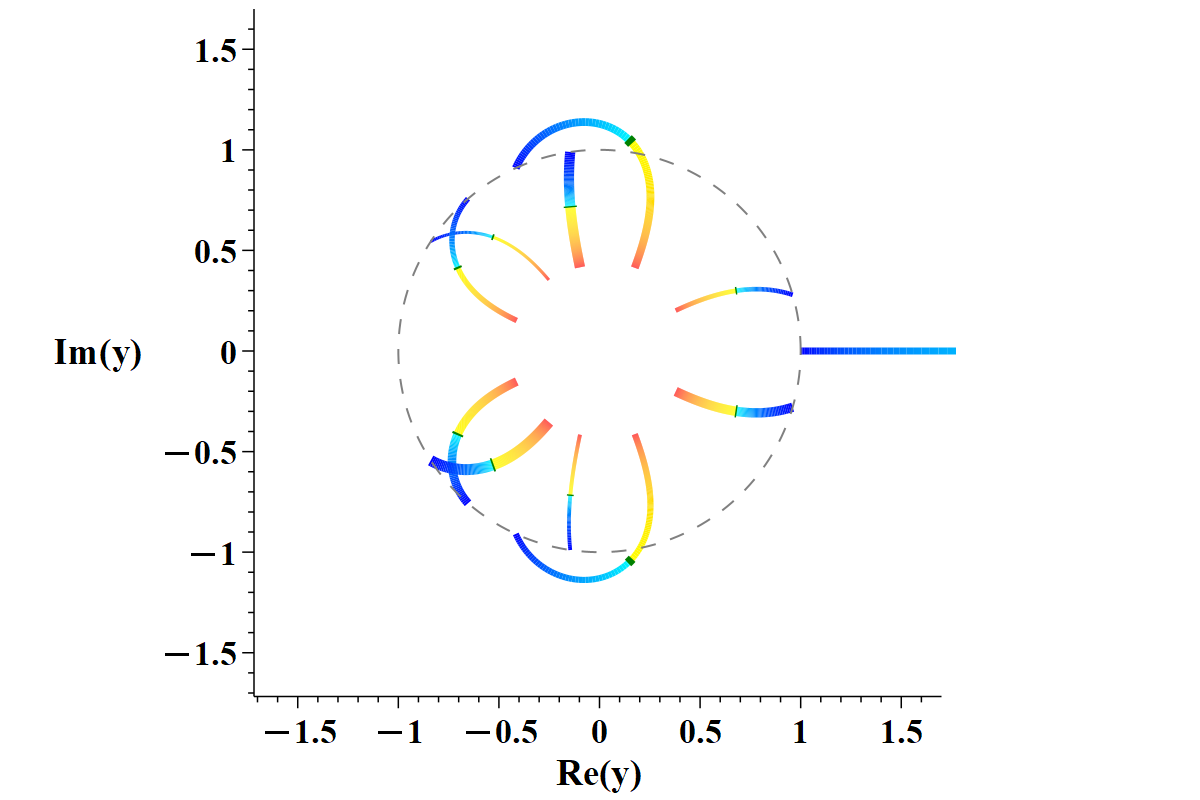}
        \textbf{(b)}
    \end{minipage}
    \caption{Phase hodographs for $a = \pi$, $b = e$ at angle $\arg(x) = 0$: 
    \textbf{(a)}~ultralogarithm ($m=0$, equation $e^{\pi w} = 1 + \pi x e^{ew}$), branches $n \in \{-5, \dots, 5\}$, $x_{\max} = 15$; 
    \textbf{(b)}~ultraradical ($m=1$, equation $y^{\pi} = 1 + \pi x y^{e}$), branches $n \in \{-5, \dots, 5\}$, $x_{\max} = 3$.}
    \label{fig:pi_e_zero}
\end{figure}

\begin{figure}[htbp]
    \centering
    \begin{minipage}[b]{0.48\textwidth}
        \centering
        \includegraphics[width=\textwidth]{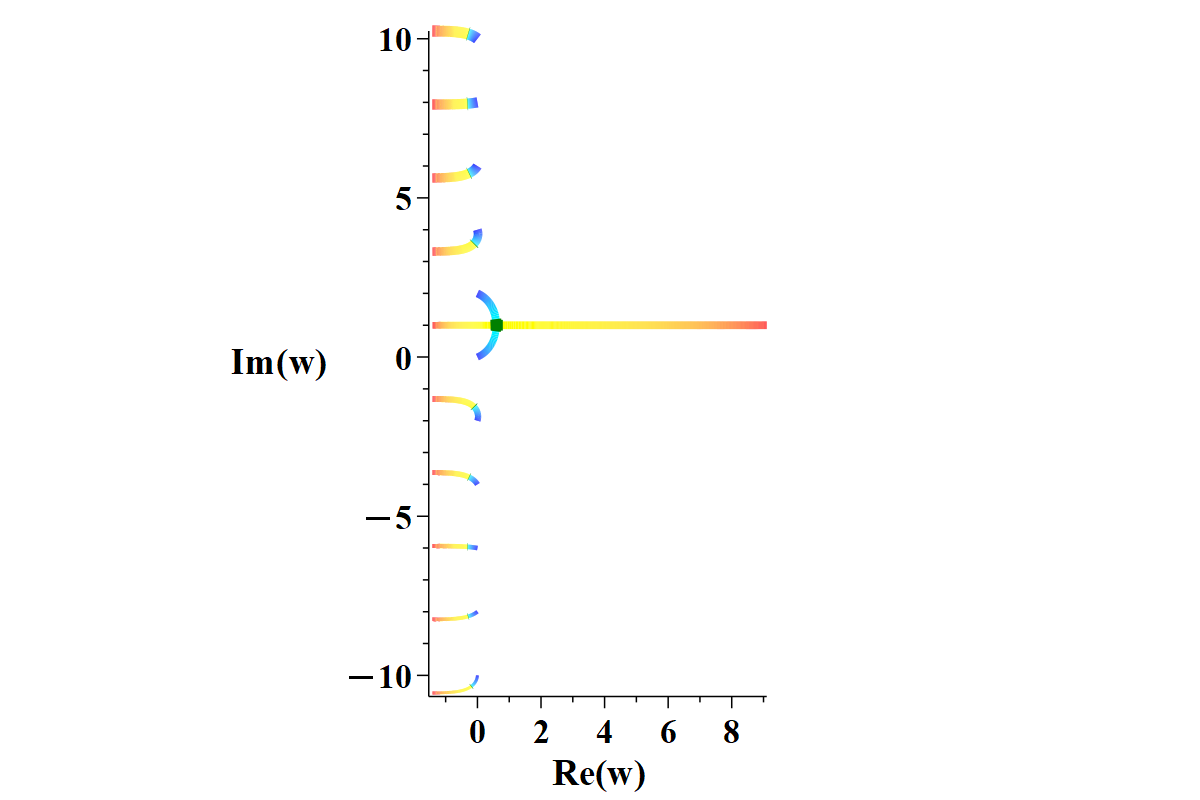}
        \textbf{(a)}
    \end{minipage}
    \hfill
    \begin{minipage}[b]{0.48\textwidth}
        \centering
        \includegraphics[width=\textwidth]{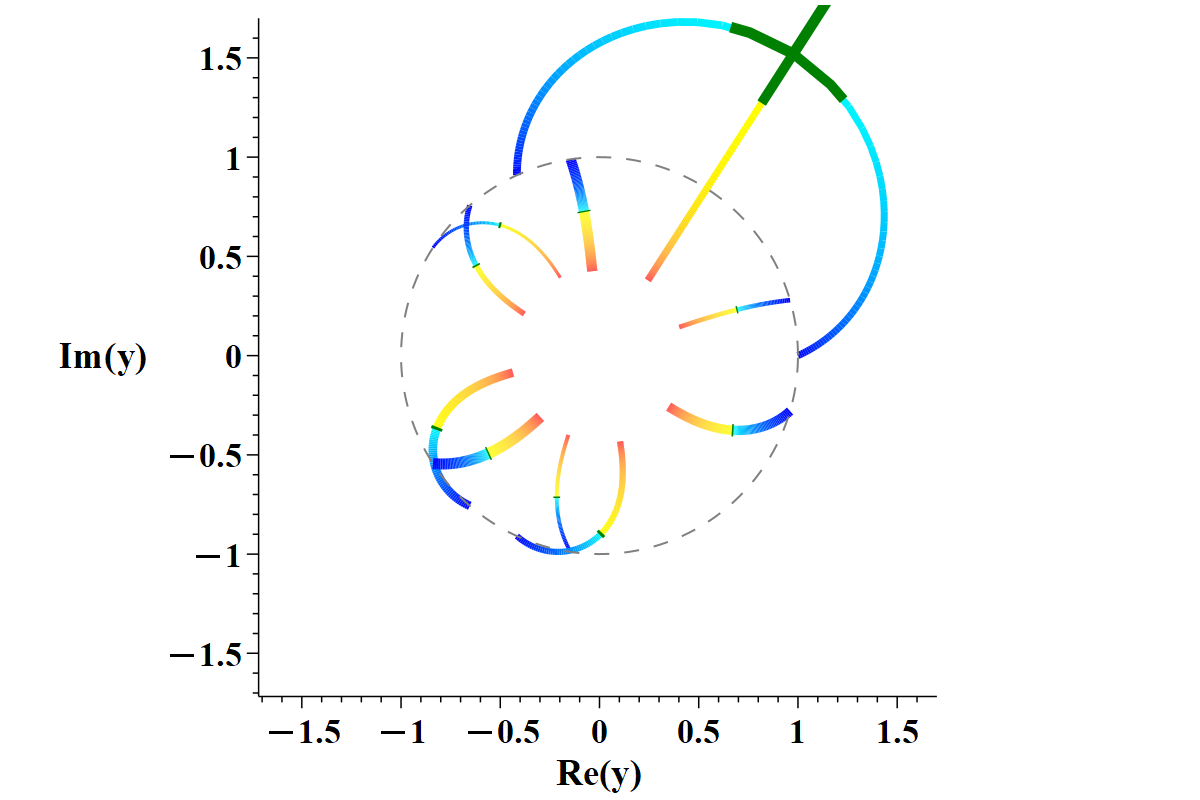}
        \textbf{(b)}
    \end{minipage}
    \caption{Phase hodographs for $a = \pi$, $b = e$ at critical angle $\arg(x) = \pi - e$: 
    \textbf{(a)}~ultralogarithm ($m=0$, equation $e^{\pi w} = 1 + \pi x e^{ew}$), branches $n \in \{-5, \dots, 5\}$, $x_{\max} = 15$; 
    \textbf{(b)}~ultraradical ($m=1$, equation $y^{\pi} = 1 + \pi x y^{e}$), branches $n \in \{-5, \dots, 5\}$, $x_{\max} = 3$.}
    \label{fig:pi_e_crit}
\end{figure}

\begin{figure}[htbp]
    \centering
    \begin{minipage}[b]{0.48\textwidth}
        \centering
        \includegraphics[width=\textwidth]{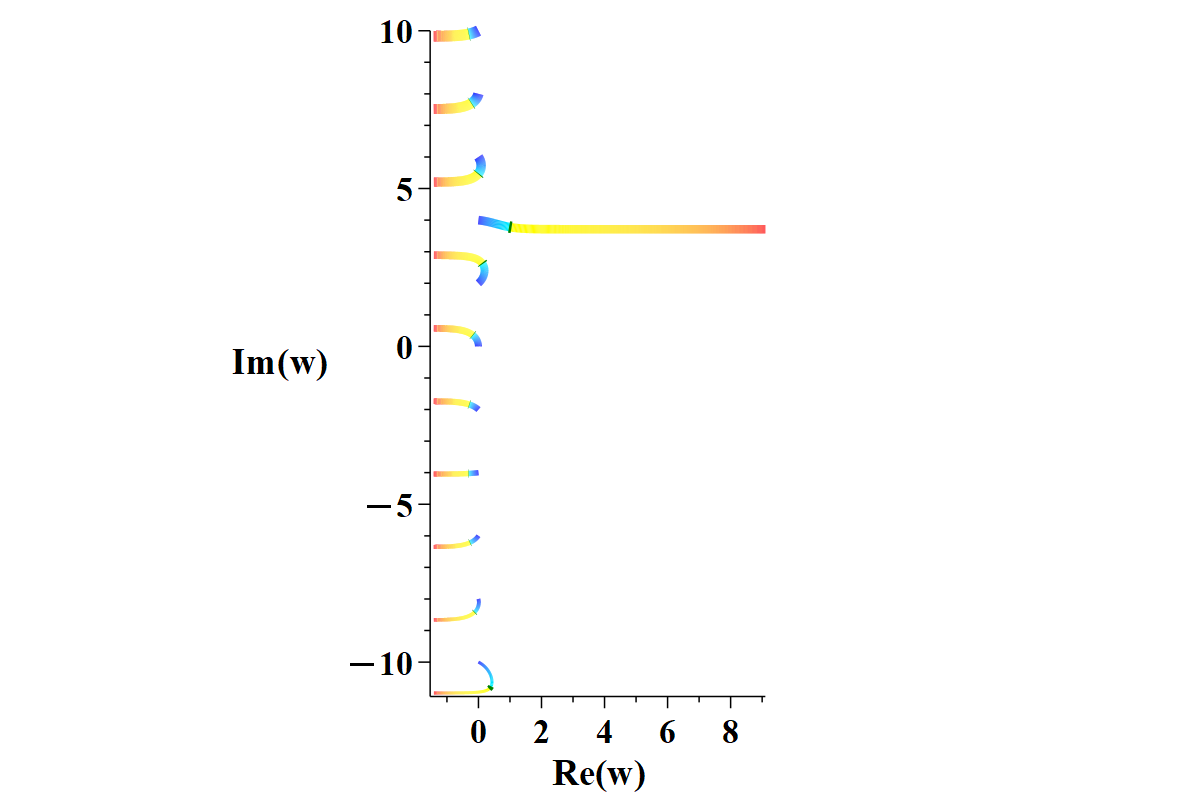}
        \textbf{(a)}
    \end{minipage}
    \hfill
    \begin{minipage}[b]{0.48\textwidth}
        \centering
        \includegraphics[width=\textwidth]{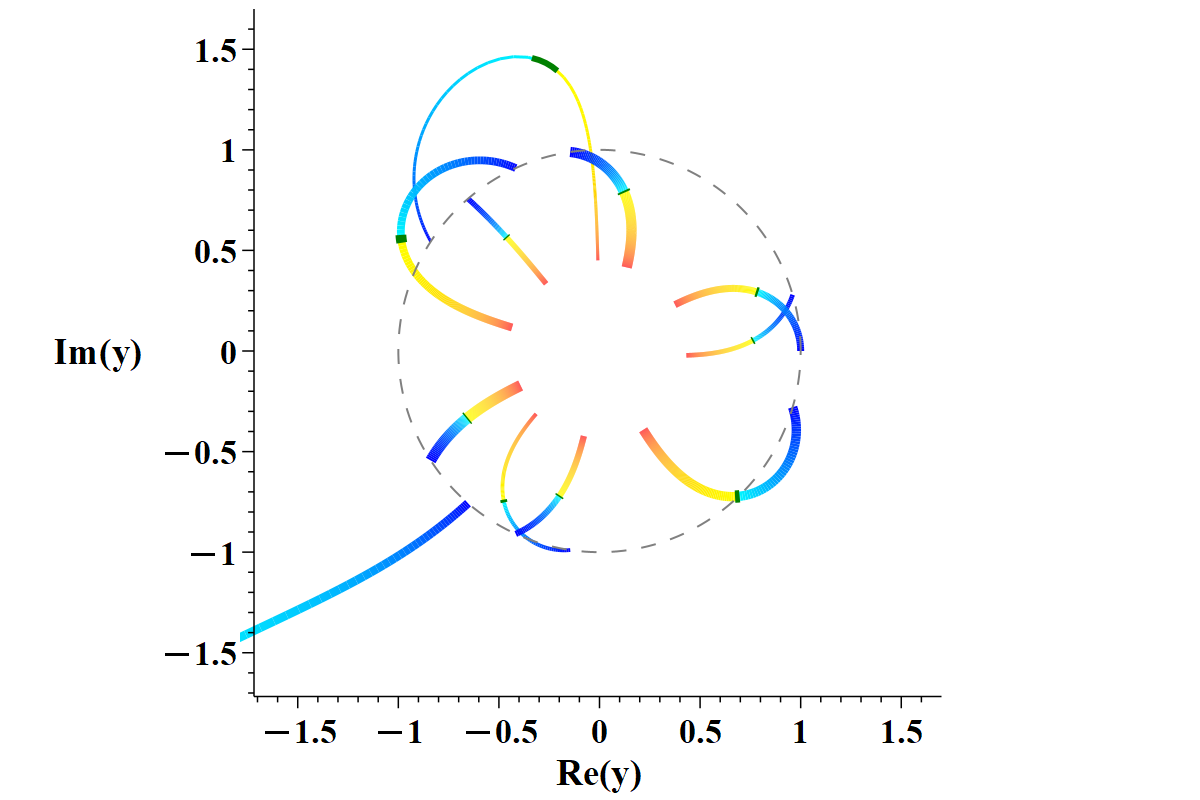}
        \textbf{(b)}
    \end{minipage}
    \caption{Phase hodographs for $a = \pi$, $b = e$ at angle $\arg(x) = \pi/2$: 
    \textbf{(a)}~ultralogarithm ($m=0$, equation $e^{\pi w} = 1 + \pi x e^{ew}$), branches $n \in \{-5, \dots, 5\}$, $x_{\max} = 15$; 
    \textbf{(b)}~ultraradical ($m=1$, equation $y^{\pi} = 1 + \pi x y^{e}$), branches $n \in \{-5, \dots, 5\}$, $x_{\max} = 3$.}
    \label{fig:pi_e_plus_pi2}
\end{figure}

Figures~\ref{fig:pi_e_minus_pi}--\ref{fig:pi_e_plus_pi2} present the paired phase hodographs of eleven branches ($n \in \{-5, \dots, 5\}$) of the ultralogarithm ($x_{\max} = 15$) and the ultraradical ($x_{\max} = 3$) along five rays $\arg(x) \in \{-\pi, -\pi/2, 0, \pi - e, \pi/2\}$. Inspecting these portraits reveals four fundamental features of the transcendental regime:
\begin{itemize}
\item \textbf{Aperiodic petal fan instead of Newton polygons.} Because the base points $y_n(0) = e^{2i n}$ lie on the unit circle with an angular step of $2$~rad ($\approx 114.59^\circ$) incommensurable with $2\pi$, the closed polygon of starting points breaks down. On the ultraradical phase hodographs (panels \textbf{b}), the inner branches converging to zero as $y \sim (\pi x)^{-1/e}$ (family~\eqref{eq:ulog_outer_k}) form an open quasi-periodic ``petal fan'' with an angular spacing of $\frac{2\pi}{e} \approx 132.4^\circ$. For real ray angles $\arg(x) = -\pi$ and $\arg(x) = 0$ (Figs.~\ref{fig:pi_e_minus_pi}b and~\ref{fig:pi_e_zero}b), this fan exhibits exact mirror symmetry about the real axis, whereas for complex angles (Figs.~\ref{fig:pi_e_minus_pi2}b, \ref{fig:pi_e_crit}b, and~\ref{fig:pi_e_plus_pi2}b) the symmetry is broken.
\item \textbf{Small-denominator effect $(a - b) = \pi - e \approx 0.4233$.} Due to the small difference between the exponents, the frequency ratio is $\frac{a - b}{a} = \frac{\pi - e}{\pi} \approx 0.1347 \approx \frac{1}{7.42}$. Consequently, on the ultralogarithm hodographs (panels \textbf{a}) within the seven central branches, the vast majority of trajectories recede to the left ($\mathrm{Re}(w) \to -\infty$, series~\eqref{eq:ulog_outer_k} with denominator $b = e$), reaching only $\mathrm{Re}(w) \approx -1.4$ at $\vert{}x\vert{} = 15$. Only a single rare branch (occurring once every $7.42$ branches on average) falls into the family~\eqref{eq:ulog_outer_h} with denominator $a - b = \pi - e$: its real part $\mathrm{Re}(w) \sim \frac{\ln(\pi\vert{}x\vert{})}{\pi - e}$ grows $\frac{e}{\pi - e} \approx 6.42$ times faster, stretching into a long horizontal ray up to $\mathrm{Re}(w) \approx +9.1$. As the ray angle $\arg(x)$ rotates, this single ``fast'' branch shifts successively upward in index $n$ (from $n = -2$ in Fig.~\ref{fig:pi_e_minus_pi2}a through $n = 0$ in Fig.~\ref{fig:pi_e_zero}a and $n = +1$ in Fig.~\ref{fig:pi_e_crit}a to $n = +2$ in Fig.~\ref{fig:pi_e_plus_pi2}a).
\item \textbf{Branch points for transcendental exponents.} By the general critical-angle formula~\eqref{eq:crit_angles} at $j = 0$, the value $\arg(x_c) = \pi - e \approx 0.4233$ corresponds to a critical ray. Figure~\ref{fig:pi_e_crit} visually confirms that branch touching occurs even in equations with irrational powers: the adjacent ultralogarithm branches $n = 0$ and $n = 1$ (Fig.~\ref{fig:pi_e_crit}a) touch on their shared sector boundary $\mathrm{Im}(w) = 1.0$ at $\mathrm{Re}(w) = \frac{1}{\pi}\ln(e/\pi) \approx -0.046$ (green marker), while in the ultraradical plane (Fig.~\ref{fig:pi_e_crit}b) this corresponds to an outer trajectory crossing at $\vert{}x\vert{} = R \approx 0.473$.
\item \textbf{Robustness of the sector-based criterion.} Across all five experiments, the sector algorithm accurately identifies the conjugate series type and continuation index even under strong frequency asymmetry and incommensurable parameters.
\end{itemize}

\begin{remark}[On projection self-intersections and branch identities for irrational exponents]
On the petal-fan hodographs of the ultraradical (Figs.~\ref{fig:pi_e_minus_pi}b--\ref{fig:pi_e_plus_pi2}b), numerous trajectory crossings appear in the complex $y$-plane without green branch-point markers, where the intersecting curves have \emph{different colors} (for instance, a blue segment $\vert{}x\vert{} < R$ of one branch crosses a yellow-to-red segment $\vert{}x\vert{} > R$ of another). Unlike true branch points, where two branches coalesce at the \emph{exact same} value of $x$, these apparent crossings indicate that distinct branches $n_1 \neq n_2$ attain the same complex value $y_0$ at \emph{different} moduli $\vert{}x_1\vert{} \neq \vert{}x_2\vert{}$ along the same ray $\arg(x_1) = \arg(x_2)$. While impossible for integer-exponent equations (where $x = (y^a - 1)/(ay^b)$ is uniquely determined by $y$), the multivaluedness of irrational powers $y^\pi$ and $y^e$ across different Riemann sheets generates an infinite family of non-trivial identities of the form $U_{n_1}(a;b;x_1) = U_{n_2}(a;b;x_2)$ with $x_1 \neq x_2$.
\end{remark}

\subsection{Connection to the Lambert W Function and Physical Branches of the Sagdeev Method}
\label{subsec:sagdeev_lambert}

To conclude the computational experiments, we demonstrate the practical significance of rigorous branch numbering using a classical problem in plasma physics---computing the stationary profile of an ion-acoustic soliton via the Sagdeev pseudopotential method~\cite{Dubinov2006, Dubinov2020}.
\begin{equation}
\label{eq:sagdeev_energy}
\mathcal{K}\left(1 - \frac{1}{\eta^2}\right) + \frac{\Theta}{\gamma - 1}\left(1 - \eta^{\gamma - 1}\right) - \Phi = 0,
\end{equation}
where $\eta = n_2 / n_1 > 0$ is the dimensionless particle density (ratio of local to unperturbed density), $\mathcal{K} = mV^2/2$ is the initial kinetic energy, $\Theta = \gamma k_B T_1$ is the thermal energy scale, $\Phi = Ze\varphi$ is the electrostatic potential energy, and $\gamma$ is the adiabatic index.

For an arbitrary adiabatic index $\gamma > 1$, multiplying Eq.~\eqref{eq:sagdeev_energy} by $\eta^2(\gamma-1)/\Theta$ reduces it to the canonical trinomial form $y^a = 1 + a x y^b$ with exponents $a = \gamma + 1$ and $b = 2$, where $\eta = v_1 y$, $v_1 = \bigl(-\mathcal{K}(\gamma-1)/\Theta\bigr)^{1/(\gamma+1)}$, and the argument $x$ depends linearly on the potential $\Phi$. In the isothermal limit $\gamma \to 1$, the thermal term approaches $-\Theta\ln\eta$, and the closed-form solution is expressed via the Lambert $W$ function:
\begin{equation}
\label{eq:sagdeev_lambert}
\eta(\varphi) = \exp\!\left(\frac{\mathcal{K} - Ze\varphi}{\Theta} + \frac{1}{2} W_n\!\left(-\frac{2\mathcal{K}}{\Theta}\exp\!\left(-\frac{2(\mathcal{K} - Ze\varphi)}{\Theta}\right)\right)\right).
\end{equation}

\begin{figure}[htbp]
    \centering
    \begin{minipage}[b]{0.32\textwidth}
        \centering
        \includegraphics[width=\textwidth]{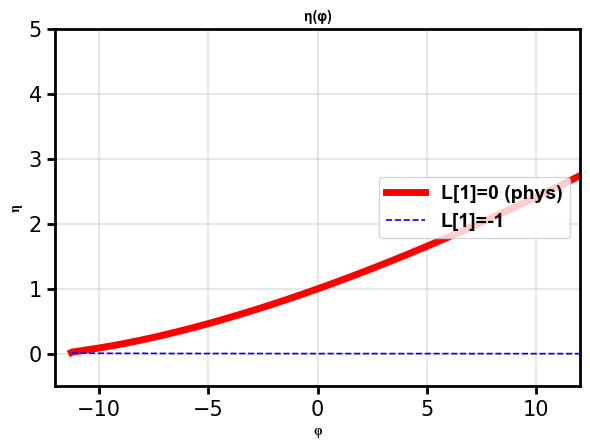}
        \textbf{(a)}
    \end{minipage}
    \hfill
    \begin{minipage}[b]{0.32\textwidth}
        \centering
        \includegraphics[width=\textwidth]{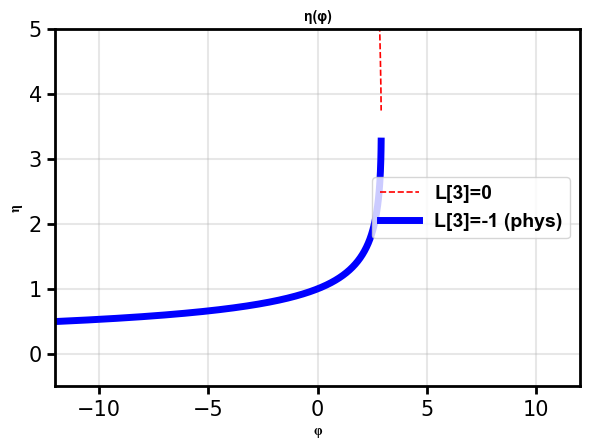}
        \textbf{(b)}
    \end{minipage}
    \hfill
    \begin{minipage}[b]{0.32\textwidth}
        \centering
        \includegraphics[width=\textwidth]{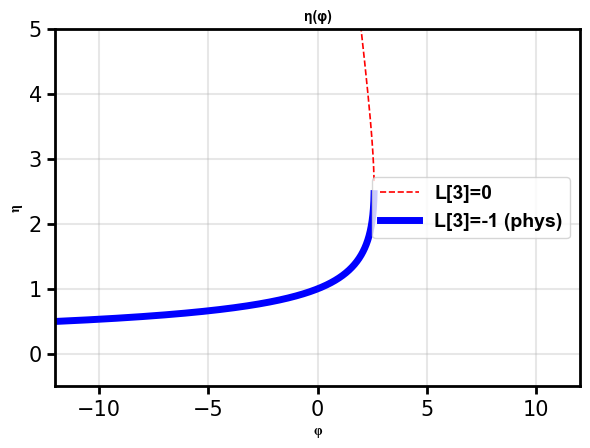}
        \textbf{(c)}
    \end{minipage}
    \caption{Real branches of the dimensionless density $\eta(\varphi)$ in the Sagdeev method (solid curve: physical branch passing through $\eta(0) = 1$; dashed curve: non-physical branch): 
    \textbf{(a)}~subsonic electron branch ($n = 0$, $\gamma = 5/3$); 
    \textbf{(b)}~supersonic ion branch in the isothermal limit ($\gamma = 1$, Lambert $W_{-1}$ solution); 
    \textbf{(c)}~supersonic ion branch at $\gamma = 5/3$ (ultraradical solution preserving the branch index $n = -1$).}
    \label{fig:sagdeev_branches}
\end{figure}

Unlike abstract complex-plane analysis, only \emph{real positive} branches $\eta(\varphi)$ passing through the unperturbed state $\eta(0) = 1$ are physically admissible here. The fundamental distinction between plasma components is governed by the ratio of the wave velocity $V$ to the species' own acoustic speed $c_s = \sqrt{\gamma k_B T_1 / m}$:
\begin{itemize}
\item \textbf{Electron branch ($n = 0$, subsonic regime, Fig.~\ref{fig:sagdeev_branches}a).} Because of their small mass, electrons have a thermal velocity far exceeding the soliton phase speed ($V < c_{s,e}$). In this subsonic regime, the physical condition $\eta(0) = 1$ is satisfied on the \textbf{principal branch $n = 0$} (solid red curve), which continuously approaches the principal branch $W_0$ of the Lambert $W$ function as $\gamma \to 1$.
\item \textbf{Ion branch in the isothermal limit ($n = -1 \leftrightarrow W_{-1}$, $\gamma = 1$, Fig.~\ref{fig:sagdeev_branches}b).} Heavy and relatively cold ions move faster than their own acoustic speed relative to the wave ($V > c_{s,i}$). In this supersonic regime, the $n = 0$ branch (dashed red curve) misses the point $(0, 1)$, and the physically realized solution switches to the \textbf{branch with index $n = -1$} (solid blue curve), which at $\gamma = 1$ coincides identically with the lower real branch $W_{-1}$ of the Lambert $W$ function.
\item \textbf{Deformation of the ion branch as $\gamma$ increases ($n = -1$, $\gamma = 5/3$, Fig.~\ref{fig:sagdeev_branches}c).} As the adiabatic index increases to $\gamma = 5/3$, the Lambert $W$ representation no longer applies and is superseded by the ultraradical (or ultralogarithm $w = \ln\eta$). The sector-based criterion automatically preserves the branch index ($n = -1$), while the turning point where branches $n = 0$ and $n = -1$ coalesce corresponds precisely to the master-series argument reaching its radius of convergence ($\vert{}x\vert{} = R$), where the local ion flow velocity matches the acoustic speed.
\end{itemize}

The topological and physical meaning of this correspondence is particularly noteworthy. Inside the disk of convergence of the inner expansion~\eqref{eq:ulog_inner} ($\vert{}x\vert{} < R$), there is only a single principal zero stationary anchor with zero index ($f^{0}_{0}$, $n = 0$). By contrast, the physical plasma process under consideration takes place exclusively \emph{outside} the disk of convergence of the inner series ($\vert{}x\vert{} > R$), where the outer master series of types $+$~\eqref{eq:ulog_outer_h} and $-$~\eqref{eq:ulog_outer_k} converge simultaneously, yielding \emph{two} principal moving anchors with zero indices in total: the positive anchor $f^{+}_{0}$ ($n_+ = 0$) and the negative anchor $f^{-}_{0}$ ($n_- = 0$). In general, which moving anchor ($f^{+}_{n}$ or $f^{-}_{n}$) continues the $n = 0$ branch depends on the argument of the variable $x$ in the master series. In this physical model, however, the structure of the energy balance equation~\eqref{eq:sagdeev_energy} rigidly fixes the correspondence: it depends neither on the sign of the particle charge $Z$ (whether electrons or positrons, positive or negative ions) nor on the adiabatic exponent $\gamma$ (integer, fractional, or irrational). Instead, the physical branch index is determined solely by the particle flow velocity relative to its own acoustic speed $c_s = \sqrt{\gamma k_B T_1 / m}$: all subsonic components ($V < c_s$, light hot particles) are governed by the principal branch $n = 0$, supported in the outer region by the positive moving anchor $f^{+}_{0}$ ($n_+ = 0$), whereas all supersonic components ($V > c_s$, heavy cold particles) universally lie on the branch $n = -1$, supported by the negative moving anchor $f^{-}_{0}$ ($n_- = 0$). In practice, this eliminates the need to track the degree $a = \gamma + 1$ of the equation or reduce indices modulo $a$: the supersonic ion branch retains the invariant universal index $n = -1$ across all parameter regimes.

Finally, at the point where the radius of convergence is reached ($\vert{}x\vert{} = R$), mathematics and physics coincide completely: the convergence boundary corresponds to the acceleration (for the subsonic branch $n = 0$) or deceleration (for the supersonic branch $n = -1$) of the particle stream in the wave frame up to the local acoustic speed ($M = 1$). At this sonic barrier, the smooth stationary soliton profile breaks down as the solution branches enter the complex plane. Nevertheless, when guided by rigid channel walls or an external force field, the physical flow itself does not stop at the sonic barrier~--- it simply leaves the domain of validity of the reversible trinomial model~\eqref{eq:sagdeev_energy} and transitions into a supercritical choked or shock-wave regime (a detailed analysis of supercritical flow beyond the master-series radius of convergence and the critical-point isomorphism theorem are presented in~\cite{Berezin_Zenodo_Thermo}).

\section{Conclusion}
\label{sec:conclusion}

In this work, we have solved the problem of global analytic continuation and unified one-dimensional branch indexing for the ultralogarithm $\mathrm{ULog}_n(a;b;x)$ with arbitrary real parameters $a > b > 0$. The main results of the study are summarized as follows:

\begin{enumerate}
\item \textbf{Sector-based conjugation criterion.} We proved a theorem establishing a one-to-one correspondence between the branch index $n_{0} \in \mathbb{Z}$ inside the disk of convergence ($\vert{}x\vert{} < R$) and the indices $n_{+}, n_{-} \in \mathbb{Z}$ of the conjugate master series in the outer region ($\vert{}x\vert{} > R$) by partitioning the imaginary axis into phase sectors of width $\frac{2\pi}{a}$. For integer exponents, a logarithmic analog of Vieta's theorem was derived, classifying the solutions into $a$ independent and a countably infinite family of dependent branches.
\item \textbf{Geometric stratification of trajectories.} Using multi-channel phase hodographs across integer ($a=2, b=1$; $a=5, b=1$; $a=7, b=3$), fractional ($a=3, b=1/3$), and transcendental incommensurable ($a=\pi, b=e$) exponents, we demonstrated that ultralogarithm trajectories remain stratified within isolated phase strips of width $\frac{2\pi}{a}$ without self-intersections, in contrast to the wrapped trajectories of the ultraradical.
\item \textbf{Uniformization and computational stability.} Passing to the logarithmic space $w = \ln y$ eliminates branch cuts and finite-order branch points of fractional powers. Iterative root refinement via Newton's and Schröder's methods in the $w$-plane followed by exponentiation $y_n = e^{w_n}$ ($m = 2$ mode) completely suppresses branch jumping near the convergence boundary without requiring an auxiliary sheet-selection parameter $u$.
\item \textbf{Physical invariance of branch indexing.} Using the Sagdeev pseudopotential method in plasma physics, we verified the continuous preservation of branch numbering in the Lambert $W$ limit ($\gamma \to 1$): subsonic components are universally governed by the principal branch $n = 0$ ($W_0$), whereas supersonic components lie on the branch $n = -1$ ($W_{-1}$), regardless of particle charge sign or adiabatic index $\gamma$.
\end{enumerate}

\section{Continuation of Research}
This research is ongoing. Updated versions and supplementary materials are available in the open Zenodo repository. To ensure you always access the latest version, please use the concept DOIs:
\begin{itemize}
\item Russian version: \url{https://doi.org/10.5281/zenodo.24218691}
\item English version: \url{https://doi.org/10.5281/zenodo.24219280}
\end{itemize}
The author would be grateful for any comments, questions, or suggestions for further development of the formalism.

\section*{Acknowledgements}

The author expresses sincere gratitude to:

\begin{itemize}
\item \textbf{Prof. A. E. Dubinov} (RFNC--VNIIEF, SarFTI NRNU MEPhI) for fruitful discussions and for pointing out the potential applications of the ultra-radical in plasma physics problems.

\item \textbf{Pavel V. Berezin, Alexey V. Berezin, and Andrey V. Gruzdov} for their extensive work on computations, numerical verification, and fruitful discussions.

\item \textbf{Prof. Tuyen Truong} (University of Oslo, Norway) and \textbf{Prof. Abdulaziz Alhaidari} (Saudi Center for Theoretical Physics) for their valuable feedback, critical questions, and constructive discussions that greatly improved this work.
\end{itemize}

\paragraph*{AI Disclosure.}
During the preparation of this manuscript, the author used AI language models (DeepSeek, Qwen, Gemini) for editing, structuring, and consistency checking. In particular, the key Theorem~1 on the sector-based conjugation criterion was formulated and rigorously proved during a collaborative analytical dialogue between the author and the Gemini language model (Google), completing an investigation initiated by the author a year earlier when the first empirical branch-matching patterns were identified with the assistance of DeepSeek. All mathematical derivations, numerical results, and conclusions have been reviewed and verified by the author, who takes full responsibility for the final version of the manuscript.

\paragraph*{Computational Tools.}
Analytical computations and numerical verification were performed using the computer algebra systems Wolfram Alpha and Maple, the Python programming language with scientific libraries, and the author's own Master-J laboratory in JavaScript with native support for arbitrary-precision complex numbers (DOI: \href{https://doi.org/10.5281/zenodo.15633781}{10.5281/zenodo.15633781}). In addition, starting from version~3.1 of the software repository~\cite{Berezin2026Py}, the \texttt{SuperMaster} module includes the universal Maple procedure \texttt{ULR := proc(m, n, a, b, x, M := 500)}, which evaluates the ultralogarithm at $m = 0$, the ultraradical directly at $m = 1$, and the uniformized ultraradical via exponentiation ($y_n = e^{w_n}$ after refining $w_n$ at $m = 0$) when $m = 2$.
Here, \texttt{M} denotes the forced series truncation limit: if the target precision is not reached within \texttt{M} terms, Newton or Schröder iterations are automatically triggered, and the procedure returns not only the computed root, but also the number of series terms actually used and the Riemann sheet index~$u$.

\section*{Data Availability Statement}
The data that support the findings of this study are openly available in the following repositories:

\begin{itemize}
\item A repository with useful implementations of the proposed algorithms in various programming languages, which will be gradually updated (Concept DOI, always resolves to the latest version): 
\href{https://doi.org/10.5281/zenodo.20355433}{10.5281/zenodo.20355433}.

\item The Master-J laboratory (JavaScript implementation with arbitrary-precision complex numbers, including the merge operation):
\href{https://doi.org/10.5281/zenodo.15633781}{10.5281/zenodo.15633781}.

\item Preprint versions of the manuscript:
\href{https://arxiv.org/abs/2512.18643}{arXiv:2512.18643} (detailed version) and
\href{https://arxiv.org/abs/2609.33459}{arXiv:2609.33459} (consolidated version).
\end{itemize}

\section*{Funding}
This research received no external funding.

\section*{Conflict of Interest}
The author declares no conflict of interest.

\end{document}